%% file: manuscript.tex
\documentclass[a4paper]{article}
\usepackage{ffmath}
\usepackage{ffmathad}
\usepackage{geometry}
\usepackage{amsmath,amssymb,amsfonts}
\usepackage{mathtools}
\usepackage{commath}
\usepackage{subfig}
\usepackage{authblk}
\usepackage{todonotes}
\usepackage{standalone}
\usepackage{multicol}
\usepackage{bm}
\usepackage{doi}
\usepackage[square,numbers]{natbib}
\hypersetup{pdfborder = 0 0 0, colorlinks=false,
 linkcolor=black,citecolor=black, filecolor=black, urlcolor=black}

\usepackage{algorithm,algpseudocode}
\usepackage[title]{appendix}
\usepackage{threeparttable}
\usepackage{booktabs}% Set rules for tables

\usepackage{pgfplots}
\pgfplotsset{compat=1.13}
\usetikzlibrary{calc}

\usepackage{cleveref} % Include cleveref last!

\newtheorem{remark}{Remark}
\usetikzlibrary{patterns} % LATEX and plain TEX when using TikZ

\title{An integral equation method for the advection-diffusion equation on time-dependent domains in the plane}

\author[1]{Fredrik Fryklund\thanks{Corresponding author. E-mail address: \href{ffr@kth.se}{ffry@kth.se}}}
\author[1]{Sara P\aa lsson}
\author[1]{Anna-Karin Tornberg}

\affil[1]{Department of Mathematics, KTH Royal Institute of Technology, Stockholm, Sweden}

\date{}

\begin{document}

\maketitle

\begin{abstract}
\label{s:abstract}
\input{abstract.tex}
\end{abstract}

% \begin{keyword}
% Heat equation, boundary integral method, modified Helmholtz \revMine{equation}, Yukawa potential, quadrature, complex domains, function extension, Rothe's method
% \end{keyword}


\section{Introduction}
\label{s:introduction}
\input{introduction.tex}
\section{Formulation}
\label{s:formulation}
\input{formulation.tex}
\section{Numerical methods}
\label{s:numericalmethods}
\input{numericalmethods.tex}
\section{Partition of unity extension}
\label{s:pux}
\input{pux.tex}
\section{Numerical experiments}
\label{s:examples}
\input{numericalexamples.tex}
\section{Conclusions}
\label{s:conclusions}
\input{conclusions.tex}
\section*{Acknowledgments}

The authors gratefully acknowledge the support from the Swedish Research Council under Grant No. $2015$-$04998$, and the G\"{o}ran Gustafsson Foundation for Research in Nature and Medicine. The authors thank Lukas Bystricky for discussions and for an implementation to solve Stokes equation in a periodic channel.

\appendix
\section{Supplementary material for spectral Ewald summation}
\label{s:appsupewald}
\input{appsupewald.tex}

\bibliographystyle{abbrvnat_mod}
\bibliography{references}

\end{document}

%% file: abstract.tex
% ABSTRACT
%!TeX spellcheck = en-US
Boundary integral methods are attractive for solving homogeneous linear constant coefficient elliptic partial differential equations on complex geometries, since they can offer accurate solutions with a computational cost that is linear or close to linear in the number of discretization points on the boundary of the domain. However, these numerical methods are not straightforward to apply to time-dependent equations, which often arise in science and engineering. We address this problem with an integral equation-based solver for the advection-diffusion equation on moving and deforming geometries in two space dimensions. In this method, an adaptive high-order accurate time-stepping scheme based on semi-implicit spectral deferred correction is applied. One time-step then involves solving a sequence of non-homogeneous modified Helmholtz equations, a method known as elliptic marching. Our solution methodology utilizes several recently developed methods, including special purpose quadrature, a function extension technique and a spectral Ewald method for the modified Helmholtz kernel. Special care is also taken to handle the time-dependent geometries. The numerical method is tested through several numerical examples to demonstrate robustness, flexibility and accuracy.
%Boundary integral methods are attractive for solving homogeneous linear constant coefficient elliptic partial differential equations on complex geometry: the ill-conditioning associated with discretizing the governing equations is avoided, boundary data is simple to incorporate, high accuracy can be attained throughout the domain, and for time-dependent geometries only the boundaries need to be updated. However, these numerical methods are not straightforward to apply to time-dependent equations, which often arise in science and engineering. We address this problem with an integral equation-based solver for the advection-diffusion equation in complicated moving and deformable geometries in two space dimensions. Our method can be explained as applying a high-order accurate time-stepping scheme based on semi-implicit spectral deferred correction. Then, one time-step involves solving a sequence of non-homogeneous modified Helmholtz equations, a method known as elliptic marching, which requires several recently developed methods. These include special purpose quadrature, function extension techniques and spectral Ewald methods. Special care is also taken to handle the time-dependent geometries. The numerical method is tested through several numerical examples to demonstrate robustness, flexibility and accuracy.
%%% Local Variables:
%%% mode: latex
%%% TeX-master: "manuscript.tex"
%%% End:

%% file: introduction.tex
% INTRODUCTION
%!TEX root = ../manuscript.tex
Many physical phenomena of interest in science and engineering involve moving domains. Yet there is no standardized numerical solution method for such problems, thus solving partial differential equations (PDEs) on time-dependent geometries is an active area of research. The main problem for classical numerical methods is that the underlying grid, where the data is represented, must be remeshed to follow the movement of the domain. Clearly, this becomes expensive as the temporal and the spatial resolutions are increased. Instead, it is desirable to keep the underlying grid fixed, which comes with its own set of challenges, depending on which method is used. We will begin by introducing the problem and our proposed method, and with this as background we will then discuss other methods.

We propose a boundary integral-based method for solving the isotropic advection-diffusion equation, also called the convection-diffusion equation, on time-dependent geometries in two spatial dimensions. The advection-diffusion equation is formulated as
\begin{align}
    \od{\adu\fp{t,\x}}{t} + \velf\fp{t,\x}\cdot \nabla \adu\fp{t,\x} &= \diffconst\Delta \adu\fp{t,\x} + \force\fp{t,\x},\quad \x\in \domain\fp{t},\quad \tO<t,\label{eq:adu}\\
    \pd{\adu\fp{t,\x}}{\normalscal} &= \adbc\fp{t,\x},\quad  \x\in \bdry\fp{t},\quad \tO<t,\label{eq:adbc}\\
    \adu\fp{\tO,\x} &= \adic\fp{\x},\quad \x\in \domain\fp{\tO},\label{eq:adic}\\
    \od{\x\fp{t}}{t} &= \velf\fp{t,\x},\quad \x\in\bdry\fp{t},\quad \tO<t, \label{eq:adbdry}
\end{align}
where $\adu$ is an unknown quantity to be solved for, $\velf$ is a velocity field that is either given explicitly or is the solution to an additional PDE, the constant $\diffconst > 0$ is the diffusion coefficient, the source term $\force\fp{t,\x}$ is a smooth explicitly given function, $\bdry(t)$ is the boundary of the domain $\domain(t)\subset\mathbb{R}^{2}$ at time $t$ with at least two continuous derivatives, $\normal$ is the outward directed unit normal at $\x\in\bdry\fp{t}$ as shown in figure \ref{fig:Omega}, $\adbc$ is given Neumann boundary data, and the domain $\bar{\domain}\fp{t}=\domain\fp{t}\cup\bdry\fp{t}$ is bounded and continuous in $t$. We assume that the velocity field $\velf$ transports $U$ and $\bdry$ with the same velocity. Multiply connected domains are not addressed in this paper, but including them requires no additional theory or computational methods. The extension to exterior problems is straight forward, assuming suitable far-field conditions are supplied. We do however consider the advection-diffusion equation above with periodic boundary conditions for parts of the boundary.

The advection-diffusion equation is often part of models of complex problems in computational physics, for example the dynamics of multiphase flows on the microscale with soluble surfactants (surface-active agents). Surfactants decrease the surface tension at interfaces between fluids. At the microscale these forces dominate, as the surface-to-volume ratio is large, thus surfactants have a significant impact on the dynamics of the multiphase flow. They can be used to control droplet dynamics, which is central for many microfluidic applications \cite{TehShia-Yen2008Dm}. Surfactants are usually present on the interface and in one of the fluid phases. This can be modeled with the advection-diffusion equation, where boundary interface conditions describe the exchange between the drop interface and the bulk. Considerable effort is given to accurately simulate microfluidic flows, as computer simulations is a much used tool to study such problems.

Here boundary integral methods enter the picture. The velocity field in drop and particle suspension problems at the microscale is often modeled with the Stokes equations, which can be recast as a second-kind Fredholm integral equation \cite{pozrikidis_1992}. Consequently, the dimensionality is reduced as the unknowns are confined to the boundary of the domain, including particle and drop surfaces, the latter being fluid-fluid interfaces. Thus only the boundary needs to be discretized \cite{PALSSON2019218,SORGENTONE2018167}. Tracking the evolution of the interfaces and the interface conditions are naturally handled with a boundary integral method, effectively avoiding remeshing the domain as it evolves in time. Another advantage of a boundary integral formulation is that high accuracy can be maintained as droplets, vesicles or particles get arbitrarily close to each other.

Time-dependent PDEs can be approximated by a sequence of inhomogeneous elliptic PDEs under appropriate temporal discretization, such as implicit-explicit (IMEX) spectral deferred correction (SDC) or IMEX Runge-Kutta methods (R-K) \cite{Dutt2000,AscherIMEX}. The resulting PDEs can be solved with a boundary integral method, which has been addressed for stationary geometry \cite{FryklundFredrik2020Aien,AFKLINTEBERG2020109353}. In the case of the advection-diffusion equation the diffusion term is treated implicitly, reducing the equation to a sequence of modified Helmholtz equations, 
\begin{align}
  \label{eq:mh}
  \mhaa\mhu\fp{\x}-\Delta\mhu\fp{\x} &= \mhrhs\fp{\x},\quad \x\in\domain,\quad \alpha\in\mathbb{R}_{+},\\
  \label{eq:mhbc}
  \pd{\mhu\fp{\x}}{\normalscal} &= \mhbc\fp{\x}, \quad\x\in\bdry,
\end{align}
where $\mhu$ corresponds to an approximation of $\adu$ at the next instance in time, $\mhaa = \ordo{\tfrac{1}{\dt}}$ with time-step $\dt$, and the right-hand side $\mhrhs$ is updated at every time-step. We write the solution to \eqref{eq:mh}, \eqref{eq:mhbc} as a sum of a single layer potential and a volume potential, the latter is due to the PDE being inhomogeneous for $f\ne 0$.

Computing the required volume potential poses several numerical challenges if $\domain$ is not geometrically simple, e.g. rectangular. We take the approach of creating a continuation $\mhrhse$ of $\mhrhs$, with a specified regularity and compact support on a box $\suppbox$ embedding $\bar{\domain}$, such that $\mhrhse\vert_{\domain} = \mhrhs$. The extension $f^{e}$ is constructed with partition of unity extension (PUX) \cite{FRYKLUNDPUX}. Now, let $\mhpu$ be the convolution of $\mhrhse$ with the modified Helmholtz Green's function over $B$. The convolution can be computed efficiently with geometrically unaware fast summation methods, such as a fast Fourier transform (FFT). The solution to \eqref{eq:mh}, \eqref{eq:mhbc} is the sum $\mhpu+\mhhu$, if $\mhhu$ satisfies the homogeneous counterpart of \eqref{eq:mh} with Neumann boundary data $\mhbc\fp{\x}-\partial\mhpu\fp{\x}/\partial\normalscal$. The solution $\mhhu$ is efficiently solved for with a boundary integral method, obtained via a Fredholm integral equation of the second kind, aided with special quadrature methods to maintain high accuracy for target points arbitrarily close to the boundary \cite{klinteberg2019adaptive}. The decomposition of $\mhu$ into $\mhpu$ and $\mhhu$ allows us to enjoy the benefits of boundary integral methods for inhomogeneous elliptic PDEs.

When discretizing the boundary integral equation, and using an iterative method to solve the linear system, matrix-vector multiplications with a dense matrix must be rapidly computed. For the modified Helmholtz equation these matrix-vector multiplications are discrete sums with modified Bessel functions of the second kind of zeroth- or first-order. In a free-space setting, the sums can be computed directly at cost $\mathcal{O}(N^{2})$ for $N$ sources and targets. As the problem size grows, the evaluation becomes prohibitively expensive. There are several ways to speed up the computations, such as with a fast multipole method (FMM) \cite{greengard1987fast,Kropinski2011modHelm} or a spectral Ewald method. The latter is an FFT based method, which is especially suitable for periodic problems. It has been implemented for the Laplace equation in $3$D \cite{AfKlintebergLudvig2014FEsf,LindboDag2011Saif}, the Stokes equation in $3$D for both free-space and periodic setting \cite{LindboDag2010Safs,afKlintebergLudvig2017FEsf}, and in $2$D for a periodic setting \cite{PalssonSara2020Aiem}. In this paper we present an Ewald method for the modified Helmholtz equation. An Ewald method is based on an Ewald decomposition, meaning that the sum is split into a ''real space`` sum which converges rapidly in $\mathbb{R}^{2}$, and a ''$\fovar$ space`` sum which converges rapidly in the frequency domain. By accelerating the evaluation of the latter using the FFT, the total computational cost of the method is $\mathcal{O}(N\log N)$.% These two sums are computed in $\mathcal{O}(N\log N)$ time. There are Ewald decompositions both for the periodic and free-space case, with only minor differences between the two Ewald methods \cite{afKlintebergLudvig2017FEsf}.

%At the core of solving the homogeneous modified Helmholtz equation is evaluating discrete sums with modified Bessel functions of the second kind of zeroth- or first-order. In a free-space setting, the sums can be computed directly at cost $\mathcal{O}(N^{2})$ for $N$ sources and targets. As the problem size grows, the evaluation becomes unfeasibly costly. There are several ways to speed up the computations, such as with a spectral Ewald method mentioned above. It is a fast Fourier transform (FFT) based method, which is especially suitable for periodic problems. It has been implemented for the Laplace equation in $3$D \cite{AfKlintebergLudvig2014FEsf,LindboDag2011Saif}, the Stokes equations in $3$D for both free-space and periodic setting \cite{LindboDag2010Safs,afKlintebergLudvig2017FEsf}, and in $2$D for a periodic setting \cite{PalssonSara2020Aiem}. The Ewald method is based on an Ewald decomposition of the sum into a ''real space`` sum which converges rapidly in $\mathbb{R}^{2}$, and a ''$\fovar$ space`` sum which converges rapidly in the frequency domain. The spectral Ewald method computes these two sums in $\mathcal{O}(N\log N)$ time, both for the periodic and free-space case, with only minor differences between the two \cite{afKlintebergLudvig2017FEsf}.

When $\domain$ evolves discretely in time the right hand side $f$ in \eqref{eq:mh} may not exist everywhere in $\domain$. The function extension technique PUX is also used to address this problem of missing data. The basic concept of PUX is to blend local extensions by a partition of unity into a global extension. The regularity of the extension is directly related to the construction of said partition of unity. The local extensions are computed through interpolation with radial basis functions (RBF). The associated interpolation matrix is precomputed such that the ill-conditioning of the problem is significantly reduced.

The main contributions in this paper are the following.
\begin{itemize}
  \item An accurate boundary integral method for the solution, and the gradient thereof, of the modified Helmholtz equation with Neumann boundary data over a wide range of values of the parameter $\mha$ in  \eqref{eq:mh}.
  \item A spectral Ewald method for the kernel associated with the modified Helmholtz equation, both in a periodic and free-space setting.
  \item A high order solver of the advection-diffusion equation on moving domains in a boundary integral framework, that is adaptive in time.
\end{itemize}
The result is a solver for the advection-diffusion equation \eqref{eq:adu}--\eqref{eq:adbdry} of order $K$ in time and tenth order in space, with $K(K-1)$ modified Helmholtz equations \eqref{eq:mh}, \eqref{eq:mhbc} to solve per time-step.

An alternative boundary integral approach to the one presented in this paper is the one of \citeauthor{WangShaobo2019FHIE} \cite{WangShaobo2019FHIE}. Here, the solution is obtained through a direct approximation of the heat kernel, thus avoiding the discretization of the differential operator with respect to time. The corresponding layer potential contains no volume potential, but  temporal integration from the initial time $\tO$ to some $t$. This integral is split into two: one ``history part'' from $\tO$ to $t-\delta$, and one ``local'' part, integrating from $t-\delta$ to $t$. The former is computed with Fourier approximation and nonuniform FFT, which tends to be inefficient for small $\delta$ due to having to overresolve the frequency domain. The latter is computed with an FMM. The result is an efficient algorithm, that potentially can be extended to inhomogenoeus problems as well. A similar approach for the inhomogeneous counterpart, but for stationary domains, is found in \cite{Liheat}, where the heat layer potential contains a volume potential. Here the layer potential is also composed of a local part and a history part

In \cite{FrachonThomas2019Acfe} Frachon and Zahedi present a Cut finite element method (CutFEM) for solving the Navier-Stokes equations on time-dependent geometry. In a CutFEM the geometry can be arbitrarily located on a fixed underlying mesh. The involved space-time integrals are in a variational formulation, and are computed by first using a quadrature rule in time. Then for each quadrature point the integrals in space are computed with CutFEM. In time, the trial and test functions are discontinuous between  one time interval $[\tn,\tnp]$ to another. Continuity in time is enforced, weakly, by adding a correction term. Stabilization terms are also added which implicitly extend data to points in the new domain where data does not exist. Despite the addition of preconditioning and stabilizers, the method suffers from large condition numbers as the order of the quadrature scheme is increased. Thus while CutFEM methods are versatile in handling time-dependent geometry, high-order in time is difficult to achieve. %Furthermore, they suffer from a severe loss of accuracy for close body interactions, such as when two drops are close to each other.
Another FEM-based approach is presented in \cite{lehrenfeld2019eulerian} where an unfitted mesh is used in the spatial domain. A time-discretization based on finite difference approximation is used, rather than space-time variational formulation. Here extension methods are applied to address the problem of missing data when the domain is updated in time. 

Yet another approach is a meshless method of high spatial order by V. Shankar et al. in \cite{ShankarVarun2020AEHM}. Here the solution is represented by an approximation with radial basis functions, centered at a scattered node set. These nodes are used when applying finite differences. It then remains to solve for the interpolation weights to obtain the solution. The corresponding differentiation matrix does not need to be recomputed at each iteration in time, instead it is updated, which with their techniques is only required for nodes neighboring the boundary. Impressive results are presented for spatial convergence, but any high order methods for temporal accuracy are not presented.

The immersed boundary smooth extension method is used in \cite{HUANG2021110162} to solve the bulk advection-diffusion equation on time-dependent geometry for simulating the Stefan problem. By applying an IMEX scheme they also obtain a sequence of modified Helmholtz equations to solve per time-step. Unlike the method presented in this paper, the function data $f$ is not extended outside $\domain$, but instead an extension of the solution that satisfies the given PDE is solved for. This is achieved by imposing the unknown solution as boundary data for a homogeneous poly-harmonic PDE of degree $k$, for which the solution is the extension of regularity $k$. Third-order temporal and pointwise spatial convergence is demonstrated for a classical Stefan problem.

This paper is organized as follows. First the mathematical formulation of the modified Helmholtz equation is presented in section \ref{s:formulation}, followed by section \ref{s:numericalmethods} covering time-stepping schemes, and numerical methods for the modified Helmholtz equation. Thereafter, PUX and its application to time-dependent geometry are introduced in section \ref{s:pux}, followed by section \ref{s:examples} where numerical results are presented. Finally, conclusions and an outlook are given section \ref{s:conclusions}.

%% file: formulation.tex
% FORMULATION
%!TEX root = ../manuscript.tex

To solve the advection-diffusion equation \eqref{eq:adu}--\eqref{eq:adbdry} we apply elliptic marching, also know as Rothe's method or method of lines \cite{HindmarshAC1993TNMo,ROMANCHAPKO1997RMFT}, meaning that we first discretize first in time, using an implicit treatment of the stiff diffusion term; the advective term and the source term are treated explicitly. Regardless of the specifics of the such discretization method, the result is a sequence of modified Helmholtz equations \eqref{eq:mh}, \eqref{eq:mhbc} to be solved for each time-step. For example, a first-order IMEX Euler method applied to the advection-diffusion equation results in
\begin{equation}
  \begin{aligned}
  \label{eq:imexad}
  \frac{\adunp\fp{\x}-\adun\fp{\x}}{\dtn} + \velfn\fp{\x}\cdot\nabla\adun\fp{\x} &= \diffconst\Delta\adunp\fp{\x}+\force_{n}(\x),\quad \x\in \domainnp\\
  \pd{\adunp\fp{\x}}{\hat{n}} &= \adbc_{n+1}\fp{\x},\quad \x\in \bdrynp,
  \end{aligned}
\end{equation}
where $\adun\fp{\x}=\adu(\tn,\x)$, $\velfn\fp{\x}=\velf(\tn,\x)$, $\force_{n}(\x) = \force\fp{\tn,\x}$,  $\domain_{n+1} = \domain(\tnp)$, $\adbc_{n+1}\fp{\x} = \adbc\fp{t_{n+1},\x}$, and
$\bdrynp = \bdry(\tnp)$. It can be formulated as the modified Helmholtz equation \eqref{eq:mh}, \eqref{eq:mhbc}, with $\mhu = \adunp$, $\mhaa = 1/\diffconst\dtn$ and $\mhrhs\fp{\x} = \adun\fp{\x}/\diffconst\dtn - \velfn\fp{\x}\cdot\nabla\adun\fp{\x} + \force_{n}(\x)$. Note that $f$ exists on $\domainn$, which may be different from $\domainnp$, and thus $f$ is not defined in all of $\domainnp$. For the remainder of this section we will assume $f$ to be defined both on $\domainn$ and $\domainnp$. We will return to the question of how such an $f$ is obtained in section \ref{sss:imexsdctimedep}.

The free-space Green's function for the differential operator $\mhaa - \Delta$ in $\mathbb{R}^{2}$ is
\begin{equation}
    \label{eq:green}
    \mathconstfont{G}\fp{\x,\y} = \frac{1}{2\pi}\besselK{0}\fp{\mha \norm{\x-\y}},
\end{equation}
where $\besselK{0}$ is the zeroth-order modified Bessel function of the second kind, and $\norm{\cdot}$ is the Euclidean norm unless specified otherwise. Note that neither the modified Helmholtz equation nor its Green's function are consistently named in the literature, but are also referred to as the screened Poisson equation, the linearized Poisson-Boltzmann equation, the Debey-Hückel equation and the Yukawa equation, and the Yukawa potential or screened Coulomb potential \cite{JUFFER1991144,CerioniAlessandro2012Eaas,ROWLINSON198915}.

By the linearity of the differential operator $\mhaa - \Delta$, the problem \eqref{eq:imexad} can be decomposed into a particular problem,
\begin{equation}
  \label{eq:mhp}
  \mhaa\mhpu\fp{\x}-\Delta\mhpu\fp{\x} = \mhrhs\fp{\x},\quad \x\in\domain,
\end{equation}
with free-space boundary conditions $\mhpu(\x)\rightarrow 0$ as $\norm{\x}\rightarrow \infty$, and a homogeneous problem
\begin{equation}
\begin{aligned}
  \label{eq:mhh}
  \mhaa\mhhu\fp{\x}-\Delta\mhhu\fp{\x} &= 0, \quad\x\in\domain,\\
  \pd{\mhhu\fp{\x}}{\normalscal} &= \mhbc\fp{\x}-\pd{\mhpu\fp{\x}}{\normalscal}, \quad\x\in\bdry.
\end{aligned}
\end{equation}
Thus one PDE includes the right-hand side $\mhrhs$, also referred to as the forcing term, while the other enforces the Neumann boundary conditions with a correction $\pdi{\mhpu\fp{\x}}{\normalscal}$. By construction $\adunp = \mhpu + \mhhu$ satisfies \eqref{eq:imexad}.

The structure outlined above can be obtained using any IMEX rule. The methods differ in the definitions of $\mhaa$ and $\mhrhs$, the required number of modified Helmholtz equations to solve per time-step, and how these results are combined to advance an approximate solution to the advection-diffusion equation in time. For example, with an IMEX R-K method the right-hand side $\mhrhs$ is a linear combination of intermediate so called stages. While this is the approach in  \cite{FryklundFredrik2020Aien} by \citeauthor{FryklundFredrik2020Aien}, where the heat equation is solved in a similar framework, here an IMEX SDC is used instead. The choice of SDC over R-K is discussed in section \ref{ss:timestepping}. The concept of SDC is to improve an initial approximation through solving a series of update equations with low-order methods, such as the first-order IMEX Euler method above. In the context of this paper, these update equations also take the form of the modified Helmholtz equation. 

\subsection{Particular solution}
Here we present the mathematical formulation of the solution to the particular problem \eqref{eq:mhp}. Leaving numerics aside, a mathematically straightforward approach is to express the solution as a convolution of the forcing term $\mhrhs$ and the free-space Green's function \eqref{eq:green} over $\domain$. Even though this formulation is simple, to numerically evaluate the integral accurately is complicated for complex geometries. Since we consider time-dependent geometries, high-order tailored quadrature rules for a specific geometry is not a feasible option. Therefore, we do as is proposed in \cite{FryklundFredrik2020Aien,AFKLINTEBERG2020109353}, and explicitly write the solution as a Fourier series, after first reformulating the problem.

Instead of the initial formulation of the particular problem \eqref{eq:mhp}, consider
  \begin{equation}
  \label{eq:mhpf}
  \mhaa\mhpu\fp{\x}-\Delta\mhpu\fp{\x} = \mhrhse\fp{\x},\quad \x\in\suppbox\subset\mathbb{R}^{2},
\end{equation}
where $\mhrhse\in C^{q}(\suppbox)$ for some $q>0$, $\bar{\domain}\subset \suppbox$ and
\begin{equation}
  \label{eq:mhrhsecond}
\begin{aligned}
 & \mhrhse\fp{\x}=\mhrhs\fp{\x},\quad \x \in \domain,\\
 & \supp{\mhrhse}\subset \suppbox.
\end{aligned}
\end{equation}
There is flexibility in choosing $\suppbox$. Thus, assuming the existence of an extension $\mhrhse$ with compact support enclosing $\bar{\domain}$, we may set $B$ to be a box $[-L/2,\,L/2]^{2}$ containing $\supp{\mhrhse}$, consider \eqref{eq:mhpf} as a periodic problem with periodic boundary conditions on $\suppbox$, and express the solution as a Fourier series. Standard results yield an explicit expression for the solution
\begin{equation}
  \label{eq:mhpfsol}
\mhpu\fp{\x} = \sum\limits_{\fovarvec}\frac{\fourier{\mhrhse}\fp{\fovarvec}}{\mhaa + \norm{\fovarvec}^{2}}e^{\I\fovarvec\cdot\x},\quad \x\in\domain,
\end{equation}
where $\{\fourier{\mhrhse}(\fovarvec)\}_{\fovarvec}$, $\fovarvec = \vec{n}2\pi/L$ with $\vec{n}\in\mathbb{Z}^{2}$, are the Fourier coefficients in the Fourier series expansion of $\mhrhse$.
Analogously, the gradient of the solution is
\begin{equation}
  \label{eq:mhpfgradsol}
 \nabla \mhpu\fp{\x} = \sum\limits_{\fovarvec}\I\fovarvec\frac{\fourier{\mhrhse}\fp{\fovarvec}}{\mhaa + \norm{\fovarvec}^{2}}e^{\I\fovarvec\cdot\x}, \quad\x\in\domain.
\end{equation}
Note that the zeroth-mode is well-defined for $\mha\ne 0$.

%%%% FIGURE %%%%
\begin{figure}
  \centering
  \includegraphics[width=0.4\textwidth]{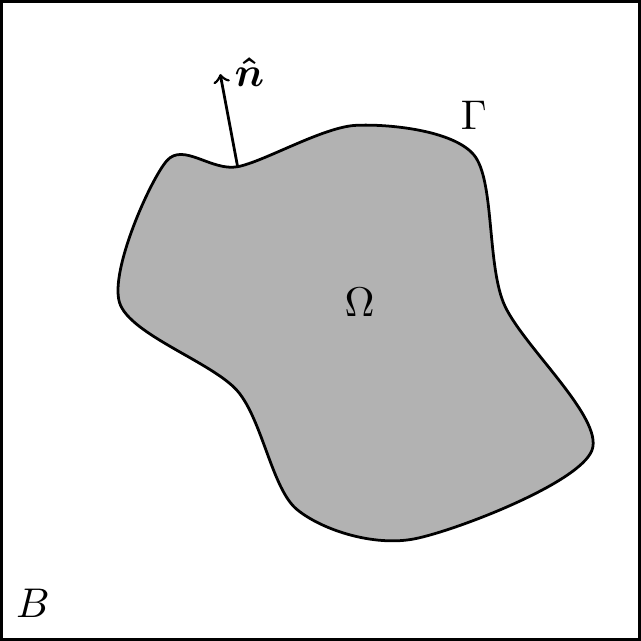}%     without .tex extension
  % or use \input{mytikz}{images/heateqimexr34/linferror.png}
\caption{A bounded domain $\domain$ in $\mathbb{R}^{2}$ with a smooth boundary $\bdry$, enclosed by the box \newline$\suppbox = [-L/2,L/2]^{2}$.}
\label{fig:Omega}
\end{figure}
%%%%%%%%%%%%%%%%%%%%%%%%%%%%%%%%%%%%%%%%%%%%%%%%5

\subsection{Homogeneous solution}
Here we consider the homogeneous problem \eqref{eq:mhh} in a boundary integral formulation which takes the geometry and the boundary data into account; the forcing term $\mhrhs$ is incorporated into the particular problem.
The homogeneous solution in the form of a single layer potential is denoted by
\begin{equation}
  \label{eq:mhhslp}
  \mhhu\fp{\x} = \frac{1}{2\pi}\int\limits_{\bdry}\kerdomain\fp{\y,\x}\bimdens(\y)\D s_{\y}, \quad \x\in\domain,
\end{equation}
and the gradient of the solution can be written as
\begin{equation}
  \label{eq:mhhslpgrad}
  \nabla\mhhu\fp{\x} = \frac{1}{2\pi}\int\limits_{\bdry}\nabla\kerdomain\fp{\y,\x}\bimdens(\y)\D s_{\y}, \quad \x\in\domain.
\end{equation}
The kernels are given by the Green's function \eqref{eq:green}, with the factor $1/2\pi$ extracted to outside the layer potential integrals, 
\begin{equation}
    \label{eq:mhhkerdomain} 
    \kerdomain\fp{\y,\x} = \besselK{0}\fp{\mha\norm{\y-\x}},
\end{equation}
and
\begin{equation}
        \label{eq:mhhkerdomaingrad} 
    \nabla\kerdomain\fp{\y,\x} = \mha\besselK{1}\fp{\mha\norm{\y-\x}}\frac{\y-\x}{\norm{\y-\x}},
\end{equation}
%\begin{multicols}{2}
%  \begin{equation}
%    \label{eq:mhhkerdomain} 
%    \kerdomain\fp{\x,\y} = \mhaa\besselK{0}\fp{\mha\abs{\y-\x}},
%  \end{equation}\break
%  \begin{equation}
%        \label{eq:mhhkerdomaingrad} 
%    \nabla\kerdomain\fp{\x,\y}=\mha^{3}\besselK{1}\fp{\mha\abs{\y-\x}}\frac{\y-\x}{\abs{\y-\x}},
%  \end{equation}
%\end{multicols}
where the gradient is with respect to $\x$, and where $\besselK{0}$ and $\besselK{1}$ are the zeroth- and first-order modified Bessel functions of the second order, respectively \cite{Kropinski2011modHelm}. Here, the layer density $\bimdens$ is assumed to be continuous and $\bimdens : \bdry\mapsto\mathbb{R}$. It is the solution to the second-kind Fredholm integral equation
\begin{equation}
  \label{eq:mhhbie}
  \frac{1}{2}\bimdens\fp{\x} + \frac{1}{2\pi}\int\limits_{\bdry}\kerbdry\fp{\y,\x}\bimdens(\y)\D s_{\y} = \mhbc\fp{\x}-\pd{\mhpu\fp{\x}}{\normalscal}, \quad \x\in\bdry,
\end{equation}
with the double layer kernel
\begin{equation}
  \label{eq:mhhkernel}
\kerbdry\fp{\y,\x} = \nabla\kerdomain\fp{\y,\x}\cdot\normal\fp{\x}  = \alpha\besselK{1}\fp{\mha\norm{\y-\x}}\frac{\y-\x}{\norm{\y-\x}}\cdot\normal\fp{\x}.
\end{equation}
The boundary integral equation \eqref{eq:mhhbie} is obtained by enforcing the Neumann boundary condition in the limit as $\x\in\domain$ goes towards the boundary and using the jump conditions given in \cite{Kropinski2011modHelm}.

We note for future reference that the limiting value of $\kerbdry$ on the boundary is
\begin{equation}
  \label{eq:mhhkerlim}
  \lim_{\x\rightarrow\y}  \kerbdry\fp{\y,\x} = -\frac{1}{2}\kappa\fp{\y},\quad \x,\,\y\in\bdry,
\end{equation}
where $\kappa$ is the curvature of $\bdry$ at $\y$, assuming $\bdry\in C^{2}$. The kernel $\kerbdry$ is thus bounded and continuous along $\bdry$, and therefore the integral operator is compact \cite{KressRainer2014LIE}. Since there is no nontrivial homogeneous solution, the Fredholm alternative states that   \eqref{eq:mhhbie} has a unique solution for any continuous data $\mhbc\fp{\x}-\partial\mhpu\fp{\x}/\partial\normalscal$ \cite{KreyszigErwin1989Ifaw}.
%%% Local Variables:
%%% mode: latex
%%% TeX-master: "manuscript.tex"
%%% End:

%% file: numericalmethods.tex
% NUMERICAL METHODS
%%%%%%%%%%%%%%%%%%%%%
In this section we present the numerical methods used to solve the advection-diffusion equation \eqref{eq:adu}--\eqref{eq:adbdry}. First we present the selected time-stepping scheme and how its application to the advection-diffusion equation results in a series of modified Helmholtz equations. Next, we demonstrate how to solve the modified Helmholtz equation on stationary domains. In section \ref{s:pux}, moving domains are introduced, in combination with the PUX method for function extension.

%%%%%%%%%%%%%%%%%%%%%%
\subsection{Time-stepping}
\label{ss:timestepping}
Here we describe an IMEX SDC method, based on the standard first-order IMEX Euler method, first for ordinary differential equations (ODE), and thereafter its application to the advection-diffusion equation. Stationary domains are considered first, thereafter time-dependent domains. In IMEX methods stiff terms are treated implicitly, in order to avoid prohibitively small time-steps \cite{AscherIMEX}. The non-stiff terms are treated explicitly, as is often more computationally efficient. 

\subsubsection{IMEX spectral deferred correction as applied to an ODE}
\label{sss:sdcgen}
Here a shortened version of the methodology in \cite{minion2003} is presented. Consider the initial value problem for a first-order ODE on standard form
\begin{align}
  \label{eq:ostform}
  \Dt{\ou} \fp{t} &= \orhs\fp{t,\ou\fp{t}}, \quad t\in [t_{a},t_{b}],\\
  \ou (t_{a}) &= \oui,
\end{align}
where the solution $\ou$ and a given initial value $\oui$ are in $\mathbb{R}$, and $F : \mathbb{R}\times\mathbb{R} \rightarrow \mathbb{R}$. It is assumed that $F$ is smooth and can be decomposed as
\begin{equation}
  \label{eq:orhsdecomp}
  \orhs\fp{t,\ou\fp{t}} = \orhse\fp{t,\ou\fp{t}} + \orhsi\fp{t,\ou\fp{t}},
\end{equation}
where $\orhse$ is a non-stiff term and $\orhsi$ is a stiff term. The subscripts denote that the stiff part is treated implicitly and the non-stiff part explicitly.

The interval $[t_{a},t_{b}]$ is subdivided into $p$ subintervals by choosing points $t_{a} = \tO < t_{1} < \ldots < t_{p} = t_{b}$, which are also referred to as instances (in time), and the length $\dtm = \tmp-\tm$ of a subinterval will be referred to as a substep.
Given $\oui$, initial (or provisional) approximate solutions $\ouakm{0}{m}\approx\ou\fp{\tm}$ are computed at the instances $\gc{\tm}_{m=1}^{p}$. To this end, introduce  $P = p + 1$ Gauss-Lobatto points and use the standard first-order IMEX Euler method for these points. Then, the initial approximate solution is the solution to 
\begin{equation}
  \label{eq:obfe}
    \ouakm{0}{m+1} = \ouakm{0}{m} + \dtm\gp{\orhse\fp{\tm,\ouakm{0}{m}}+\orhsi\fp{\tmp,\ouakm{0}{m+1}}}, \quad m = 0,1,\ldots,p-1,\quad  \ouakm{0}{0} = \oui.
\end{equation}
An increasingly accurate approximation  $\gc{\ouakm{k+1}{m}}_{m=0}^{p}$, for $k = 0,1,\ldots K-1$, is obtained through a sequence of corrections. The superscript denotes the number of times it has been corrected; e.g. $\ouakm{k}{m}$ is the $k$th corrected approximation of $\ou\fp{\tm}$. Each correction revises the current approximate solution and increases its formal order of accuracy.

A corrected approximation is computed as follows. Given $\gc{\ouakm{k}{m}}_{m=0}^{p}$ we proceed as in \cite{minion2003} and obtain
\begin{equation}
  \label{eq:ocoreqaimexfin}
  \begin{aligned}
    \ouakm{k+1}{m+1} &= \ouakm{k+1}{m} + \dtm\gp{\orhse(\tm,\ouakm{k+1}{m})-\orhse\fp{\tm,\ouakm{k}{m}}+\orhsi\fp{\tmp,\ouakm{k+1}{m+1}}-\orhsi\fp{\tmp,\ouakm{k}{m+1}}} \\ & + I_{m}^{m+1}\fp{\gc{\ouakm{k}{m}}_{m=0}^{p}},
      \end{aligned}
\end{equation}
for $m = 0,1,\ldots,p-1$, initialized at $\ouakm{k+1}{0} =\oui$. We refer to \eqref{eq:ocoreqaimexfin} as the update equation and an application of the update equation means sequentially solving it for $\gc{\ouakm{k+1}{m}}_{m=1}^{p}$. If the update equation is applied $K-1$ times, then we have in total $K$ provisional solutions at each instance $\gc{\tm}_{m=1}^{p}$ on $[t_{a},t_{b}]$. To obtain a $K$th-order method, it is sufficient to set $P=K$, i.e. $p = K-1$, as explained below. Finally, we accept $\ouakm{K-1}{K-1}$ as the approximation of  $\ou\fp{t_{b}}$.

In the right-hand side of \eqref{eq:ocoreqaimexfin} the expression  $I_{m}^{m+1}\fp{\gc{\ouakm{k}{m}}_{m=0}^{p}}$ appears, which is a $p$th-order numerical quadrature approximation of
\begin{equation}
  \label{eq:sdcintegral}
  \int\limits_{\tm}^{\tmp} \orhs\fp{\tau,\ouakm{k}{}\fp{\tau}}\D \tau.
\end{equation}
Note that $\ouakm{k}{}$ is a continuous function. This integral over $[\tm,\tmp]$ is computed by integrating the $p$th-degree polynomial that is interpolating $\ouakm{k}{}$ over $[t_{a},t_{b}]$, giving
\begin{equation}
  \label{eq:oquad}
  I_{m}^{m+1}\fp{\ouakm{k}{0},\ouakm{k}{1},\ldots,\ouakm{k}{p}} = \sum\limits_{j=0}^{p}q_{m}^{j}\orhs\fp{t_{j},\ouakm{k}{j}}, \quad m = 0,\ldots, p.
\end{equation}
The coefficients $q_{m}^{j}$ can be precomputed, reducing the quadrature to a simple matrix-vector multiplication. Each of the $p$ different quadrature rules only require a rescaling of nodes and weights from the basic interval to $[t_{a},t_{b}]$.%Note that there are $p$ quadrature rules; one for each subinterval. They can be precomputed and rescaled with $t_{b}-t_{a}$

The error in computing the integral \eqref{eq:sdcintegral} with \eqref{eq:oquad} is $\mathcal{O}\fp{\dt^{P+1}}$, assuming $\orhs$ is smooth enough. Now, each application of the update equation picks up at least a single order of accuracy \cite{CausleyMathew2019Otco}, thus to obtain a $K$th-order method $K-1$ corrections are sufficient. This implies that the number of Gauss-Lobatto points $P$ should be set to be equal to $K$, as proposed in \cite{minion2003}. Thus, the total number of implicit systems, i.e. \eqref{eq:obfe} and \eqref{eq:ocoreqaimexfin}, to solve over $[t_{a},t_{b}]$ is $K^{2} - K$.

To summarize. First we solve the given ODE \eqref{eq:ostform} for $\gc{\ouakm{0}{m}}_{m=1}^{K-1}$ at $\gc{\tm}_{m=1}^{K-1}$, given an initial value $\ouakm{0}{0}$ at $\tO$, with a standard first-order IMEX method. Then the corrected terms $\gc{\ouakm{1}{m}}_{m=1}^{K-1}$ are computed, initialized at $\ouakm{1}{0} = \ouakm{0}{0}$. This process is repeated $K-1$ times. Finally we obtain the approximation $\ouakm{K-1}{K-1}$ of $\ou\fp{t_{K-1}} = \ou\fp{t_{b}}$.
  
\subsubsection{IMEX spectral deferred correction for the advection-diffusion equation on stationary domains}

Consider the interval $[t_{a},t_{b}]$ as a time-step $[\tn,\tnp] = [\tn,\tn + \dtn]$, which is subdivided as described above. Assume that the domain $\domain$ is stationary. Writing the advection-diffusion equation \eqref{eq:adu}--\eqref{eq:adbdry} with the split into a stiff and a non-stiff term, as in \eqref{eq:orhsdecomp}, yields
\begin{equation}
  \orhse = -\velf \cdot \nabla\ou + \force,\quad \orhsi = \diffconst\Delta \ou.
\end{equation}
Inserted into \eqref{eq:obfe}, and collecting unknown terms on the left-hand side results in the modified Helmholtz equation \eqref{eq:mh}, with $\mhu = \adunp$, $\mhaa = 1/\diffconst\dtn$ and $\mhrhs\fp{\x} = \adun\fp{\x}/\diffconst\dtn - \velfn\fp{\x}\cdot\nabla\adun\fp{\x} + \force_{n}(\x)$. Analogously \eqref{eq:ocoreqaimexfin} becomes
\begin{equation}
  \begin{aligned}
  \label{eq:sdcadeq}
  \frac{1}{\diffconst\dtm}\ouakm{k+1}{m+1} - \Delta\ouakm{k+1}{m+1} &= \frac{1}{\diffconst\dtm}\ouakm{k+1}{m} + \frac{1}{\diffconst}\fp{\velf^{k+1}_{m}\cdot\nabla\ouakm{k+1}{m} -\velf^{k}_{m}\cdot\nabla\ouakm{k}{m}-\diffconst\Delta\ouakm{k}{m+1}}\\ &+ \frac{1}{\diffconst\dtm}I_{m}^{m+1}\fp{\gc{\ouakm{k}{m}}_{m=0}^{p}}.
  \end{aligned}
\end{equation}
Recast as the modified Helmholtz equation \eqref{eq:mh} we have
\begin{equation}
  \label{eq:mhsdcterms_u_a}
  \mhu =  \ouakm{k+1}{m+1}, \quad  \mha = \sqrt{\frac{1}{\diffconst\dtm}},
\end{equation}
and
\begin{equation}
  \label{eq:mhsdcterms_f}\mhrhs = \frac{1}{\diffconst\dtm}\ouakm{k+1}{m} + \frac{1}{\diffconst}\fp{\velf^{k+1}_{m}\cdot\nabla\ouakm{k+1}{m}-\velf^{k}_{m}\cdot\nabla\ouakm{k}{m}-\diffconst\Delta\ouakm{k}{m+1}} + \frac{1}{\diffconst\dtm}I_{m}^{m+1}\fp{\gc{\ouakm{k}{m}}_{m=0}^{p}}.
\end{equation}
For both formulations of the modified Helmholtz equation we enforce the Neumann boundary conditions $\mhbc$, as evaluated at $\tmp$.

Thus at each substep $(m,k)$ both the potential $\nabla\ouakm{k+1}{m+1}$ and the Laplacian $\Delta\ouakm{k+1}{m+1}$ are computed. The former can be computed alongside $\mhu$ by \eqref{eq:mhpfgradsol} and \eqref{eq:mhhslpgrad}, and the latter can be obtained directly from \eqref{eq:mh} once $\mhu$ is computed. These values have to be given as initial conditions as well at the start of each substep, as we need to initialize IMEX SDC with $\orhse\fp{\tO,\ouakm{0}{0}}$ and $\orhsi\fp{\tO,\ouakm{0}{0}}$ for each time-step.
\iffalse
As with IMEX Runge-Kutta methods, IMEX SDC methods may suffer from order reduction, which means that the observed order of convergence is lower than the theory suggests \cite{KENNEDY2003139,minion2003}. Especially using a method of lines (MOL) approach, with time dependent boundary conditions imposed at intermediate stages, may give order reduction \cite{minion2003}. \todo{Look more into this.}
\fi
\begin{remark}
High order IMEX methods can also be used to compute provisional solutions in SDC, such as R-K methods. However, even though the solution in R-K methods is a linear combination of stage values, it is not necessarily true for its gradient and Laplacian. Some R-K methods have a property sometimes called first-same-as-last (FSAL), meaning that the final implicit stages are the same as the first ones in the next iteration, thus reducing the number of implicit problems to solve by one per time-step. The explicit stages are assumed to have closed forms, in which case FSAL is superfluous for the explicit part. In the method proposed in this paper, the explicit part contains $-\velf\cdot\nabla \adu$ which is not explicit in a true sense; it is the gradient of the solution, which has to be solved for; or the differential operator has to be discretized, which is hard to do to high order. To the extent of our knowledge, only first-order IMEX Euler is FSAL also for the explicit term.
\end{remark}

\subsubsection{IMEX spectral deferred correction for the advection-diffusion equation on time-dependent domains}
\label{sss:imexsdctimedep}
Here we discuss how to handle time-dependent and deformable domains when solving the advection-diffusion equation \eqref{eq:adu}--\eqref{eq:adbdry}. The problem we address is the following. Assume there is a time-dependent domain $\domain(t)$ superimposed on a discrete stationary grid $\X$. When the domain $\domain(t)$ changes over time there may be discrete points from $\X$ entering the domain where data from previous time-levels do not exist. These points are in the parts of the new domain that do not overlap the previous one, see figure \ref{fig:slice}. To solve the modified Helmholtz equation at this instance in time, all data need to be known at all discrete points inside the new domain. We take the approach of accommodating for the missing data by extrapolating existing data to these points using PUX, as described in section \ref{s:pux}.

First we present how we update the domain in time. We consider a Lagrangian specification of the flow field $\velf$. The position of a point $\x$ in a time-dependent domain $\domain(t)$, or on its boundary $\bdry(t)$, changes in time as
\begin{equation}
  \label{eq:odebdry}
  \od{\x\fp{t}}{t} = \velf(t,\x(t)),\quad \x\in\bar{\domain}\fp{t}, \quad \tO < t.
\end{equation}
In our examples, the velocity field $\velf$ will either be given explicitly, or obtained by numerical solution of the Stokes equations. In the latter case, the Stokes equations are solved with a boundary integral method, similar to the one presented in this paper, using the same discretization of the boundary. See the papers \cite{OjalaRikard2015Aaie,BystrickyLukas2020Aaie}.

The update of $\x$ in \eqref{eq:odebdry} can be computed analytically in some simple cases, otherwise it is computed numerically, using an explicit R-K method of the same order as the IMEX SDC for the advection-diffusion equation. Clustering of the discretization nodes on the boundary is avoided by modifying the tangential velocity components \cite{OjalaRikard2015Aaie}.
%\begin{remark}
%  If updating $\x\in\bar{\domain}\fp{t}$ is done numerically, then it is preferable to apply the same time-marching scheme as for the governing equaions, albeit not necessary, as long as the accuracy in updating the domain is greater than in updating the solution to the governing equations. For example, applying the first-order IMEX Euler to \eqref{eq:odebdry} gives
%\begin{equation}
%  \x_{n+1} = \x_{n} + \dtn\velf_{n}\fp{\xn},\quad \xn\in\bar{\domain}\fp{t}, \quad \tO < t,
%\end{equation}
%with $\velf_{n}\fp{\xn} = \velf\fp{\tn,\xn}$.
%\end{remark}

We consider an IMEX SDC applied to an interval $[\tn,\tnp]$, with the instances $\{t_{j}\}_{j =0}^{p}$, where $t_{0} = \tn$ and $t_{p} = \tnp$, and $\dtm = \tmp -\tm$ for $m = 0,\ldots,p-1$. When addressing moving domains it is useful to see $\mhrhs$ from \eqref{eq:mhsdcterms_f} as a sum over data $\mhrhs_{j}$ on different time-levels, i.e. $\mhrhs = \sum_{j=0}^{p}  \mhrhs_{j}$, where
\begin{equation}
  \label{eq:mhrhsj}
\mhrhs_{j}  =  \frac{ q_{m}^{j}}{\diffconst\dtm}(-\velf^{k}_{j}\cdot\nabla\ouakm{k}{j}+\force_{j}+\diffconst\Delta\ouakm{k}{j}) +
 \frac{1}{\diffconst}\begin{dcases}
    \frac{1}{\dtm}\ouakm{k+1}{m} +\velf^{k+1}_{m}\cdot\nabla\ouakm{k+1}{m}-\velf^{k}_{m}\cdot\nabla\ouakm{k}{m}, \quad j = m,\\
    \diffconst\Delta\ouakm{k}{m+1}, \quad j = m+1,\\
    0, \quad \text{ otherwise},
    \end{dcases}
\end{equation}
for a given $m$ and $k$. Thus the different $\mhrhs_{j}$ correspond to different domains $\domain_{j}$, with $j= 0,\ldots,p$ at $\tm$ for $m = 0,\ldots,p$. Solving for $\ouakm{k+1}{m+1}$ in \eqref{eq:sdcadeq} means solving the modified Helmholtz equation in $\domain_{m+1}\cap\X$, where only $\mhrhs_{m+1}$ is entirely known. The other terms $\mhrhs_{j}$ are only partially known and need to be extrapolated to the points in $\domain_{m+1}\cap\X$ where data does not exist, otherwise $f$ will be incomplete. This is visualized in figure \ref{fig:domaintimelevels}, which demonstrates the application of IMEX SDC of order three to an one-dimensional domain translated from the left to the right. Here we solve the governing equations in $\bar{\Omega}_{1}\cap\X$, thus the data $\mhrhs_{0}$ and $\mhrhs_{2}$ need to be extended from $\Omega_{0}\cap\X$ and $\Omega_{2}\cap\X$, see figure \ref{fig:domaintimelevels}. 

This motivates the definition of a slice $\slice{t_a}{t_b}$. It consists of the points in $\X$ that are in the relative complement of $\domain(t_a)$ in $\domain(t_b)$, i.e. $\slice{t_a}{t_b} = (\domain(t_b)\setminus\domain(t_a))\cap\X$ for $\tO<t_a<t_b$, see figure \ref{fig:slice}. Then for each $\mhrhs_{j}$, with $j \neq m+1$, no data exists in $\slice{t_{j}}{t_{m+1}}$, we say that data is missing, and must be created artificially by extending it from $\domain_{j}\cap\X$ to $\slice{t_{j}}{t_{m+1}}$. A slice is naturally thin, typically of the size of $\dx$, thus the data is only extended a short distance.

We expect the procedure of extending data to the slices to impose a CFL-type condition of the form
\begin{equation}
  \label{eq:CFL}
\max_{\x\in\bar{\domain}}\norm{\velf(\x)}\dt/\dx\leq C.
\end{equation}
The intuition is that the method becomes unstable if the slices are too large, i.e. the domain moves too much in a single time-step, relative the resolution of the underlying grid $\X$.

\def\firstcircle{(0,0) circle (2cm)}
\def\secondcircle{(0:0.6cm) circle (2cm)}

\colorlet{circle edge}{blue!90}
\colorlet{circle area}{blue!20}

\tikzset{filled/.style={fill=circle area, draw=circle edge, thick},
    outline/.style={draw=circle edge, thick}}
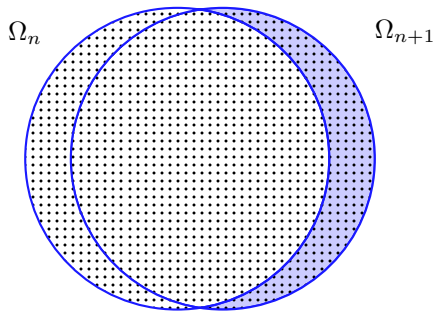
\begin{figure}
    \centering
\begin{tikzpicture}
    \begin{scope}
        \clip \secondcircle;
        \draw[filled, even odd rule] \firstcircle
                                     \secondcircle node {};
    \end{scope}
    \draw[outline,pattern=dots] \firstcircle node {}
                   \secondcircle;
    \node[anchor=south] at (-2,1.4) {$\domain_{n}$};
    \node[anchor=south] at (3,1.4) {$\domain_{n+1}$};
\end{tikzpicture}
\caption{The domain $\domain_{n}$ (left circle) at time $\tn$ is advanced in time to $\domain_{n+1}$ (right circle) at time $\tnp$. The shaded blue area is the slice $\slice{n}{n+1}$ where data in $\domain_{n+1}$ is missing. Note that the size of the blue area is larger than for a typical time-step.}
\label{fig:slice}
\end{figure}

\usetikzlibrary{decorations.markings}
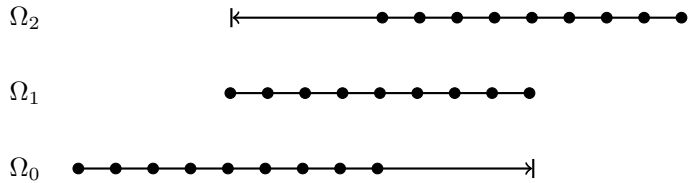
\begin{figure}
    \centering
    \begin{tikzpicture}[decoration={markings,
      mark=between positions 0 and 1 step 14pt
      with { \draw [fill] (0,0) circle [radius=2pt];}}]
    \path[postaction={decorate}] (0,0) to (4,0);
\draw [thick] (0, 0) -- (4,0);
\draw[thick,->|]        (4,0)   -- (6,0);

\path[postaction={decorate}] (2,1) to (6,1);
\draw [thick] (2, 1) -- (6,1);

\path[postaction={decorate}] (4,2) to (8,2);
\draw [thick] (4, 2) -- (8,2);
\draw[thick,->|]        (4,2)   -- (2,2);
\node[] at (-0.7,2) {$\domain_{2}$};
\node[] at (-0.7,1) {$\domain_{1}$};
\node[] at (-0.7,0) {$\domain_{0}$};
\end{tikzpicture}
\caption{A visualization of a three overlapping one-dimensional domains $\domain_{0}$,  $\domain_{1}$ and  $\domain_{2}$ at different time-levels. Data is available at the black dots.}
\label{fig:domaintimelevels}
\end{figure}

%%%%%%%%%%%%%%%%%%%%%%%%%%%%%%%%%%%%%%%%%%%%%%
%%%%%%%%%%%%%%%%%%%%%%%%%%%%%%%%%%%%%%%%%%%%%%
\subsubsection{Adaptivity}
\label{sss:adaptivity}
For an IMEX SDC it is by construction simple to obtain an approximation of the $\ell_{2}$ error at $\tm$ after $k+1$ corrections, namely $r = \|\ouakm{k+1}{m}-\ouakm{k}{m}\|_{2}$ \cite{Dutt2000}. If the estimate $r$ of the error is less than some tolerance $\epsilon_{\mathrm{TOL}}$, then $\ouakm{k+1}{m}$ is accepted as solution at time $\tnp$. The time-step is then updated to $\dt_{n,\mathrm{NEW}}$ by
\begin{equation}
\label{eq:adaptiveScheme}
\delta t_{n,\mathrm{NEW}} = \delta t_{n,\mathrm{OLD}}  (0.9\,\epsilon_{\mathrm{TOL}}/r)^{\frac{1}{K+1}},
\end{equation}
with $K$ being the order of the method and the value $0.9$ is a safety factor. If the solution is accepted, then $\delta t_{n+1} = \delta t_{n,\mathrm{NEW}}$, otherwise the computations start over at $t_{N}$ with $\delta t_{n} = \delta t_{n,\mathrm{NEW}}$. Thus even if the solution is accepted the step size is updated by the scheme \eqref{eq:adaptiveScheme}, meaning that the time-step size may increase if appropriate.

%However, the SDC is known to be slow to get into the asymptotic region, 
\subsection{Evaluating the particular solution}
Recall the periodic formulation of the particular problem \eqref{eq:mhpf}. It assumes the existence of an extension $\mhrhse$ satisfying \eqref{eq:mhrhsecond}, that is $\mhrhse = \mhrhs$ in $\domain$ and $\supp{\mhrhse}\subset \suppbox$, which is constructed with PUX as described in section \ref{s:pux}. The approximate solution is expressed as a truncated Fourier series. If the data is represented on a uniform grid on a geometrically simple domain, then it allows for efficient use of FFTs to compute the Fourier coefficients. To this end, consider a box $\suppbox=[-\tfrac{L}{2},\tfrac{L}{2}]^{2}$ in $\mathbb{R}^{2}$ that contains $\bar{\domain}$. Let $\X$ denote a set with $\Nu^{2}$ uniformly distributed elements $\x = (x_{1},x_{2})$, where $x_{i} = -L/2+n\dx$ with $n = 0,\ldots,\Nu-1$ and  $\dx = L/(\Nu-1)$. These elements are referred to as nodes or points. The solution to a modified Helmholtz equation is computed at all points that fall inside $\domain$. It is straightforward to compute the particular solution \eqref{eq:mhpfsol} and its gradient \eqref{eq:mhpfgradsol} on $\X$ with standard FFT libraries. It involves two applications of a FFT, each of complexity $\mathcal{O}\fp{\Nu^{2}\log \Nu^{2}}$. The order of convergence is $q + 2$ for the Fourier coefficients of $\mhrhse\in C^{q}(\suppbox)$ \cite{FRYKLUNDPUX,spectralMethodsInMatlab}.

The gradient of the particular solution $\nabla\mhpu$ is a correction in the boundary data for the homogeneous problem, see \eqref{eq:mhh} and \eqref{eq:mhhbie}. It is evaluated on the boundary $\bdry$ with a nonuniform FFT (NUFFT) \cite{LeslieGreengard2004AtNF}. If $\Nbdry$ is the number of evaluation points on the boundary, then computing the particular solution on these points by NUFFT is of complexity  $\mathcal{O}\fp{\Nu^{2}\log \Nu^{2} + \Nbdry}$, where $\Nbdry$ is the number of discretization points on the boundary $\bdry$, introduced in the next section.
%\begin{remark}
%  A volume FMM is an option to FFT. Then, assuming the existance of an extension $\mhrhse$, the volume potential is
%  \begin{equation}
%    \label{eq:mhpuvolpot}
%    \mhpu\fp{\x} = \frac{\mhaa}{2\pi}\int\limits_{\Omega}\besselK{0}\fp{\mha\abs{\y-\x}}\mhrhs\fp{\y}\D A_{\y}= \frac{\mhaa}{2\pi}\int\limits_{\mathbb{R}^{2}}\besselK{0}\fp{\mha\abs{\y-\x}}\mhrhse\fp{\y}\D \y.
%  \end{equation}
%  This formulation is favorable over the one proposed in this paper, as it allows an spatially adaptive grid. However, the current version of PUX requires a uniform grid for efficiency, making FFTs an obvious choice. See section \ref{s:movingdomains}.
%\end{remark}  
\subsection{Computing the homogeneous solution}
\label{ss:comphomo}
Here we present the layout for solving the boundary integral equation \eqref{eq:mhhbie} for the layer densities $\bimdens$ and how to compute the layer potentials \eqref{eq:mhhslp} and \eqref{eq:mhhslpgrad} as a post-processing step.

Partition the boundary into $\Npan$ segments, which we refer to as panels, such that
\begin{equation}
  \bdry = \bigcup\limits_{i=1}^{\Npan}\panel_{i}.  
\end{equation}
For each panel there exists a parametrization
\begin{equation}
\panel_{i} = \gcv{\bm{\gamma}_{i}\fp{t}\,\vert\, t\in \gb{-1,1}},
\end{equation}
and each panel is discretized in the parameter $t$ with a $\Nq$-point Gauss-Legandre quadrature rule, using the canonical quadrature nodes $\{\tGj\}_{j = 1}^{\Nq}$ and quadrature weights $\{\wGj\}_{j = 1}^{\Nq}$. Then, a layer potential for some continuous kernel $\lpker$ over a single panel is approximated as
\begin{equation}
  \label{eq:genlayerpot}
  \int\limits_{\panel_{i}}\lpker\fp{\y,\x}\bimdens(\y)\D s_{\y} = \int\limits_{-1}^{1}\lpker\fp{\y\fp{\bm{\gamma}_{i}\fp{t}},\x}\bimdens(\y\fp{\bm{\gamma}_{i}\fp{t}})\abs{\bm{\gamma}^{\prime}_{i}(t)}\D t \approx  \sum\limits_{j = 1}^{\Nq}\lpker\fp{\y_{ij},\x}\bimdens_{ij}s_{ij}\wGj,
\end{equation}
where $\y_{ij}= \y(\bm{\gamma}_{i}(\tGj))$, $\bimdens_{ij}$ is an approximation of $\bimdens(\y_{ij})$ and $s_{ij}= |\bm{\gamma}^{\prime}_{i}(\tGj)|$. Also, it will be useful to introduce $\normal_{ij} = \normal\fp{\x_{ij}}$. If the geometry and integrand are well-resolved, then the approximation \eqref{eq:genlayerpot} is accurate to order $2\Nq$. The kernels  \eqref{eq:mhhslp} and \eqref{eq:mhhslpgrad} are singular and nearly singular unless the target point $\x$ is well-separated from the source points $\{\y_{j}\}_{j = 1}^{\Nq}$, resulting in a significant loss of accuracy. It can be regained by applying corrections obtained by product integration methods, as explained below.

We note that a parametrization of the boundary is only needed for initialization of the advection-diffusion equation solver. As the domain evolves in time it is the discretization that is updated; this applies to all methods we employ to solve the advection-diffusion equation in this paper. We use the approach in \cite{PALSSON2019218} to numerically update the discretization of the moving boundaries.

Discretizing \eqref{eq:mhhbie} and collocating at the $\Nbdry=\Npan\Nq$ nodes, known as the Nystr\"{o}m method, gives the linear system
\begin{equation}
  \label{eq:mhhbiedisc}
  \frac{1}{2\pi}\bimdens_{i^\prime j^\prime} +   \frac{1}{2\pi}\sum\limits_{i = 1}^{\Npan} \sum\limits_{j = 1}^{\Nq}\kerbdry\fp{\y_{ij},\x_{i^\prime j^\prime}}\bimdens_{ij}s_{ij}\wGj =  \mhbc\fp{\x_{i^\prime j^\prime}}-\pd{\mhpu\fp{\x_{i^\prime j^\prime}}}{\normalscal_{i^\prime j^\prime}}, \quad i^\prime = 1,\ldots, \Npan,\,j^\prime = 1,\ldots,\Nq,
\end{equation}
with $\kerbdry$ given by \eqref{eq:mhhkernel}. As for the original boundary integral formulation, the Fredholm alternative states that there exists a unique solution. The corresponding system matrix is dense, but since its eigenvalues have only one accumulation point, which is away from the origin, the layer density  $\bimdens$ can be solved for efficiently with GMRES \cite{TrefethenLloydN1997Nla}. The condition number and the dimensionality of the Krylov subspace are independent of the resolution of the boundary, assuming the density $\bimdens$ and the geometry are sufficiently resolved. However, the condition number does depend on the curvature \cite{KressRainer2014LIE,KreyszigErwin1989Ifaw}.

After solving for the density $\bimdens$, the solution $\mhhu$ and its gradient $\nabla\mhhu$ can be computed by
\begin{equation}
  \label{eq:mhhslpdisc}
  \mhhu\fp{\x} = \frac{1}{2\pi}\sum\limits_{i = 1}^{\Npan} \sum\limits_{j = 1}^{\Nq}\kerdomain\fp{\y_{ij},\x}\bimdens_{ij}s_{ij}\wGj, \quad \x\in\domain,
\end{equation}
and
\begin{equation}
  \label{eq:mhhslpgraddisc}
  \nabla\mhhu\fp{\x} = \frac{1}{2\pi}\sum\limits_{i = 1}^{\Npan} \sum\limits_{j = 1}^{\Nq}\nabla\kerdomain\fp{\y_{ij},\x}\bimdens_{ij}s_{ij}\wGj, \quad \x\in\domain,
\end{equation}
with the kernels \eqref{eq:mhhkerdomain} and \eqref{eq:mhhkerdomaingrad}. Computing the sums in \eqref{eq:mhhbiedisc}, \eqref{eq:mhhslpdisc} and \eqref{eq:mhhslpgraddisc} can be accelerated with spectral Ewald summation, see section \ref{ss:ewald}. Then, the complexity is of order $\mathcal{O}(\Nbdry\log \Nbdry)$ for solving the boundary integral equation, and $\mathcal{O}(\Nu^{2}\log \Nu^{2})$ to evaluate the solution on the uniform grid. Note that using standard Gauss-Legendre quadrature introduces errors due to the presence of near singularities when a target point is not well-separated from the source points. The numerical solution is corrected with kernel-split quadrature with product integration as explained below.
%\todo{To A-K: Is this correct when targets $\neq$ sources?} 
\subsubsection{Kernel-split quadrature with product integration}
Standard Gauss-Legendre quadrature methods are only accurate if sources and targets are well-separated, otherwise the kernels \eqref{eq:mhhbiedisc}, \eqref{eq:mhhslpdisc} and \eqref{eq:mhhslpgraddisc} are singular or nearly singular, which in turn degrades the accuracy in the approximate solution. To regain the lost accuracy we apply a kernel-split product integration method for the modified Helmholtz equation, introduced in \cite{klinteberg2019adaptive}. It is an extension of \cite{specialQuad}, which is a quadrature method based on decomposing the integrand into explicit singular terms multiplying smooth functions, and applying product integration; meaning that the smooth functions are approximated by polynomials and each term is integrated separately. This original method \cite{specialQuad} is efficient and accurate for several PDEs, such as the Helmholtz equation and the Stokes equations \cite{HelsingJohan2015Voae,OjalaRikard2015Aaie}. However, as we explain later in this section, it may fail completely for the modified Helmholtz equation, hence the need for a modified version.

%However, the modified Helmholtz equation (and the modified Stokes equations) has layer potentials that depends on $\alpha$ and contain a log-type singularity multiplying $\besselI{0}\fp{\mha r}$ or $\besselI{1}\fp{\mha r}$, with $r \geq 0$. Both have the asymptotic behavior $\exp{\mha r }/\sqrt{2\pi \mha r}$ as $\mha r$ goes to infinity, and thus are increasingly poor candidates for polynomial interpolation as $\mha$ grows. It is bound to happen; recall that $\mha$ is inverse proportional to the time-step. Also, there are numerical cancellation errors, as discussed below.  Thus high resolution in time requires unfeasible high resolution of the boundary. On the other hand, only panels close to a target point will introduce these errors and local corrections for these points are sufficient to restore the lost accuracy. Therefore, we take the approach of local refinement of the panels, on which standard product integration is used. This local refinement is only part of the numerical integration procedure and does not affect the degrees of freedom in the discretization.

To apply the kernel-split quadrature in \cite{klinteberg2019adaptive,specialQuad} the kernels  \eqref{eq:mhhkernel},  \eqref{eq:mhhkerdomain} and \eqref{eq:mhhkerdomaingrad} must be split in order to identify the character of singularities. All three contain $\besselK{0}$ or $\besselK{1}$ which by standard decompositions \cite[\S10.31]{NIST:DLMF} can be written as 
\begin{align}
\besselK{0}\fp{\norm{\y - \x}} &=  -\besselI{0}\fp{\norm{\y - \x}}\log\fp{\norm{\y - \x}} + \besselK{0}^{S}\fp{\norm{\y - \x}},\quad \y\ne \x,\\
\besselK{1}\fp{\norm{\y - \x}} &= \frac{1}{\norm{\y - \x}} + \besselI{1}\fp{\norm{\y - \x}}\log\fp{\norm{\y - \x}} + \besselK{1}^{S}\fp{\norm{\y - \x}},\quad \y\ne\x,
\end{align}
where smooth remainder terms are collected in $\besselK{0}^{S}$ and $ \besselK{1}^{S}$. For on surface evaluation we have
\begin{equation}
\kerbdry\fp{\y,\x}   =  \kerbdrylog\fp{\y,\x}\log\fp{\norm{\y_{ij}-\x}} + \kerbdryself\fp{\y,\x}, \quad \x,\y\in\bdry,
\end{equation}
with
\begin{equation}
\kerbdrylog\fp{\y,\x} = \alpha\besselI{1}\fp{\mha\norm{\y-\x}}\frac{\y-\x}{\norm{\y-\x}}\cdot\normal\fp{\x},
\end{equation}
which goes to zero in the limit as $\x$ goes to $\y$. By \eqref{eq:mhhkerlim} we have
\begin{equation}
  \label{eq:limitkappa}
  \lim\limits_{\x\rightarrow\y}\kerbdry\fp{\y,\x} = \kerbdryself\fp{\y,\y} = -\frac{1}{2}\kappa\fp{\x}, \quad \x,\,\y\in\bdry.
\end{equation}
%The term is nearly-singular, as it has a limit value if a target point $\x$ goes to a point $\y$ along the boundary. Still, direct computation of the log-term will be increasingly inaccurate as its argument decreases.
For $\kerdomain$ we have
\begin{align}
  \kerdomain\fp{\y,\x} = \kerdomainlog\fp{\y,\x}\log\fp{\,\norm{\y-\x}} + \kerdomainself\fp{\y,\x}, \quad \x\in\domain,\,\y\in\bdry.
\end{align}
% It has a log-type near singularity as a target point $\x$ in the domain goes to the boundary.
Here
\begin{equation}
  \label{eq:kerdomainlog}
\kerdomainlog\fp{\y,\x} = -\besselI{0}\fp{\mha\norm{\y-\x}},
\end{equation}
and we have no finite limit for $\kerdomain\fp{\y,\x}$ corresponding to \eqref{eq:limitkappa}. Again the remaining terms are collected in $\kerdomainself$, which is smooth.
For the kernel $\nabla\kerdomain$ the decomposition is
\begin{align}
\nabla\kerdomain\fp{\y,\x} = \nabla\kerdomainlog\fp{\y,\x}\log\fp{\norm{\y-\x}} + \nabla\kerdomaincau\fp{\y,\x}\frac{\y-\x}{\norm{\y-\x}^{2}}  + \nabla\kerdomainself\fp{\y,\x}, \quad \x\in\domain,\,\y\in\bdry,
\end{align}
where the symbol $\nabla$ on the right-hand side is not the gradient, but part of the factors' names. Here, $\nabla\kerdomainself$ is a smooth function, $\nabla\kerdomaincau = 1$, and 
\begin{equation}
  \nabla\kerdomainlog\fp{\y,\x} = \alpha\besselI{1}\fp{\mha\norm{\y-\x}}\frac{\y-\x}{\norm{\y-\x}}.
\end{equation}
Again, this has a log-type near singularity, but also a Cauchy-type near singularity.

In all three cases above the smooth terms are computed using regular Gauss-Legendre quadrature, but product integration is needed for all singular or nearly-singular evaluation in order to compute the layer potentials accurately. The method in its entirety is explained in \cite{klinteberg2019adaptive}. 

We now demonstrate why product integration from \cite{specialQuad} may be insufficient for the modified Helmholtz equation, and to that end only present a simple case with the log-kernel. Thus, consider a panel with end points $-1$ and $1$, following the real line. The goal is to compute 
\begin{equation}
  \label{eq:unitpanelint}
  \int\limits_{-1}^{1}\sigma_{\mha}(y)\log\fp{\,\abs{y-x}}\D y = \Re\gpv{\int\limits_{-1}^{1}\sigma_{\mha}(y)\log\fp{y-x}\D y}, \quad x \in \mathbb{C},
\end{equation}
using the definition of the complex logarithm. Here, $\sigma_{\mha}$ is a smooth function that depends on the parameter $\mha$. Instead of applying standard quadrature methods, we expand the smooth function in a monomial basis. That is,
\begin{equation}
  \label{eq:polyinterp}
  p_{\Nq-1}(y) = \sum\limits_{k = 0}^{\Nq-1} c_{k}y^{k} \quad \text{where}  \quad p_{\Nq-1}(y_{j}) = \sigma_{\mha}(y_{j}), \quad j = 1,\ldots,\Nq,
\end{equation}
and we have
\begin{equation}
  \int\limits_{-1}^{1}\sigma_{\mha}(y)\log(y-x)\D y = \sum\limits_{k=0}^{\Nq-1}c_{k}\int\limits_{-1}^{1}y^{k}\log(y-x)\D y + R_{\Nq} (x),
\end{equation}
with the remainder $R_{\Nq} = \mathcal{O}(\Nq)$. The involved integrals can be computed analytically, given by recursion formulas.

Now let $\sigma_{\mha}(r) = \besselI{0}\fp{r\mha}$, as in \eqref{eq:kerdomainlog}, which grows asymptotically as $\exp{\mha r }/\sqrt{2\pi \mha r}$, where $\mha$ is proportional to the inverse square root of the time-step. For large $\mha$, i.e. small time-steps, the function $\sigma_{\alpha}$ may not be accurately represented by  \eqref{eq:polyinterp} with e.g. a $15$th degree polynomial, or even a $31$st degree polynomial, unless the maximum value of $r$ is small. Consequently, product integration does not produce an accurate correction. A naive approach is global refinement of the boundary, however both $\besselK{0}$ and $\besselK{1}$ decay asymptotically as  $\sqrt{\pi/(\mha r)}e^{-\mha r}$, i.e. they are very localized and only a small portion of the boundary needs to be refined. This would carry an unnecessary computational cost. Hence, the method in \cite{specialQuad} needs to be modified.

In \cite{klinteberg2019adaptive}, new panels are created  through recursive bisection such that they are of suitable length for product integration. It does not add degrees of freedom, but merely interpolates known quantities to shorther panels. The Cauchy-type singularity is also computed on the refined panels, as it increases the accuracy. Note that its multiplying function is a constant, thus product integration is not affected by increasing $\mha$.

The recursive bisection algorithm also addresses another problem: $\kerbdrylog\fp{\y,\x}\log\fp{\|\y-\x\|}$ and $\kerdomainlog\fp{\y,\x}\log\fp{\|\y-\x\|}$ grow with opposing signs to $\kerbdryself\fp{\y,\x}$ and $\kerdomainself\fp{\y,\x}$, both faster than $\besselK{0}$ and $\besselK{1}$ decrease. Thus they must cancel, potentially introducing cancellation errors when in finite precision.

%While the above case features the dynamics of product integration, several complicating factors are left out. While a panel can always be shifted and rotated such that its endpoints are at $-1$ and $1$, it does not have to follow the real axis. Product integration is formulated in terms of complex analysis, but here these reformulations are left out. Moreover, we have not discussed possible branch cuts and how to obtain the polynomial coefficients in a stable way.

\subsection{Fast summation}
\label{ss:ewald}
Discretizing the integrals \eqref{eq:mhhbie}, \eqref{eq:mhhslp} and \eqref{eq:mhhslpgrad} that arise from a boundary integral formulation yield the discrete sums \eqref{eq:mhhbiedisc},  \eqref{eq:mhhslpdisc} and \eqref{eq:mhhslpgraddisc}. Each can be computed directly at a cost $\mathcal{O}(N^{2})$ for $N$ sources and targets, which becomes costly for large $N$, e.g. when solving equation \eqref{eq:mhhbiedisc} with $N = \Npan\Nq$. In addition we need to compute several such sums in each time-step and inside each GMRES solve, where the number of iterations is independent of $N$. When solving the advection-diffusion equation with an IMEX SDC the $\mathcal{O}(N^{2})$ cost can be decreased by fast summation methods, such as an FMM \cite{doi:10.1137/0909044} or a spectral Ewald summation \cite{AfKlintebergLudvig2014FEsf,LindboDag2011Saif,LindboDag2010Safs,afKlintebergLudvig2017FEsf,PalssonSara2020Aiem}, to $\mathcal{O}(N)$ or $\mathcal{O}(N\log N)$.

The references above are all for the kernels associated with more common PDEs, such as the Laplace equation and the Stokes equation. There is an FMM implementation for the modified Helmholtz equation described in \cite{Kropinski2011modHelm}, but to the best of our knowledge there are no Ewald summation formulas derived, upon which a spectral Ewald method can be built. In this paper we present a derivation of such formulas, along with truncation error estimates, and a spectral Ewald summation method for an $\mathcal{O}(N\log N)$ evaluation for the sums \eqref{eq:mhhbiedisc}, \eqref{eq:mhhslpdisc} and \eqref{eq:mhhslpgraddisc}, and their periodically extended counterparts. For the interest of space, the main results are presented in this paper, while supplementary material is found in \cite{PalssonSara2019SaEs}.

We consider fast Ewald summation over an FMM since periodicity is simple to incorporate and considerable speedups are possible for uniformly distributed target points, such as the uniform grid, on which we evaluate the solution used in our computations. The basic concept of Ewald summation is as follows. A sum is decomposed into a sum in ``real space''  and a sum in ``$\fovar$  space''; the terms in the former decays rapidly in $\mathbb{R}^{2}$ and the Fourier transform of the latter converges rapidly in the frequency domain. Hence, the sums can be truncated at a cut-off radius and a maximum frequency, respectively. The evaluation of the Fourier part is then accelerated using FFTs, as will be explained below.
%Here we present a summation of the main results in \cite{PalssonSara2019SaEs}. The novelties in this paper are: the Ewald decompositon, error estimates and improvements in terms of computational efficiency when the target points are located on a uniform grid.

\subsubsection{Ewald decomposition}
\label{sss:ewalddecomp}
We present the basic Ewald decomposition for $\mhhu$ and  $\nabla \mhhu$ in free-space for a target points $\x$ in $\domain$. The derivations are presented in appendix \ref{s:appsupewald}. Initially we consider $\mhhu$, the principle for $\nabla\mhhu$ is analogous. Associate each pair of integers $(i,j)$, $i = 1,\ldots,\Npan$ and $j = 1,\ldots,\Nq$ with some $n = 1,\ldots,\Nbdry$, such that $\bimdens_{n} = \bimdens_{ij}$, $s_{n} = s_{ij}$ etc., and write \eqref{eq:mhhslpdisc} as

\begin{equation}
  \label{eq:mhhslpdiscewald}
  \mhhu\fp{\x} = \frac{1}{2\pi}\sum\limits_{n = 1}^{\Nbdry}\kerdomain\fp{\y_{n},\x}\ewfunn, \quad \x\in\domain,
\end{equation}
where 
\begin{equation}
  \label{eq:ewalddata}
  \ewfunn = \bimdens_{n}s_{n}\wG_{n}.
\end{equation}
By using the Ewald screening function
\begin{equation}
  \label{eq:ewaldscreening}
  \ewscreen(\x,\ewparam) = \frac{\ewparam^{2}}{\pi}e^{-\mhaa/4\ewparam^{2}}e^{-\ewparam^{2}\norm{\x}}\Leftrightarrow \widehat{\ewscreen}(\fovarvec,\ewparam) = e^{-({\mhaa}+\norm{\fovarvec}^{2})/4\ewparam^{2}},
\end{equation}
the kernel $\kerdomain$ can be written as
\begin{equation}
  \label{eq:ewalddecompkernel}
  \kerdomain\fp{\y_{n},\x} = \underbrace{\kerdomain\fp{\y_{n},\x} - \fp{\kerdomain\fp{\y_{n},\cdot}\ast\ewscreen\fp{\cdot,\ewparam}}\fp{\x}}_{\ewrealker\fp{\y_{n}-\x,\ewparam}}
  + \underbrace{\fp{\kerdomain\fp{\cdot,\y_{n}}\ast\ewscreen\fp{\cdot,\ewparam}}\fp{\x}}_{\ewfourierker\fp{\y_{n}-\x,\ewparam}}.
\end{equation}
With this decomposition, an alternative formulation of the homogeneous solution $\mhhu$ \eqref{eq:mhhslpdiscewald} is
\begin{equation}
    \label{eq:ewalddecomp}
    \mhhu \fp{\x} = \frac{1}{2\pi}\bigg(\underbrace{\sum\limits_{n = 1}^{\Nbdry}\ewrealker(\y_{n}-\x,\ewparam)\ewfunn}_{= \ewur\fp{\x}} + \underbrace{\frac{1}{4\pi^{2}}\int\limits_{\mathbb{R}^{2}}\widehat{\ewfourierker}(\fovarvec,\ewparam)\sum\limits_{n = 1}^{\Nbdry}\ewfunn e^{i\fovarvec\cdot(\x-\y_{n})}\D \fovarvec}_{=\ewuf\fp{\x}}\bigg),
\end{equation}
with the closed form expressions
\begin{align}
  \label{eq:ewalddecompreal}
\ewrealker(\y-\x,\ewparam) &= \frac{1}{2}\besselK{0}\fp{\norm{\y-\x}^{2}\ewparam^{2},\frac{\mhaa}{4\ewparam}},\\
  \label{eq:ewalddecompfourier}
\widehat{\ewfourierker}(\fovarvec,\ewparam) &= \frac{2\pi}{\mhaa+\norm{\fovarvec}^{2}}e^{-(\mhaa + \norm{\fovarvec}^{2})/4\ewparam^{2}}.
\end{align}
Here $\besselK{0}\fp{\cdot,\cdot}$ denotes the incomplete modified Bessel function of the second kind of zeroth order \cite{HarrisFrankE2008IBgi}, defined as 
\begin{equation}
\label{eq:incompbesselK}
\besselK{\nu}\fp{\rho_1,\rho_2} = \int\limits_{1}^{\infty}\frac{e^{-\rho_{1}t-\rho_{2}/t}}{t^{\nu+1}}\,\D t,\quad \nu\in\mathbb{Z}.
\end{equation}
If a target point and a source point coincide, i.e. $\x = \y_{n}$, then the self-interaction term
\begin{equation}
  \label{eq:ewalddecompself}
\ewselfker\fp{\ewparam} = -\frac{1}{4\pi}\int\limits_{1}^{\infty}\frac{e^{-\frac{\mhaa}{4\ewparam^{2}}t}}{t}\D t = \frac{-1}{4\pi} E_{1}\left(\frac{\mhaa}{4\ewparam^{2}}\right)
\end{equation}
is added to \eqref{eq:ewalddecomp}.

Analogously, the decomposition $\nabla\mhhu$ from \eqref{eq:mhhslpgraddisc} is
\begin{equation}
    \label{eq:ewalddecompgrad}
\nabla    \mhhu \fp{\x} = \frac{1}{2\pi}\bigg(\underbrace{\sum\limits_{n = 1}^{\Nbdry}\ewrealkergrad(\y_{n}-\x,\ewparam)\ewfungradn}_{= \nabla\ewur\fp{\x}} + \underbrace{\frac{1}{4\pi^{2}}\int\limits_{\mathbb{R}^{2}}\widehat{\ewfourierkergrad}(\fovarvec,\ewparam)\sum\limits_{n = 1}^{\Nbdry}\ewfungradn e^{i\fovarvec\cdot(\y_{n}-\x)}\D\fovarvec}_{= \nabla\ewuf\fp{\x}}\bigg),
\end{equation}
%\begin{equation}
%    \label{eq:ewalddecompgrad}
%\nabla    \mhhu \fp{\x} = \frac{\mhaa}{2\pi}\bigg(\underbrace{\sum\limits_{n = 1}^{\Nbdry}\ewrealkergrad(\x-\y_{n},\ewparam)\ewfungrad\fp{\y_{n}}}_{= \nabla\ewur\fp{\x}} + \underbrace{\frac{1}{4\pi^{2}}\int\limits_{\mathbb{R}^{2}}\widehat{\ewfourierkergrad}(\fovarvec,\ewparam)\sum\limits_{n = 1}^{\Nbdry}\ewfungrad\fp{\y_{n}}e^{i\fovarvec\cdot(\y_{n}-\x)}\D \fovarvec}_{=\nabla\ewuf\fp{\x}}\bigg),
%\end{equation}
with $\ewfungradn= \ewfunn$ for $n = 1,\ldots,\Nbdry$, and
\begin{align}
  \label{eq:ewalddecomprealgrad}
\ewrealkergrad(\y-\x,\ewparam) &= \ewparam^{2}(\y-\x)\besselK{-1}\fp{\norm{\y-\x}^{2}\ewparam^{2},\frac{\mhaa}{4\ewparam^{2}}},\\
  \label{eq:ewalddecompfouriergrad}
\widehat{\ewfourierkergrad}(\fovarvec,\ewparam) &=\frac{-i2\pi\fovarvec}{\mhaa+\norm{\fovarvec}^{2}}e^{-(\mhaa + \norm{\fovarvec}^{2})/4\ewparam^{2}},
\end{align}
with $\besselK{-1}\fp{\cdot,\cdot}$ from \eqref{eq:incompbesselK}. The self-interaction term $\ewselfkergrad$ is identically equal to zero.

The sum in \eqref{eq:mhhbiedisc} can be treated analogously to $\nabla\mhhu$, since
\begin{equation}
\sum\limits_{n = 1}^{\Nbdry}\kerbdry\fp{\y_{n},\x_{n^\prime}}\bimdens_{n}s_{n}\wG_{n} = \sum\limits_{n = 1}^{\Nbdry}\nabla\kerdomain\fp{\y_{n},\x_{n^\prime}} f_{\mathcal{D},n}\cdot\normal\fp{\x_{n^\prime}} \quad n^\prime = 1,\ldots, \Nbdry,
\end{equation}
with $f_{\mathcal{D},n} = \bimdens_{n}s_{n}\wG_{n}$.

For the periodic case the sums over  $\ewrealker$ and $\ewrealkergrad$ also include periodic images, and the inverse Fourier transforms are replaced by inverse Fourier series, see \cite{PalssonSara2019SaEs}.

\subsubsection{The spectral Ewald method}
\label{sss:ewalddecomp}
The spectral Ewald method is a fast summation method that makes the real space sum cost $\mathcal{O}(N)$ and the $\fovar$ space sum cost $\mathcal{O}(N\log N)$ for $N$ sources and targets \cite{AfKlintebergLudvig2014FEsf,LindboDag2010Safs,LINDBO20118744}. The balance, in terms of computational effort, between the two is controlled by the splitting parameter $\ewparam$.

%If the target points are on a uniform grid, then the final step in the spectral Ewald method is superfluous and can be skipped. This is the case for \eqref{eq:mhhslpdisc}  and \eqref{eq:mhhslpgraddisc} when evaluating on the uniform grid $\X$, but not for  \eqref{eq:mhhbiedisc} which is evaluated at the $\Nbdry$ boundary nodes, see the end of this section.

The summation in real space is accelerated as follows. Introduce a cut-off radius $\ewcutoffrad>0$, and include only contributions from source points within a distance $\ewcutoffrad$ of the target point $\x$ in the real space sum. That is
\begin{equation}
\label{eq:ewurtrunc}
\ewur\fp{\x} \approx \sum\limits_{\{\y_{n}:\,\, \norm{\y_{n}-\x} < \ewcutoffrad\}}\ewrealker(\,\norm{\y_{n}-\x},\ewparam)\ewfunn,
\end{equation}
with $\ewrealker$ as in \eqref{eq:ewalddecompreal}. The complexity can now be reduced to $\mathcal{O}(N)$ by modifying $\ewcutoffrad$ with increasing $N$, such that the number of elements in $\{\y_{n}:\,\, \norm{\y_{n}-\x} < \ewcutoffrad\}$ remains constant. By first creating a linked list, the creation of $\{\y_{n}:\,\, \norm{\y_{n}-\x} < \ewcutoffrad\}$ also costs $\mathcal{O}(N)$. The summation of $\nabla\ewur$ is treated analogously.

We now study the $\fovar$ space sum $\ewuf$. The goal is to compute it for any target point in the plane, discretizing the involved integral with the trapezoidal rule and accelerate the computations with FFTs. If $\mha$ is small, which is not the case in this paper, then we need to take the approach discussed in appendix \ref{s:appsupewald}.

 The idea is now to spread the data at the source points $\y_{n}$ to a uniform grid of $\ewN$ points with spacing $h$ on the box $[-\mathcal{R}/2,\mathcal{R}/2]^{2}$, apply an FFT, scale the results with the Green's function and a window function, apply an inverse FFT and evaluate the integral in the $\fovar$ space sum at the target point $\x$ with the trapezoidal rule. The spreading means convolving the data $\ewfun$ with a window function $\ewwind(\x,\ewparam,\ewwindparam)$,
\begin{equation}
  \label{eq:convwind}
  \ewH (\x) = \sum\limits_{n=1}^{\Nbdry}\ewfunn  \ewwind\fp{\y_{n}-\x,\ewparam,\ewwindparam},
\end{equation}
 where $\ewwindparam > 0$ is a scaling parameter. The Fourier transform of the window function is denoted $\widehat{\ewwind}\fp{\,\norm{\fovarvec},\ewparam,\ewwindparam} := \widehat{\ewwindk}$. We use
\begin{equation}
  \label{eq:ewwind}
  \ewwind\fp{\x,\ewparam,\ewwindparam} =
\begin{cases}
  e^{-4\ewwindparam\,\abs{x_{1}}/p^{2}h^{2}}e^{-4\ewwindparam\,\abs{x_{2}}/p^{2}h^{2}},\quad &\max\fp{\,\abs{x_{1}},\abs{x_{2}}}\leq ph/2,\\
  0,\quad &\text{otherwise},
\end{cases}
\end{equation}
for $\x = (x_{1},x_{2})$. Here $p$ denotes the number of grid points in the truncation, such that the domain of support is a square with sides of length $ph$. The parameter $\ewwindparam$ balances the errors from truncating the window function and how well the Gaussians \eqref{eq:ewwind} are approximated on the grid. We use $\eta=0.95^2\pi p/2$.

In \cite{PalssonSara2019SaEs} it is shown that if
\begin{equation}
\widehat{\tilde{\ewH}}\fp{\fovarvec} = \widehat{\ewfourierkertrun}(\fovarvec,\ewparam)\frac{\widehat{\ewH}\fp{-\fovarvec}}{\widehat{\ewwindk}^{2}},
\end{equation}
where
\begin{equation}
  \widehat{\ewH}(\fovarvec) = \sum\limits_{n = 1}^{\Nbdry}\ewfunn \widehat{\ewwindk}e^{-i\fovarvec\cdot\y_{n}},
\end{equation}
then the k space sum approximation \eqref{eq:ewaldsumfourierapprox} can be written as
\begin{equation}
\label{eq:ewaldsumfourierapprox2}
\ewuf\fp{\x} \approx \int\limits_{\mathbb{R}^{2}}\ewwindk\fp{\x-\bm{\rho}}\tilde{\ewH}\fp{\bm{\rho}}\D \bm{\rho},
\end{equation}
where $\tilde{\ewH}$ is the inverse Fourier transform of $\widehat{\tilde{\ewH}}$. This integral is computed with the trapezoidal rule.

To summarize, the steps for computing the spectral Ewald approximation of the $\fovar$ space sum for $u$ are
\begin{itemize}
   \item \emph{Spreading}: evaluate $\ewH\fp{\x}$ in \eqref{eq:convwind} on the uniform grid with $\ewN$ points. This involves evaluating $\mathcal{O}(\Nbdry)$ window functions on $p^{2}$ points each, at a cost of $\mathcal{O}(p^{2}\Nbdry)$.
   \item \emph{FFT}: compute $\widehat{\ewH}$ using a $2$D FFT, at a cost $\mathcal{O}(\ewN\log\ewN)$.
   \item \emph{Scaling}: compute the tensor product $\widehat{\tilde{\ewH}}(\fovarvec)$ for all $\fovarvec\in [-\sqrt{\ewN}/2,\,\sqrt{\ewN}/2-1]^{2}$, which is of $\mathcal{O}(\ewN)$.
   \item \emph{IFFT}: compute $\tilde{\ewH}\fp{\x}$ on the uniform grid using the $2$D inverse FFT, at a cost of $\mathcal{O}(\ewN\log\ewN)$.
   \item \emph{Quadrature}: The integral \eqref{eq:ewaldsumfourierapprox2} needs to be computed for each target point $\x$ to obtain $\ewuf\fp{\x}$. This is done with the trapezoidal rule. Using the compactly supported window function \eqref{eq:ewwind}, this is of cost $\mathcal{O}\fp{p^{2}}$ for each target point.
\end{itemize}
The steps for evaluating spectral Ewald approximation of the $\fovar$ space sum for $\nabla u$ are analogous, except now three inverse FFTs are required, since this is a vector quantity.

To further reduce the cost of the spreading and quadrature steps fast Gaussian gridding can be applied \cite{LeslieGreengard2004AtNF}. See \cite{LINDBO20118744} for implementation details.

For the sums \eqref{eq:mhhslpdisc} and \eqref{eq:mhhslpgraddisc} the target points $\x$ belong to the uniform grid $\X$, and the final quadrature step can be avoided. It is achieved by constructing the uniform grid in the spreading step such that $\X$ is a subset of it. See \cite{PalssonSara2019SaEs} for how to modify the scheme accordingly.

To understand the computational complexity of the full scheme, we have to consider also the evaluation of the real space sum. We will comment on this in the parameter selection section \ref{sss:ewaldestK1} below.

\subsubsection{Error estimates}
\label{sss:ewalderrest}
The basic idea of Ewald decomposition is to split a slowly converging sum into two rapidly converging sums in their respective spaces. The resulting sums are then truncated based on carefully selected rules such that a certain tolerance is satisfied with optimal efficiency. These rules are presented here. For their derivations, see \cite{PalssonSara2019SaEs}.

Truncation errors for the real space sums $\ewur$ \eqref{eq:ewurtrunc}, and the corresponding for $\nabla\ewur$, are introduced when they are limited to source points within a distance $\ewcutoffrad$ of the target point. Thus these errors are essentially the same in both the free-space and periodic setting. Let $\delta\ewur$ and $\delta\nabla\ewur$ be the RMS errors. They can be approximated as
\begin{equation}
  \label{eq:realestk0}
  (\delta\ewur)^2 \approx \frac{\pi\ewfunsum}{4L^{2}\ewcutoffrad^{4}\ewparam^{6}}e^{-2\ewcutoffrad^{2}\ewparam^{2}},
\end{equation}
and
\begin{equation}
    \label{eq:realestk1}
  (\delta\nabla\ewur)^2 \approx \frac{\pi\ewfungradsum}{L^{2}\mhaa\ewcutoffrad^{2}\ewparam^{2}}e^{-2\ewparam^{2}\ewparam^{2}},
\end{equation}
where $\ewfunsum = \sum_{n}\ewfunn^{2}$ and $\ewfungradsum = \sum_{n}\ewfungradn^{2}$. How well the error estimates perform, for a range of values for $\ewparam$ and $\ewcutoffrad$, is demonstrated in figure \ref{fig:ewaldestK0} and in figure \ref{fig:ewaldestK1}. They predict the errors well for error levels below $10^{-2}$.
% see \cite{PalssonSara2019SaEs} and the references therein for guidance on setting the parameters.

For the k space sums $\ewuf$ and $\nabla\ewuf$ truncation errors are introduced when wave numbers $\fovarvec$ that satisfy $\norm{\fovarvec}\geq\fovarlim$ are excluded. The error estimates for the corresponding errors $\delta\ewuf$ and $\delta\nabla\ewuf$ are different for the free-space and periodic setting. We present here only the former, the latter is given in \cite{PalssonSara2019SaEs}. The RMS errors can be approximated by
\begin{equation}
  \label{eq:fourirestk0}
  (\delta\ewuf)^2 \approx \frac{64\pi\ewfunsum\ewparam^{4}}{L(\mhaa + \fovarlim^{2})^{2}}e^{-2(\mhaa + \fovarlim^{2})/4\ewparam}\left(\frac{1}{\sqrt{2\pi\fovarlim}} - \frac{\mha\besselK{0}\fp{\mha\mathcal{R}}}{\sqrt{\mathcal{R}\pi}} -\frac{\mha\sqrt{\mathcal{R}}\besselK{1}\fp{\mha\mathcal{R}}}{\pi\fovarlim} \right)^{2},
\end{equation}
and
\begin{equation}
    \label{eq:fourirestk1}
  (\delta\nabla\ewuf)^2 \approx \frac{8\ewfungradsum\ewparam^{2}}{L\pi^{2}\mhaa}\frac{e^{-2(\mhaa + \fovarlim^{2})/4\ewparam}}{(\mhaa + \fovarlim^{2})^{2}}\left(\sqrt{2\pi\fovarlim} - \frac{8\mha\besselK{0}\fp{\mha\mathcal{R}}\fovarlim}{\sqrt{\mathcal{R}}} -2\mha\sqrt{\mathcal{R}}\besselK{1}\fp{\mha\mathcal{R}} \right)^{2}.
\end{equation}

The truncation errors and estimates for $\mhhu$ are shown together in figure  \ref{fig:ewaldestK0}, for a test domain with $500$ randomly distributed sources and targets, together with random point sources $\ewfunn\in [0,1]$. In this example, $\mha = 1$ and $L = 2\pi$. The estimates hold for varying $L$ and $\alpha$. The estimates follow the error well for error levels below $10^{-2}$. For larger errors, the estimate is not as sharp. That is, however, far from the region of interest for an accurate computation. The corresponding plots for $\nabla\mhhu$ are shown in figure \ref{fig:ewaldestK1}, where we get similar results.
\begin{figure}[htbp]
% Generated by: est_plots.jl
  \centering
  \includestandalone[width=0.49\textwidth]{fredrik_real_K0_est}
  \hfill
  \includestandalone[width=0.49\textwidth]{fredrik_kspace_K0_est}
  \caption{Truncation error (absolute) and estimates for different values of $\xi$ for $\ewur$ and $\ewuf$. The test domain contains $500$ randomly distributed sources and targets, and random point sources $\ewfunn\in [0,1]$. Left: estimate as derived in \eqref{eq:realestk0} for the real space sum for different cut-off radii $r_{c}$. Right: estimate from \eqref{eq:fourirestk0} for the k space sum, when varying $k_{\infty}$.}
  \label{fig:ewaldestK0}
\end{figure}

\begin{figure}[htbp]
% Generated by: est_plots.jl
  \centering
  \includestandalone[width=0.49\textwidth]{fredrik_real_K1_est}
  \hfill
  \includestandalone[width=0.49\textwidth]{fredrik_kspace_K1_est}
  \caption{Truncation error (absolute) and estimates for different values of $\xi$ for \linebreak $\nabla\ewur$ and $\nabla\ewuf$. The test domain contains $500$ randomly distributed sources and targets, and random point sources $f_{\mathcal{D},n}\in [0,1]$. Left: estimate as derived in \eqref{eq:realestk1} for the real space sum for different cut-off radii $r_{c}$. Right: estimate from \eqref{eq:fourirestk1} for the k space sum, when varying $k_{\infty}$.}
  \label{fig:ewaldestK1}
\end{figure}

\subsubsection{Parameter selection}
\label{sss:ewaldestK1}
The parameters of the Ewald summation impact the split of computational work between the real space sum and the k space sum and the level of the truncation errors. To keep the evaluation of the real space sum $\mathcal{O}(N)$, for $N$ sources/targets, the cut-off radius $r_c$ is set initially, such that the average number of neighbors of a target point is constant as $N$ is increased. Using the truncation error estimate in \eqref{eq:realestk0} and \eqref{eq:realestk1} respectively, a suitable $\xi$ is chosen given the desired truncation error level.  Given $\xi$, an appropriate $k_\infty=\ewN/2$ is computed from the truncation error estimates for the k space sum \eqref{eq:fourirestk0}, \eqref{eq:fourirestk1}. This couples the grid size for the k space evaluations to $N$, the number of sources/targets, and yields a total computational complexity of the method of $\mathcal{O}(N\log N)$. The support of the window functions, $p$, is chosen large enough to keep the approximation errors low. Here, $p=24$ suffices. The parameter $\eta$ is chosen in order to balance the errors from truncating the Gaussians, and resolving them on the grid.%on and approximation, in other words how well the Gaussians are resolved on the grid. Here, we have chosen $\eta=0.95^2\pi p/2$.
%%%%%%%%%%%%%%%%%%%%%%%End of file %%%%%%%%%%%%%%%%%%%%%%%%%%%%%%

%%% Local Variables:
%%% mode: latex
%%% TeX-master: "manuscript.tex"
%%% End:

%% file: fredrik_real_K0_est.tex

\definecolor{mycolor2}{rgb}{0.85000,0.32500,0.09800}%

\begin{tikzpicture}

\begin{axis}[%
  width=4.0187in,
  height=3.1693in,
  at={(1.011in,0.683in)},
scale only axis,
xmin=0,
xmax=1.6,
xtick distance=0.4,
xlabel style={font=\color{white!15!black}},
xlabel={\Large $r_c$},
ymode=log,
ymin=1e-16,
ymax=50,
ytick distance=1e4,
yminorticks=true,
axis background/.style={fill=white},
axis x line*=bottom,
axis y line*=left,
ylabel style={font=\color{white!15!black}},
ylabel={\Large Error, $\ewur$},
xmajorgrids,
ymajorgrids,
yminorgrids,
legend style={legend cell align=left, align=left, draw=white!15!black},
execute at begin axis={
    \draw (rel axis cs:0,0) -- (rel axis cs:1,0)
          (rel axis cs:0,1) -- (rel axis cs:1,1);
}
]
\addplot [color=black, line width=2.0pt]
  table[row sep=crcr]{%
0.01	0.699184931674334\\
0.0418529862611203	0.638467913690762\\
0.0737059725222407	0.534482402861655\\
0.105558958783361	0.429736504303604\\
0.137411945044481	0.352726154732132\\
0.169264931305602	0.282139089683418\\
0.201117917566722	0.231953619310251\\
0.232970903827842	0.181371473384394\\
0.264823890088963	0.142809571209477\\
0.296676876350083	0.109598596927031\\
0.328529862611203	0.0858866394711645\\
0.360382848872324	0.058739298366476\\
0.392235835133444	0.0438276184528128\\
0.424088821394564	0.0329912144952847\\
0.455941807655685	0.0217768032720603\\
0.487794793916805	0.0149459487211414\\
0.519647780177925	0.0101641899271533\\
0.551500766439046	0.00659860357026272\\
0.583353752700166	0.00443141119433352\\
0.615206738961286	0.00281330326356445\\
0.647059725222407	0.00184133916883955\\
0.678912711483527	0.00114260502413556\\
0.710765697744647	0.00075902787589625\\
0.742618684005768	0.000469919908657813\\
0.774471670266888	0.000280484094321824\\
0.806324656528008	0.000168235501375859\\
0.838177642789129	9.81122356000742e-05\\
0.870030629050249	5.43893508361195e-05\\
0.901883615311369	2.91025856556736e-05\\
0.93373660157249	1.69502120678917e-05\\
0.96558958783361	9.48454427380166e-06\\
0.997442574094731	5.20129793065709e-06\\
1.02929556035585	2.7910300503632e-06\\
1.06114854661697	1.4399586046179e-06\\
1.09300153287809	7.62031545645318e-07\\
1.12485451913921	3.63649369758254e-07\\
1.15670750540033	1.79360798323667e-07\\
1.18856049166145	9.08283865645464e-08\\
1.22041347792257	4.29553113416806e-08\\
1.25226646418369	2.10329662405543e-08\\
1.28411945044481	9.65350023205675e-09\\
1.31597243670593	4.42164068669301e-09\\
1.34782542296705	1.89292242972623e-09\\
1.37967840922817	8.00768976333651e-10\\
1.41153139548929	3.55074762415991e-10\\
1.44338438175042	1.48559147172973e-10\\
1.47523736801154	6.39421284157744e-11\\
1.50709035427266	2.59105586066568e-11\\
1.53894334053378	1.00618681414134e-11\\
1.5707963267949	4.01444872293127e-12\\
};
\addlegendentry{\large error}

\addplot [color=mycolor2, dashed, line width=2.0pt]
  table[row sep=crcr]{%
0.01	680.306235430977\\
0.0418529862611203	38.2644732948679\\
0.0737059725222407	11.9359210276676\\
0.105558958783361	5.5277988001407\\
0.137411945044481	3.04258724428055\\
0.169264931305602	1.83643725404045\\
0.201117917566722	1.16975543482617\\
0.232970903827842	0.769747385874788\\
0.264823890088963	0.516488588495305\\
0.296676876350083	0.3503479526533\\
0.328529862611203	0.238823836127074\\
0.360382848872324	0.162902596909362\\
0.392235835133444	0.110830432312994\\
0.424088821394564	0.075024949384061\\
0.455941807655685	0.0504353212015486\\
0.487794793916805	0.0336187865066656\\
0.519647780177925	0.0221926620615247\\
0.551500766439046	0.0144935568092197\\
0.583353752700166	0.0093564555325505\\
0.615206738961286	0.00596636448877697\\
0.647059725222407	0.00375585033426669\\
0.678912711483527	0.00233282773079352\\
0.710765697744647	0.00142902435367276\\
0.742618684005768	0.000863000952762709\\
0.774471670266888	0.000513630590922009\\
0.806324656528008	0.00030118232109832\\
0.838177642789129	0.000173953463905171\\
0.870030629050249	9.89373967382305e-05\\
0.901883615311369	5.54014807962692e-05\\
0.93373660157249	3.05375978961607e-05\\
0.96558958783361	1.65664032606488e-05\\
0.997442574094731	8.84373620417004e-06\\
1.02929556035585	4.64513196224266e-06\\
1.06114854661697	2.4002751463007e-06\\
1.09300153287809	1.22004390772351e-06\\
1.12485451913921	6.09952986784787e-07\\
1.15670750540033	2.99904175137132e-07\\
1.18856049166145	1.4500936411734e-07\\
1.22041347792257	6.89448949793537e-08\\
1.25226646418369	3.22306117787824e-08\\
1.28411945044481	1.48137705999399e-08\\
1.31597243670593	6.69369258497106e-09\\
1.34782542296705	2.97333213132777e-09\\
1.37967840922817	1.29829931991876e-09\\
1.41153139548929	5.57234323050297e-10\\
1.44338438175042	2.35077965667588e-10\\
1.47523736801154	9.74715585820957e-11\\
1.50709035427266	3.97208248340332e-11\\
1.53894334053378	1.59080651040721e-11\\
1.5707963267949	6.26127562131296e-12\\
};
\addlegendentry{\large estimate}

\addplot [color=black, line width=2.0pt, forget plot]
  table[row sep=crcr,forget plot]{%
0.01	0.409171441350409\\
0.0418529862611203	0.334059864640036\\
0.0737059725222407	0.238712425193285\\
0.105558958783361	0.160850907569177\\
0.137411945044481	0.106314079583839\\
0.169264931305602	0.068439576389602\\
0.201117917566722	0.043927546537931\\
0.232970903827842	0.0265101959721853\\
0.264823890088963	0.0157439656243007\\
0.296676876350083	0.00821327336912609\\
0.328529862611203	0.00475399677715057\\
0.360382848872324	0.00213672436197735\\
0.392235835133444	0.0010601436415521\\
0.424088821394564	0.000539897101370534\\
0.455941807655685	0.000226051616218944\\
0.487794793916805	9.27615525574086e-05\\
0.519647780177925	3.64752603528223e-05\\
0.551500766439046	1.41550917767399e-05\\
0.583353752700166	5.10695113007247e-06\\
0.615206738961286	1.75568659455004e-06\\
0.647059725222407	6.16963552602457e-07\\
0.678912711483527	1.84903132966127e-07\\
0.710765697744647	6.01386663420241e-08\\
0.742618684005768	1.82542140534175e-08\\
0.774471670266888	4.94779427672836e-09\\
0.806324656528008	1.30811634688978e-09\\
0.838177642789129	3.24533936660868e-10\\
0.870030629050249	7.65437025865566e-11\\
0.901883615311369	1.54408390953157e-11\\
0.93373660157249	3.61916038004314e-12\\
0.96558958783361	7.90248826522441e-13\\
0.997442574094731	1.54020419627745e-13\\
1.02929556035585	2.91909402323653e-14\\
1.06114854661697	5.23243674098856e-15\\
1.09300153287809	1.15808605145789e-15\\
1.12485451913921	4.28215888875103e-16\\
1.15670750540033	3.7558280817917e-16\\
1.18856049166145	3.72948263913496e-16\\
1.22041347792257	3.72944410143515e-16\\
1.25226646418369	3.72943135087342e-16\\
1.28411945044481	3.72943135087342e-16\\
1.31597243670593	3.72943135087342e-16\\
1.34782542296705	3.72943135087342e-16\\
1.37967840922817	3.72943135087342e-16\\
1.41153139548929	3.72943135087342e-16\\
1.44338438175042	3.72943135087342e-16\\
1.47523736801154	3.72943135087342e-16\\
1.50709035427266	3.72943135087342e-16\\
1.53894334053378	3.72943135087342e-16\\
1.5707963267949	3.72943135087342e-16\\
};

\addplot [color=mycolor2, dashed, line width=2.0pt, forget plot]
  table[row sep=crcr]{%
0.01	146.711221008919\\
0.0418529862611203	8.03669753821256\\
0.0737059725222407	2.3635256639878\\
0.105558958783361	0.999030188104695\\
0.137411945044481	0.485837373800759\\
0.169264931305602	0.250810069081337\\
0.201117917566722	0.132277030109939\\
0.232970903827842	0.0697682089568664\\
0.264823890088963	0.0363236570700078\\
0.296676876350083	0.0185074848437449\\
0.328529862611203	0.00917368804622615\\
0.360382848872324	0.00440465890037228\\
0.392235835133444	0.00204202368656144\\
0.424088821394564	0.000911852157972332\\
0.455941807655685	0.000391444304423554\\
0.487794793916805	0.00016129956125455\\
0.519647780177925	6.37200364487229e-05\\
0.551500766439046	2.41077429771176e-05\\
0.583353752700166	8.72785234880306e-06\\
0.615206738961286	3.02148030858951e-06\\
0.647059725222407	9.99612588336435e-07\\
0.678912711483527	3.15877765416425e-07\\
0.710765697744647	9.52990266250379e-08\\
0.742618684005768	2.74392385687394e-08\\
0.774471670266888	7.53742687695434e-09\\
0.806324656528008	1.97475253513243e-09\\
0.838177642789129	4.93317891924841e-10\\
0.870030629050249	1.17480239255201e-10\\
0.901883615311369	2.66647761057159e-11\\
0.93373660157249	5.76731233643297e-12\\
0.96558958783361	1.18861079344806e-12\\
0.997442574094731	2.33500162935274e-13\\
1.02929556035585	4.38378301185962e-14\\
1.06114854661697	7.98833429396485e-15\\
1.09300153287809	1.53888891411918e-15\\
1.12485451913921	4.34643617417901e-16\\
1.15670750540033	2.54722404117067e-16\\
1.18856049166145	2.26826182875346e-16\\
1.22041347792257	2.22710616393243e-16\\
1.25226646418369	2.22132903163493e-16\\
1.28411945044481	2.22055746651654e-16\\
1.31597243670593	2.22045942923239e-16\\
1.34782542296705	2.22044757834918e-16\\
1.37967840922817	2.22044621554134e-16\\
1.41153139548929	2.22044606645842e-16\\
1.44338438175042	2.22044605094468e-16\\
1.47523736801154	2.22044604940905e-16\\
1.50709035427266	2.22044604926446e-16\\
1.53894334053378	2.22044604925151e-16\\
1.5707963267949	2.22044604925041e-16\\
};

\addplot [color=black, line width=2.0pt, forget plot]
  table[row sep=crcr,forget plot]{%
0.01	0.187553455161797\\
0.0418529862611203	0.110961722739225\\
0.0737059725222407	0.0482068741457576\\
0.105558958783361	0.0194514935980138\\
0.137411945044481	0.00577192277287467\\
0.169264931305602	0.00160963670734627\\
0.201117917566722	0.000409747899280907\\
0.232970903827842	8.45757469744773e-05\\
0.264823890088963	1.55069027625681e-05\\
0.296676876350083	1.71741071309136e-06\\
0.328529862611203	2.24634546848288e-07\\
0.360382848872324	1.88154241866397e-08\\
0.392235835133444	1.52740417323858e-09\\
0.424088821394564	1.09969606712736e-10\\
0.455941807655685	5.83894052678065e-12\\
0.487794793916805	2.59580214258231e-13\\
0.519647780177925	8.41402654131882e-15\\
0.551500766439046	5.2422995841482e-16\\
0.583353752700166	3.02905946955615e-16\\
0.615206738961286	2.97894364988146e-16\\
0.647059725222407	2.97894206136889e-16\\
0.678912711483527	2.97894206136889e-16\\
0.710765697744647	2.97894206136889e-16\\
0.742618684005768	2.97894206136889e-16\\
0.774471670266888	2.97894206136889e-16\\
0.806324656528008	2.97894206136889e-16\\
0.838177642789129	2.97894206136889e-16\\
0.870030629050249	2.97894206136889e-16\\
0.901883615311369	2.97894206136889e-16\\
0.93373660157249	2.97894206136889e-16\\
0.96558958783361	2.97894206136889e-16\\
0.997442574094731	2.97894206136889e-16\\
1.02929556035585	2.97894206136889e-16\\
1.06114854661697	2.97894206136889e-16\\
1.09300153287809	2.97894206136889e-16\\
1.12485451913921	2.97894206136889e-16\\
1.15670750540033	2.97894206136889e-16\\
1.18856049166145	2.97894206136889e-16\\
1.22041347792257	2.97894206136889e-16\\
1.25226646418369	2.97894206136889e-16\\
1.28411945044481	2.97894206136889e-16\\
1.31597243670593	2.97894206136889e-16\\
1.34782542296705	2.97894206136889e-16\\
1.37967840922817	2.97894206136889e-16\\
1.41153139548929	2.97894206136889e-16\\
1.44338438175042	2.97894206136889e-16\\
1.47523736801154	2.97894206136889e-16\\
1.50709035427266	2.97894206136889e-16\\
1.53894334053378	2.97894206136889e-16\\
1.5707963267949	2.97894206136889e-16\\
};

\addplot [color=mycolor2, dashed, line width=2.0pt, forget plot]
  table[row sep=crcr]{%
0.01	18.2018753510154\\
0.0418529862611203	0.880910954604218\\
0.0737059725222407	0.196571333944968\\
0.105558958783361	0.0541436890732535\\
0.137411945044481	0.0147357896720995\\
0.169264931305602	0.0036563235256801\\
0.201117917566722	0.000795985457450841\\
0.232970903827842	0.000148834363311618\\
0.264823890088963	2.35919610160828e-05\\
0.296676876350083	3.14307399936159e-06\\
0.328529862611203	3.49854765044135e-07\\
0.360382848872324	3.23963905164664e-08\\
0.392235835133444	2.48765163675777e-09\\
0.424088821394564	1.5801637368018e-10\\
0.455941807655685	8.28725801104046e-12\\
0.487794793916805	3.58498706812949e-13\\
0.519647780177925	1.29753080334282e-14\\
0.551500766439046	5.95437235452454e-16\\
0.583353752700166	2.31028907053189e-16\\
0.615206738961286	2.2222213354907e-16\\
0.647059725222407	2.22047484019332e-16\\
0.678912711483527	2.22044643227139e-16\\
0.710765697744647	2.22044605342838e-16\\
0.742618684005768	2.22044604928767e-16\\
0.774471670266888	2.22044604925059e-16\\
0.806324656528008	2.22044604925031e-16\\
0.838177642789129	2.22044604925031e-16\\
0.870030629050249	2.22044604925031e-16\\
0.901883615311369	2.22044604925031e-16\\
0.93373660157249	2.22044604925031e-16\\
0.96558958783361	2.22044604925031e-16\\
0.997442574094731	2.22044604925031e-16\\
1.02929556035585	2.22044604925031e-16\\
1.06114854661697	2.22044604925031e-16\\
1.09300153287809	2.22044604925031e-16\\
1.12485451913921	2.22044604925031e-16\\
1.15670750540033	2.22044604925031e-16\\
1.18856049166145	2.22044604925031e-16\\
1.22041347792257	2.22044604925031e-16\\
1.25226646418369	2.22044604925031e-16\\
1.28411945044481	2.22044604925031e-16\\
1.31597243670593	2.22044604925031e-16\\
1.34782542296705	2.22044604925031e-16\\
1.37967840922817	2.22044604925031e-16\\
1.41153139548929	2.22044604925031e-16\\
1.44338438175042	2.22044604925031e-16\\
1.47523736801154	2.22044604925031e-16\\
1.50709035427266	2.22044604925031e-16\\
1.53894334053378	2.22044604925031e-16\\
1.5707963267949	2.22044604925031e-16\\
};

\addplot [color=black, line width=2.0pt, forget plot]
  table[row sep=crcr,forget plot]{%
0.01	0.113095007791622\\
0.0418529862611203	0.0452799341915577\\
0.0737059725222407	0.0102373823514063\\
0.105558958783361	0.00198272856930001\\
0.137411945044481	0.000196036201539599\\
0.169264931305602	1.43951088080358e-05\\
0.201117917566722	9.03151999282312e-07\\
0.232970903827842	3.25538539110444e-08\\
0.264823890088963	8.23480089835219e-10\\
0.296676876350083	8.98573383408922e-12\\
0.328529862611203	9.24505848333757e-14\\
0.360382848872324	7.44295703480745e-16\\
0.392235835133444	2.59635966587112e-16\\
0.424088821394564	2.53568878926248e-16\\
0.455941807655685	2.53568878926248e-16\\
0.487794793916805	2.53568878926248e-16\\
0.519647780177925	2.53568878926248e-16\\
0.551500766439046	2.53568878926248e-16\\
0.583353752700166	2.53568878926248e-16\\
0.615206738961286	2.53568878926248e-16\\
0.647059725222407	2.53568878926248e-16\\
0.678912711483527	2.53568878926248e-16\\
0.710765697744647	2.53568878926248e-16\\
0.742618684005768	2.53568878926248e-16\\
0.774471670266888	2.53568878926248e-16\\
0.806324656528008	2.53568878926248e-16\\
0.838177642789129	2.53568878926248e-16\\
0.870030629050249	2.53568878926248e-16\\
0.901883615311369	2.53568878926248e-16\\
0.93373660157249	2.53568878926248e-16\\
0.96558958783361	2.53568878926248e-16\\
0.997442574094731	2.53568878926248e-16\\
1.02929556035585	2.53568878926248e-16\\
1.06114854661697	2.53568878926248e-16\\
1.09300153287809	2.53568878926248e-16\\
1.12485451913921	2.53568878926248e-16\\
1.15670750540033	2.53568878926248e-16\\
1.18856049166145	2.53568878926248e-16\\
1.22041347792257	2.53568878926248e-16\\
1.25226646418369	2.53568878926248e-16\\
1.28411945044481	2.53568878926248e-16\\
1.31597243670593	2.53568878926248e-16\\
1.34782542296705	2.53568878926248e-16\\
1.37967840922817	2.53568878926248e-16\\
1.41153139548929	2.53568878926248e-16\\
1.44338438175042	2.53568878926248e-16\\
1.47523736801154	2.53568878926248e-16\\
1.50709035427266	2.53568878926248e-16\\
1.53894334053378	2.53568878926248e-16\\
1.5707963267949	2.53568878926248e-16\\
};

\addplot [color=mycolor2, dashed, line width=2.0pt, forget plot]
  table[row sep=crcr]{%
0.01	5.32615348859839\\
0.0418529862611203	0.209684111240814\\
0.0737059725222407	0.0295344909624121\\
0.105558958783361	0.00398449020120855\\
0.137411945044481	0.000412148545166887\\
0.169264931305602	3.01592288691005e-05\\
0.201117917566722	1.50250467089258e-06\\
0.232970903827842	4.98871641040201e-08\\
0.264823890088963	1.08959338166004e-09\\
0.296676876350083	1.55208568322328e-11\\
0.328529862611203	1.43552023392255e-13\\
0.360382848872324	1.07648162702875e-15\\
0.392235835133444	2.25322137172666e-16\\
0.424088821394564	2.22052674910745e-16\\
0.455941807655685	2.22044617654907e-16\\
0.487794793916805	2.22044604937876e-16\\
0.519647780177925	2.2204460492504e-16\\
0.551500766439046	2.22044604925031e-16\\
0.583353752700166	2.22044604925031e-16\\
0.615206738961286	2.22044604925031e-16\\
0.647059725222407	2.22044604925031e-16\\
0.678912711483527	2.22044604925031e-16\\
0.710765697744647	2.22044604925031e-16\\
0.742618684005768	2.22044604925031e-16\\
0.774471670266888	2.22044604925031e-16\\
0.806324656528008	2.22044604925031e-16\\
0.838177642789129	2.22044604925031e-16\\
0.870030629050249	2.22044604925031e-16\\
0.901883615311369	2.22044604925031e-16\\
0.93373660157249	2.22044604925031e-16\\
0.96558958783361	2.22044604925031e-16\\
0.997442574094731	2.22044604925031e-16\\
1.02929556035585	2.22044604925031e-16\\
1.06114854661697	2.22044604925031e-16\\
1.09300153287809	2.22044604925031e-16\\
1.12485451913921	2.22044604925031e-16\\
1.15670750540033	2.22044604925031e-16\\
1.18856049166145	2.22044604925031e-16\\
1.22041347792257	2.22044604925031e-16\\
1.25226646418369	2.22044604925031e-16\\
1.28411945044481	2.22044604925031e-16\\
1.31597243670593	2.22044604925031e-16\\
1.34782542296705	2.22044604925031e-16\\
1.37967840922817	2.22044604925031e-16\\
1.41153139548929	2.22044604925031e-16\\
1.44338438175042	2.22044604925031e-16\\
1.47523736801154	2.22044604925031e-16\\
1.50709035427266	2.22044604925031e-16\\
1.53894334053378	2.22044604925031e-16\\
1.5707963267949	2.22044604925031e-16\\
};

\node at (1.45,1e-8) {\large $\xi=3$};
\node at (1,1e-10) {\large $\xi=5$};
\node at (0.6,1e-12) {\large $\xi=10$};
\node at (0.22,1e-14) {\large $\xi=15$};

\end{axis}
\end{tikzpicture}%

%% file: fredrik_kspace_K0_est.tex
\definecolor{mycolor2}{rgb}{0.85000,0.32500,0.09800}%

\begin{tikzpicture}

\begin{axis}[%
  width=4.0187in,
  height=3.1693in,
  at={(1.011in,0.683in)},
scale only axis,
xmin=10,
xmax=100,
xlabel style={font=\color{white!15!black}},
xlabel={\Large $k_\infty$},
ymode=log,
ymin=1e-16,
ymax=50,
yminorticks=true,
ytick distance=1e4,
ylabel style={font=\color{white!15!black}},
ylabel={\Large Error, $\ewuf$},
axis background/.style={fill=white},
xmajorgrids,
ymajorgrids,
yminorgrids,
legend style={legend cell align=left, align=left, draw=white!15!black}
]
\addplot [color=black, line width=2.0pt]
  table[row sep=crcr]{%
15	0.000164513749962934\\
20	4.35237385260644e-07\\
25	1.06147057140914e-09\\
30	7.257351531141e-12\\
35	2.34257081790779e-13\\
40	1.79721773772723e-13\\
45	2.29175147143437e-13\\
50	2.14798921727684e-13\\
55	2.52934325192188e-13\\
60	1.79453703280291e-13\\
65	1.74791437706403e-13\\
70	2.04559144532665e-13\\
75	2.62631775625065e-13\\
80	2.71341890422365e-13\\
85	1.9183524518014e-13\\
90	2.50818616454282e-13\\
95	3.03601296273926e-13\\
100	2.10101307660692e-13\\
};
\addlegendentry{\large error}

\addplot [color=mycolor2, dashed, line width=2.0pt]
  table[row sep=crcr]{%
15	0.000327156241436366\\
20	1.23622117542392e-06\\
25	1.36733165516655e-09\\
30	4.2744793229672e-13\\
35	1.00340917395878e-14\\
40	1.00000007308891e-14\\
45	1.00000000000041e-14\\
50	1e-14\\
55	1e-14\\
60	1e-14\\
65	1e-14\\
70	1e-14\\
75	1e-14\\
80	1e-14\\
85	1e-14\\
90	1e-14\\
95	1e-14\\
100	1e-14\\
};
\addlegendentry{\large estimate}

\addplot [color=black, line width=2.0pt]
  table[row sep=crcr]{%
15	0.0120486799516313\\
20	0.000800014434421489\\
25	6.33674834247509e-05\\
30	4.54592191323645e-06\\
35	8.66281254659321e-08\\
40	2.60834149307192e-09\\
45	6.6484107527657e-11\\
50	2.05760168012603e-12\\
55	2.26019708843095e-13\\
60	1.73011754186149e-13\\
65	1.72422428442427e-13\\
70	1.93917447513588e-13\\
75	2.16082834447632e-13\\
80	2.18746586605151e-13\\
85	2.41727592341634e-13\\
90	1.99544258549233e-13\\
95	2.04076247786438e-13\\
100	2.02224905002687e-13\\
};

\addplot [color=mycolor2, dashed, line width=2.0pt]
  table[row sep=crcr]{%
15	0.0505067827359263\\
20	0.00428379829406637\\
25	0.000258691047130645\\
30	1.04889115440688e-05\\
35	2.76715970695472e-07\\
40	4.66158841244222e-09\\
45	4.95501141174086e-11\\
50	3.39386007231616e-13\\
55	1.13620952256383e-14\\
60	1.00034879130695e-14\\
65	1.0000005512361e-14\\
70	1.00000000053629e-14\\
75	1.00000000000032e-14\\
80	1e-14\\
85	1e-14\\
90	1e-14\\
95	1e-14\\
100	1e-14\\
};

\addplot [color=black, line width=2.0pt]
  table[row sep=crcr]{%
15	0.148785685994904\\
20	0.0500019848177853\\
25	0.01951834705372\\
30	0.00740057001147764\\
35	0.00189588040604518\\
40	0.000621855103368363\\
45	0.000152624793992617\\
50	3.85319160648467e-05\\
55	1.02118460004368e-05\\
60	2.54632472259935e-06\\
65	3.49705967578478e-07\\
70	6.92608338932532e-08\\
75	1.23587151558135e-08\\
80	2.06164894348472e-09\\
85	3.42552652207408e-10\\
90	6.10489218455591e-11\\
95	1.0023441208826e-11\\
100	1.89248087288288e-12\\
};

\addplot [color=mycolor2, dashed, line width=2.0pt]
  table[row sep=crcr]{%
15	1.10034085374538\\
20	0.346750911576057\\
25	0.113198695056166\\
30	0.0361012530979148\\
35	0.010899794884456\\
40	0.00305752353380076\\
45	0.000787239645504703\\
50	0.000184513727590428\\
55	3.91349725027229e-05\\
60	7.47855485302541e-06\\
65	1.28334845701105e-06\\
70	1.97253443444146e-07\\
75	2.70999712404202e-08\\
80	3.32247293860736e-09\\
85	3.63015063575012e-10\\
90	3.53153526369483e-11\\
95	3.06377366623959e-12\\
100	2.44723885842078e-13\\
};

\addplot [color=black, line width=2.0pt]
  table[row sep=crcr]{%
15	0.310659605101951\\
20	0.15511911499648\\
25	0.0876661845175065\\
30	0.0510999197136731\\
35	0.0237028815138398\\
40	0.0126961838751578\\
45	0.00594167030508805\\
50	0.00284011306411412\\
55	0.00134690993034753\\
60	0.000631712296357399\\
65	0.000209828214446134\\
70	9.14034756134113e-05\\
75	3.47459898376489e-05\\
80	1.26016477398661e-05\\
85	4.74605326167538e-06\\
90	1.66641316564815e-06\\
95	5.43957298331901e-07\\
100	1.92780843304415e-07\\
};

\addplot [color=mycolor2, dashed, line width=2.0pt]
  table[row sep=crcr]{%
15	3.38867538773454\\
20	1.36168931156309\\
25	0.607601535935349\\
30	0.283906989197564\\
35	0.134619457553696\\
40	0.0635703147993302\\
45	0.0295356064337676\\
50	0.0133900234805325\\
55	0.00588831111813903\\
60	0.00250078178413469\\
65	0.00102233785377252\\
70	0.000401259659565919\\
75	0.000150896794624903\\
80	5.4280175710187e-05\\
85	1.86518913405398e-05\\
90	6.11560153067792e-06\\
95	1.91153936563423e-06\\
100	5.69129114462404e-07\\
};

\node at (20,1e-10) {\large $\xi=3$};
\node at (52,1e-10) {\large $\xi=5$};
\node at (78,1e-10) {\large $\xi=10$};
\node at (90,1e-7) {\large $\xi=15$};

\end{axis}
\end{tikzpicture}%

%% file: fredrik_real_K1_est.tex
\definecolor{mycolor2}{rgb}{0.85000,0.32500,0.09800}%

\begin{tikzpicture}

\begin{axis}[%
    width=4.0187in,
    height=3.1693in,
    at={(1.011in,0.683in)},
  scale only axis,
  xmin=0,
  xmax=1.6,
  xtick distance=0.4,
  xlabel style={font=\color{white!15!black}},
  xlabel={\Large $r_c$},
  ymode=log,
  ymin=1e-16,
  ymax=50,
  ytick distance=1e4,
  yminorticks=true,
  axis background/.style={fill=white},
  axis x line*=bottom,
  axis y line*=left,
  ylabel style={font=\color{white!15!black}},
  ylabel={\Large Error, $\nabla\ewur$},
  xmajorgrids,
  ymajorgrids,
  yminorgrids,
  legend style={legend cell align=left, align=left, draw=white!15!black},
  execute at begin axis={
      \draw (rel axis cs:0,0) -- (rel axis cs:1,0)
            (rel axis cs:0,1) -- (rel axis cs:1,1);
  }
  ]
  \addplot [color=black, line width=2.0pt]
    table[row sep=crcr]{%
0.05	5.7604916431882\\
0.1024412526481	4.21447201503005\\
0.1548825052962	2.70022967740854\\
0.2073237579443	1.90966352289744\\
0.2597650105924	1.2914592432038\\
0.312206263240499	0.926646197544137\\
0.364647515888599	0.61421255701028\\
0.417088768536699	0.371944554607997\\
0.469530021184799	0.238746161027128\\
0.521971273832899	0.129390540532561\\
0.574412526480999	0.070772646832267\\
0.626853779129099	0.0382211315929039\\
0.679295031777199	0.0190307389817253\\
0.731736284425298	0.00862059172604254\\
0.784177537073398	0.00385554797175823\\
0.836618789721498	0.00186808926446057\\
0.889060042369598	0.000772848954495503\\
0.941501295017698	0.000306638922368827\\
0.993942547665798	0.000114244843624529\\
1.0463838003139	3.84174717544087e-05\\
1.098825052962	1.446989550999e-05\\
1.1512663056101	4.67945457888267e-06\\
1.2037075582582	1.48302259431247e-06\\
1.2561488109063	4.2139601601862e-07\\
1.3085900635544	1.1983043703736e-07\\
1.3610313162025	3.42513745492808e-08\\
1.4134725688506	8.70277711817837e-09\\
1.4659138214987	2.25361498937927e-09\\
1.5183550741468	5.08542293091845e-10\\
1.5707963267949	1.19616154128577e-10\\
};
\addlegendentry{\large error}

\addplot [color=mycolor2, dashed, line width=2.0pt]
  table[row sep=crcr]{%
0.05	34.1683683381184\\
0.1024412526481	15.5193236964791\\
0.1548825052962	9.09076629915173\\
0.2073237579443	5.72414637881964\\
0.2597650105924	3.66469515908466\\
0.312206263240499	2.32775548549283\\
0.364647515888599	1.44799790108105\\
0.417088768536699	0.875339802959652\\
0.469530021184799	0.511690725855993\\
0.521971273832899	0.288264928044564\\
0.574412526480999	0.156129066954895\\
0.626853779129099	0.081154467767802\\
0.679295031777199	0.0404290018053059\\
0.731736284425298	0.0192828578460404\\
0.784177537073398	0.00879807956291174\\
0.836618789721498	0.00383753827146212\\
0.889060042369598	0.00159929674012334\\
0.941501295017698	0.000636532805306694\\
0.993942547665798	0.000241859838515907\\
1.0463838003139	8.77039431191546e-05\\
1.098825052962	3.03436918671867e-05\\
1.1512663056101	1.00140506152035e-05\\
1.2037075582582	3.15177315346145e-06\\
1.2561488109063	9.45860239717093e-07\\
1.3085900635544	2.70619024397474e-07\\
1.3610313162025	7.38056254144018e-08\\
1.4134725688506	1.91852560734564e-08\\
1.4659138214987	4.75275817670848e-09\\
1.5183550741468	1.12197141897522e-09\\
1.5707963267949	2.52370293635306e-10\\
};
\addlegendentry{\large estimate}

\addplot [color=black, line width=2.0pt, forget plot]
  table[row sep=crcr,forget plot]{%
0.05	4.54102723049103\\
0.1024412526481	2.85415418303376\\
0.1548825052962	1.36297075764668\\
0.2073237579443	0.683592579831275\\
0.2597650105924	0.285583485224781\\
0.312206263240499	0.125955053148722\\
0.364647515888599	0.0466926016993664\\
0.417088768536699	0.0138927198888\\
0.469530021184799	0.00432771850786486\\
0.521971273832899	0.00102018500455027\\
0.574412526480999	0.000232538605340891\\
0.626853779129099	4.46666290613793e-05\\
0.679295031777199	7.52297199864677e-06\\
0.731736284425298	1.01656393682973e-06\\
0.784177537073398	1.21196815415518e-07\\
0.836618789721498	1.53330866508591e-08\\
0.889060042369598	1.56732987454583e-09\\
0.941501295017698	1.29837696998836e-10\\
0.993942547665798	9.56806309324953e-12\\
1.0463838003139	5.85967164726229e-13\\
1.098825052962	3.5487673738186e-14\\
1.1512663056101	2.48263601696848e-15\\
1.2037075582582	1.54461765292022e-15\\
1.2561488109063	1.54544307542172e-15\\
1.3085900635544	1.54283979330286e-15\\
1.3610313162025	1.54468547167284e-15\\
1.4134725688506	1.54682934007579e-15\\
1.4659138214987	1.54985112865269e-15\\
1.5183550741468	1.54984365235432e-15\\
1.5707963267949	1.54900527160605e-15\\
};

\addplot [color=mycolor2, dashed, line width=2.0pt, forget plot]
  table[row sep=crcr]{%
0.05	19.6971644713985\\
0.1024412526481	7.87232110642142\\
0.1548825052962	3.71588225708417\\
0.2073237579443	1.72657134536188\\
0.2597650105924	0.746972999481753\\
0.312206263240499	0.293615220007324\\
0.364647515888599	0.1035056836975\\
0.417088768536699	0.0324719902738707\\
0.469530021184799	0.00902104079597177\\
0.521971273832899	0.00221176347819312\\
0.574412526480999	0.000477429575615314\\
0.626853779129099	9.05727537696517e-05\\
0.679295031777199	1.50806008502651e-05\\
0.731736284425298	2.20149862784371e-06\\
0.784177537073398	2.8153811744967e-07\\
0.836618789721498	3.15199298080365e-08\\
0.889060042369598	3.08763324808874e-09\\
0.941501295017698	2.64522605881193e-10\\
0.993942547665798	1.98123585336351e-11\\
1.0463838003139	1.2970832311784e-12\\
1.098825052962	7.43924269400179e-14\\
1.1512663056101	3.92749078721973e-15\\
1.2037075582582	3.83716754968973e-16\\
1.2561488109063	2.2820399445699e-16\\
1.3085900635544	2.22249475868774e-16\\
1.3610313162025	2.22050553374214e-16\\
1.4134725688506	2.2204475567424e-16\\
1.4659138214987	2.22044608259197e-16\\
1.5183550741468	2.22044604989383e-16\\
1.5707963267949	2.22044604926115e-16\\
};

\addplot [color=black, line width=2.0pt, forget plot]
  table[row sep=crcr,forget plot]{%
0.05	2.61206417704889\\
0.1024412526481	0.85447673065652\\
0.1548825052962	0.129636260544633\\
0.2073237579443	0.0156588749774309\\
0.2597650105924	0.00101778097123938\\
0.312206263240499	4.60711832580211e-05\\
0.364647515888599	1.07139055721104e-06\\
0.417088768536699	1.56844543493233e-08\\
0.469530021184799	1.38741654742976e-10\\
0.521971273832899	6.65467346449667e-13\\
0.574412526480999	2.44476053000651e-15\\
0.626853779129099	9.4735184027551e-16\\
0.679295031777199	9.57037849319632e-16\\
0.731736284425298	9.32159293032068e-16\\
0.784177537073398	8.96308276411185e-16\\
0.836618789721498	8.92441076042083e-16\\
0.889060042369598	8.8254760860919e-16\\
0.941501295017698	8.75091417468303e-16\\
0.993942547665798	8.68673650567359e-16\\
1.0463838003139	8.82238992629828e-16\\
1.098825052962	8.79615107443443e-16\\
1.1512663056101	8.46137112305942e-16\\
1.2037075582582	8.36488923178684e-16\\
1.2561488109063	8.33240365864172e-16\\
1.3085900635544	8.26395420994308e-16\\
1.3610313162025	8.39970118053283e-16\\
1.4134725688506	8.36458490388611e-16\\
1.4659138214987	8.3478600713527e-16\\
1.5183550741468	8.2539122015036e-16\\
1.5707963267949	8.10940708051831e-16\\
};

\addplot [color=mycolor2, dashed, line width=2.0pt, forget plot]
  table[row sep=crcr]{%
0.05	8.1647614461889\\
0.1024412526481	1.79165540726327\\
0.1548825052962	0.307378519289188\\
0.2073237579443	0.0343640209182643\\
0.2597650105924	0.00236799191807871\\
0.312206263240499	9.81426392091781e-05\\
0.364647515888599	2.41487810421136e-06\\
0.417088768536699	3.50057803531697e-08\\
0.469530021184799	2.97464847489887e-10\\
0.521971273832899	1.47699312896931e-12\\
0.574412526480999	4.49499277278426e-15\\
0.626853779129099	2.29237565440241e-16\\
0.679295031777199	2.22051639996894e-16\\
0.731736284425298	2.22044608918553e-16\\
0.784177537073398	2.22044604926346e-16\\
0.836618789721498	2.22044604925032e-16\\
0.889060042369598	2.22044604925031e-16\\
0.941501295017698	2.22044604925031e-16\\
0.993942547665798	2.22044604925031e-16\\
1.0463838003139	2.22044604925031e-16\\
1.098825052962	2.22044604925031e-16\\
1.1512663056101	2.22044604925031e-16\\
1.2037075582582	2.22044604925031e-16\\
1.2561488109063	2.22044604925031e-16\\
1.3085900635544	2.22044604925031e-16\\
1.3610313162025	2.22044604925031e-16\\
1.4134725688506	2.22044604925031e-16\\
1.4659138214987	2.22044604925031e-16\\
1.5183550741468	2.22044604925031e-16\\
1.5707963267949	2.22044604925031e-16\\
};

\addplot [color=black, line width=2.0pt, forget plot]
  table[row sep=crcr,forget plot]{%
0.05	1.48673484813805\\
0.1024412526481	0.169979326853669\\
0.1548825052962	0.00438942792051825\\
0.2073237579443	5.09994845993657e-05\\
0.2597650105924	1.4505478065103e-07\\
0.312206263240499	1.74476740061884e-10\\
0.364647515888599	4.30242007756412e-14\\
0.417088768536699	7.76967612598198e-16\\
0.469530021184799	7.48272924795623e-16\\
0.521971273832899	7.3707341071198e-16\\
0.574412526480999	7.15306543377005e-16\\
0.626853779129099	6.89019061126477e-16\\
0.679295031777199	6.69929628734263e-16\\
0.731736284425298	6.59280893426719e-16\\
0.784177537073398	6.33149524294782e-16\\
0.836618789721498	6.10535906548813e-16\\
0.889060042369598	5.99985454034106e-16\\
0.941501295017698	5.80355342995224e-16\\
0.993942547665798	5.77325542426119e-16\\
1.0463838003139	5.65679410815925e-16\\
1.098825052962	5.65055942079972e-16\\
1.1512663056101	5.56391824887679e-16\\
1.2037075582582	5.45047669874988e-16\\
1.2561488109063	5.35953699287807e-16\\
1.3085900635544	5.35229346707737e-16\\
1.3610313162025	5.35947861491264e-16\\
1.4134725688506	5.29198925247322e-16\\
1.4659138214987	5.291147223627e-16\\
1.5183550741468	5.37818503330621e-16\\
1.5707963267949	5.39986393674738e-16\\
};

\addplot [color=mycolor2, dashed, line width=2.0pt, forget plot]
  table[row sep=crcr]{%
0.05	3.98231138710185\\
0.1024412526481	0.321711085047354\\
0.1548825052962	0.0102168777686005\\
0.2073237579443	0.000106315140644293\\
0.2597650105924	3.42871301866058e-07\\
0.312206263240499	3.34412140409197e-10\\
0.364647515888599	9.75876659427255e-14\\
0.417088768536699	2.30442103042809e-16\\
0.469530021184799	2.22044818405686e-16\\
0.521971273832899	2.22044604926626e-16\\
0.574412526480999	2.22044604925031e-16\\
0.626853779129099	2.22044604925031e-16\\
0.679295031777199	2.22044604925031e-16\\
0.731736284425298	2.22044604925031e-16\\
0.784177537073398	2.22044604925031e-16\\
0.836618789721498	2.22044604925031e-16\\
0.889060042369598	2.22044604925031e-16\\
0.941501295017698	2.22044604925031e-16\\
0.993942547665798	2.22044604925031e-16\\
1.0463838003139	2.22044604925031e-16\\
1.098825052962	2.22044604925031e-16\\
1.1512663056101	2.22044604925031e-16\\
1.2037075582582	2.22044604925031e-16\\
1.2561488109063	2.22044604925031e-16\\
1.3085900635544	2.22044604925031e-16\\
1.3610313162025	2.22044604925031e-16\\
1.4134725688506	2.22044604925031e-16\\
1.4659138214987	2.22044604925031e-16\\
1.5183550741468	2.22044604925031e-16\\
1.5707963267949	2.22044604925031e-16\\
};

\node at (1.5,5e-8) {\large $\xi=3$};
\node at (1.05,5e-10) {\large $\xi=5$};
\node at (0.62,5e-12) {\large $\xi=10$};
\node at (0.25,1e-14) {\large $\xi=15$};

\end{axis}
\end{tikzpicture}%

%% file: fredrik_kspace_K1_est.tex
\definecolor{mycolor2}{rgb}{0.85000,0.32500,0.09800}%

\begin{tikzpicture}

\begin{axis}[%
  width=4.0187in,
  height=3.1693in,
  at={(1.011in,0.683in)},
scale only axis,
xmin=10,
xmax=100,
xlabel style={font=\color{white!15!black}},
xlabel={\Large $k_\infty$},
ymode=log,
ymin=1e-16,
ymax=50,
yminorticks=true,
ytick distance=1e4,
axis background/.style={fill=white},
xmajorgrids,
ymajorgrids,
yminorgrids,
ylabel={\Large Error, $\nabla\ewuf$},
legend style={legend cell align=left, align=left, draw=white!15!black}
]
\addplot [color=black, line width=2.0pt]
  table[row sep=crcr]{%
15	0.000223508792078935\\
20	7.06575806615659e-07\\
25	1.65383736926523e-09\\
30	1.08591507809574e-11\\
35	2.59217122834068e-14\\
40	1.45949555026672e-14\\
45	1.21927993741473e-14\\
50	1.49853924526281e-14\\
55	1.19232024485522e-14\\
60	1.27405808090365e-14\\
65	1.41788229237255e-14\\
70	1.61572242835034e-14\\
75	8.87960041605652e-15\\
80	1.89522231140416e-14\\
85	1.1950906512021e-14\\
90	9.69893477724826e-15\\
95	1.08154585323494e-14\\
100	3.99172644535282e-15\\
};
\addlegendentry{\large error}

\addplot [color=mycolor2,dashed, line width=2.0pt]
  table[row sep=crcr]{%
15	8.29890322968294e-05\\
20	4.18119492092127e-07\\
25	5.78085998227111e-10\\
30	2.21786274463146e-13\\
35	1.00201786196169e-14\\
40	1.00000004944083e-14\\
45	1.00000000000031e-14\\
50	1e-14\\
55	1e-14\\
60	1e-14\\
65	1e-14\\
70	1e-14\\
75	1e-14\\
80	1e-14\\
85	1e-14\\
90	1e-14\\
95	1e-14\\
100	1e-14\\
};
\addlegendentry{\large estimate}

\addplot [color=black, line width=2.0pt]
  table[row sep=crcr]{%
15	0.0127717651488208\\
20	0.00122278714031825\\
25	0.000107833592800726\\
30	1.08503643490291e-05\\
35	2.41212061299357e-07\\
40	6.1937241184203e-09\\
45	1.31454503099482e-10\\
50	3.84606269750877e-12\\
55	1.8991999362553e-13\\
60	3.30236093816703e-14\\
65	1.55079328676946e-14\\
70	1.06263050574535e-14\\
75	1.50109592879191e-14\\
80	1.98250118073726e-14\\
85	1.94372616971614e-14\\
90	1.28065091882696e-14\\
95	1.79347304552381e-14\\
100	1.01916679164348e-14\\
};

\addplot [color=mycolor2, dashed, line width=2.0pt]
  table[row sep=crcr]{%
15	0.0128119996479793\\
20	0.00144888864736854\\
25	0.000109369813420615\\
30	5.32142336337704e-06\\
35	1.63786645085305e-07\\
40	3.1533388575592e-09\\
45	3.77103922916693e-11\\
50	2.88516914205476e-13\\
55	1.12669126000952e-14\\
60	1.00035391051594e-14\\
65	1.00000060593723e-14\\
70	1.00000000063486e-14\\
75	1.00000000000041e-14\\
80	1e-14\\
85	1e-14\\
90	1e-14\\
95	1e-14\\
100	1e-14\\
};

\addplot [color=black, line width=2.0pt]
  table[row sep=crcr]{%
15	0.109778576688432\\
20	0.0484007219129095\\
25	0.022569519996301\\
30	0.0113500725006431\\
35	0.00363941963796222\\
40	0.00124352125592895\\
45	0.000410857145449872\\
50	0.000112833464206897\\
55	2.92743710502828e-05\\
60	1.00916386334846e-05\\
65	1.48442950072973e-06\\
70	2.77107429651655e-07\\
75	4.45642729745857e-08\\
80	6.96086573397181e-09\\
85	1.0146631466553e-09\\
90	1.38882728159246e-10\\
95	2.16168752475032e-11\\
100	3.5947017648926e-12\\
};

\addplot [color=mycolor2, dashed, line width=2.0pt]
  table[row sep=crcr]{%
15	0.27912970475569\\
20	0.11728315903136\\
25	0.0478597091540044\\
30	0.0183160897951388\\
35	0.00645173042880306\\
40	0.00206833198516923\\
45	0.000599115441460403\\
50	0.000156023605793966\\
55	3.64015531116941e-05\\
60	7.58860030891223e-06\\
65	1.41075382184909e-06\\
70	2.335158934628e-07\\
75	3.43735502556189e-08\\
80	4.49516806651878e-09\\
85	5.21837470828104e-10\\
90	5.37476754301166e-11\\
95	4.91632979551865e-12\\
100	4.06966663471713e-13\\
};

\addplot [color=black, line width=2.0pt]
  table[row sep=crcr]{%
15	0.196779091114269\\
20	0.121639758818382\\
25	0.078830173603953\\
30	0.0562707107398258\\
35	0.02995330458782\\
40	0.0170131504600528\\
45	0.00978852569161039\\
50	0.00523416880183543\\
55	0.00276270635692042\\
60	0.00161024489276607\\
65	0.000643544612191186\\
70	0.000299729730387568\\
75	0.000117841705362472\\
80	4.94268562839396e-05\\
85	1.89180701907717e-05\\
90	6.90628325891267e-06\\
95	2.44487816889969e-06\\
100	8.26691825693381e-07\\
};

\addplot [color=mycolor2, dashed, line width=2.0pt]
  table[row sep=crcr]{%
15	0.859624503872212\\
20	0.460570452009674\\
25	0.256890176843124\\
30	0.144041147089966\\
35	0.0796830087001205\\
40	0.0430035988122354\\
45	0.0224775746347914\\
50	0.011322516608171\\
55	0.00547703642531289\\
60	0.00253758029364407\\
65	0.00112383119920783\\
70	0.000475025971741511\\
75	0.000191397212652612\\
80	7.34388838017409e-05\\
85	2.68124890241236e-05\\
90	9.30845283746691e-06\\
95	3.07116490985441e-06\\
100	9.62515098196378e-07\\
};

\node at (20,1e-10) {\large $\xi=3$};
\node at (52,1e-10) {\large $\xi=5$};
\node at (78,1e-10) {\large $\xi=10$};
\node at (93,1.1e-7) {\large $\xi=15$};

\end{axis}
\end{tikzpicture}%

%% file: pux.tex
PUX is a numerical method to construct  high-order function extensions \cite{FRYKLUNDPUX}. Discs, called partitions, are distributed along the boundary $\bdry$; in each partition there is a Gaussian basis that interpolates $\mhrhs$ at the nodes in $\domain$ that falls in the partition, then the interpolant is evaluated at the other nodes inside the partition, i.e. those not in $\domain$. Multiple partitions overlap and their respective extensions are weighted together with a high-order partition of unity to form a global extension. This has two purposes: to accommodate for missing data in the slices $\gc{\slice{n}{m}}$ and to construct the extended function $\mhrhse$ \eqref{eq:mhrhsecond}.
\subsection{Interpolation with radial basis functions}
\label{ss:interprbf}
We start out with a brief introduction of interpolation with RBFs and how to construct a local extension, a notion specified below. Recall that the computational domain is a box $\suppbox=[-\tfrac{L}{2},\tfrac{L}{2}]^{2}$ in $\mathbb{R}^{2}$ that contains $\bar{\domain}$. $\suppbox$ has an underlying grid $\X$ consisting of $\Nu^{2}$ uniformly distributed points. Let $\rbf^{q}$ be a compactly supported univariate symmetric positive function, such that $\rbf^{q}(\x) = \rbf^{q}(\|\x-\rbfcenter\|)$ for some point $\rbfcenter$ in $\X$. We refer to this as a compactly supported radial basis function (RBF), centered at $\rbfcenter$. The superscript $q$ indicates the highest regularity subset $C^{q}_{0}$ of $C_{0}$ that $\psi^{q}$ is a member of. The support of $\rbf^{q}$ defines a partition $\partition$, i.e.
$\partition = \text{supp}(\rbf^{q})\subset\mathbb{R}^{2}$, which is a disc with a radius denoted $R$. In the partition $\partition$ there are $N_{\rbfgauss}$ points $\{\gausscenteri\}_{i=1}^{N_{\rbfgauss}}$ drawn from the quasi-uniform Vogel-node distribution, defined as
\begin{equation}
  \gausscenteri = \sqrt{i/N_{\rbfgauss}}\, (\cos{(i\pi(3-\sqrt{5}))},\sin{(i\pi(3-\sqrt{5}))}),\quad i = 1,\dots,N_{\rbfgauss},
\end{equation}
in local coordinates of the partition. This provides a near-optimal distribution of basis functions for interpolation \cite{LarssonElisabeth2017ALSR}. Each point $\gausscenteri$ is the center for a Gaussian $\rbfgauss_{i}(\x) = \text{exp}(-\varepsilon^2 \|\x-\gausscenteri\|^2)$, where $\varepsilon$ is a shape-parameter setting the width of the Gaussian. This forms a collection of Gaussians  $\{\rbfgauss_{i}\}_{i=1}^{N_{\rbfgauss}}$, distributed over the partition $\partition$, and it defines the interpolation basis in $\partition$. The standard form of RBF approximation $I_{f}$ of a function $f$ in this basis is then
\begin{equation}
  \label{eq:stdrbfinterp}
 I_{f}(\x) = \sum\limits_{i = 1}^{N_{\rbfgauss}}\lambda_{i}\rbfgauss_{i}(\x).
\end{equation}
The unknown coefficients $\{\lambda_{j}\}$ can be determined by collocating at the centers $\{\gausscenteri\}$ of the RBFs and solving the resulting linear system $f_{\gausscenter} = \rbfmat_{\gausscenter}\rbfweights$ where $f_{\gausscenter,i} = f(\gausscenteri)$, $\rbfmat_{\gausscenter,i,j} = \rbfgauss_{i}(\x_{j})$, and $\rbfweights_{i} = \lambda_{i}$.

The condition number of $\rbfmat_{\gausscenter}$ is highly dependent on the shape parameter $\varepsilon$ in the Gaussian basis: the smaller the shape parameter, to a limit, the smaller interpolation error, but the higher condition number; values of $10^{18}$ are not uncommon for reasonable values of $\varepsilon$ \cite{LarssonE2005Taca,Fasshauer:2007:MAM:1506263}. Also, the interpolation weights $\{\lambda_{j}\}$ are large in magnitude and are of varying sign. We can circumvent these drawbacks by proceeding as follows. Let $\Xpart$ denote the points $\{\x_{\partition,i}\}$ in $\X$ within Euclidean distance $R$ of $\rbfcenter$. Furthermore, let $f_{\partition}$ be a vector corresponding to the values of $f$ at $\Xpart$. Then using \eqref{eq:stdrbfinterp} and $\rbfweights = \rbfmat_{\gausscenter}^{-1}f_{\gausscenter}$ gives 
\begin{equation}
  \label{eq:rbfinterpreform}
  f_{\partition} = \rbfmat_{\partition}\rbfweights  = \rbfmat_{\partition}\rbfmat_{\gausscenter}^{-1}f_{\gausscenter} = Af_{\gausscenter},
\end{equation}
where $\rbfmat_{\partition,i,j} = \rbfgauss_{i}(\x_{\partition,j})$ and $\rbfmat = \rbfmat_{\partition}\rbfmat_{\gausscenter}^{-1}$. We now avoid using the interpolation weights, and are left with dealing with the ill-conditioning of the problem. It can be reduced significantly by using the algorithm RBF-QR \cite{MR2801193}. The drawback is that RBF-QR is an $\ordo{N^{3}}$ algorithm, and in PUX we consider several partitions; it is unfeasible to apply RBF-QR for every partition every time PUX is called. We return to how to circumvent this bottleneck later in this section.

We now describe how to create a local extension $\mhrhse$ for the partition $\partition$. The setting is that $f$ is only partially known in $\Xpart$, and $f_{\gausscenter}$ is unknown and must be solved for. To do so, we introduce the following notation. Let the set $\Xpart$ be split into two disjoint sets: one set $\Xint$ where function data $f$ is known and $\Xext$ where function data does not exist. Furthermore, let $\mhrhse$ be a vector with unknown elements, corresponding to the unknown values of the function $\mhrhs$ in $\Xext$. We refer to them as the local extension of $\mhrhs$.  Let $f_{\partition}= [\mhrhs\,\mhrhse]^\intercal$ where $f$ is evaluated at $\Xint$, and suppose that both $\Xint$ and $\Xext$ are nonempty, and the interpolation matrix $\rbfmat$ is decomposed as 
\begin{equation}
  \rbfmat = \begin{bmatrix}
    \rbfmat_{I}\\
    \rbfmat_{E}
  \end{bmatrix},
  \label{eq:rbfmatdecomp}
\end{equation}
then \eqref{eq:rbfinterpreform} takes the form
\begin{equation}
\begin{bmatrix}
    \mhrhs\\
    \mhrhse
    \end{bmatrix} =  \begin{bmatrix}
    \rbfmatint\\
    \rbfmatext
    \end{bmatrix} f_{\gausscenter}.
\end{equation}
We can obtain $f_{\gausscenter}$ by solving $ \rbfmatint  f_{\gausscenter} = \mhrhs$ in the least squares sense, as the basis locations and data locations are decoupled \cite{LarssonElisabeth2017ALSR}. %Furthermore, we set $N_{\rbfgauss}$ to be one at least one fourth of the number of elements in $\Xint$ \cite{puxapprox}.
It is now straightforward to compute $\mhrhse$ through
\begin{equation}
  \label{eq:localsystem}
\mhrhse = \rbfmatext f_{\gausscenter}.
\end{equation}

The bottleneck of having to compute $\rbfmat$ with RBF--QR for every partition can be circumvented entirely, making it sufficient to call RBF--QR once as a precomputation step. It is based on the following observation. Assume that the number of Gaussian basis-functions is fixed, and that the center $\rbfcenter$ belongs to $\X$, then $\rbfmat$ is independent of the position of $\rbfcenter\in\X$. It is also independent of the distribution of $\Xpart$ into $\Xint$ and $\Xext$, as it merely affects the decomposition \eqref{eq:rbfmatdecomp}, i.e. the arrangement of the rows of the matrix $\rbfmat$. Thus, RBF-QR is called once as a precomputation step to compute $\rbfmat$, which is reused for all partitions for all time-steps.

\subsection{Partition of unity}
We now consider a collection of $N_{\partition}$ partitions $\{\partition_{i}\}_{i=1}^{N_{\partition}}$, each associated with a compactly supported RBF $\rbf^{q}_{i}=\rbf^{q}(\,\|\cdot- \rbfcenteri\|)$ where $\rbfcenteri\in\X$ for $i = 1,\ldots,N_{\partition}$. Thus all partitions are discs with radius $R$. Now, for each partition associate a weight function $w_{i}$, defined as
\begin{equation}
\label{eq:weightfunction}
w_{i}(\x) = \frac{\rbf^{q}_{i}(\x)}{\sum\limits_{j=1}^{N_{\partition}}\rbf^{q}_{j}(\x)}.
\end{equation}
If $\rbf^{q}_{i}\in C^{q}_{0}$, for $i = 1,\ldots,N_{\partition}$, then $w_{i}$ belongs to $C^{q}_{0}$. Also, if $\rbf^{q}_{i}$ is a positive function, then so is $w_{i}$. The set of weights $\{w_{i}\}_{i = 1}^{N_{\partition}}$ forms a partition of unity, meaning
\begin{equation}
\label{eq:POU}
\sum\limits_{i=1}^{N_{\partition}}w_{i}(\x)= 1,\quad \forall\x\in \bigcup\limits_{i = 1}^{N_{\partition}}\bar{\partitioni},
\end{equation}
which is referred to in the literature as Shepard's method \cite{shepard}. When a single partition is considered, we refer to the associated interpolation error, interpolant and so on as local, as opposed to global.

Assume that there exists a local extension $\mhrhse_{i}$ for each partition $\partitioni$. By weighting them together through a partition of unity we construct
%\begin{equation}
%\mhrhse(\x) = \sum\limits_{i=1}^{N_{\partition}}w_{i}(\x)\mhrhse_{i}(\x),
%\label{eq:fe_notcompact}
%\end{equation}
\begin{equation}
 \label{eq:fe_notcompact}
      \tilde{\mhrhse}(\x) = \begin{cases}\mhrhs(\x),\quad \x\in\X_{\domain},\\\sum\limits_{i = 1}^{N_{\partition}}w_{i}(\x)\mhrhse_{i}(\x) ,\quad \x \in \bigcup\limits_{i=1}^{N_{\partition}}\X_{E,i},\\
      0, \quad \text{otherwise},
      \end{cases}
\end{equation}
and we refer to $\tilde{\mhrhse}$ as a \emph{global extension}, or simply an extension, of $\mhrhs$. This kind of extension has no imposed global regularity such that $\tilde{f}^{e}\in C^{q}(\mathbb{R}^2)$; in most events $\tilde{f}^{e}$ will be discontinuous, as it is extended with zero outside of the partition. The global extension $\tilde{f}^{e}$ is used for extending data into slices arising from time-dependent geometry, see section \ref{sss:imexsdctimedep}, as here, only data from within the partitions is used.

\subsection{Function extension with compact support}
 We want an extension with a specified regularity as it is extended by zero outside its support, i.e. the global extension should have regularity $q\geq 0$ in $\mathbb{R}^{2}$. To achieve this, the boundary of  $\domain$ is covered with overlapping partitions, i.e. the centers $\{ \rbfcenteri\}_{i=1}^{N_{\partition}}$ for the partitions $\{ \partitioni\}_{i=1}^{N_{\partition}}$ are distributed uniformly with respect to arc length along the boundary, as in the right image in figure \ref{fig:nemhext}. The centers are then shifted to the closest point in $\X\cap\domain$, which is motivated below.

Refer to the set of partitions $\{ \partitioni\}_{i=1}^{N_{\partition}}$ as extension partitions and now introduce the zero partitions
$\{\partitioni^{0}\}_{i = N_{\partition} + 1}^{N_{\partition} + N_{\partition}^{0}}$. The corresponding RBFs $\rbf^{q}$ are included in \eqref{eq:weightfunction} to form the set of weights $\{w^{0}_{i}\}_{i=1}^{N_{\partition} + N_{\partition}^{0}}$, that is
\begin{equation}
\label{eq:weightfunctionzero}
w^{0}_{i}(\x) = \frac{\rbf^{q}_{i}(\x)}{\sum\limits_{j=1}^{N_{\partition} + N^{0}_{\partition}}\rbf^{q}_{j}(\x)}.
\end{equation}
The weight functions \eqref{eq:weightfunctionzero} also form a partition of unity \eqref{eq:POU}. The zero partitions are distributed such that they overlap the extension partitions, but do not intersect $\bar{\Omega}$. The associated local extension $\mhrhse_i$ is set to be identically equal to zero for $i = N_{\partition} + 1,\ldots, N_{\partition}  + N^{0}_{\partition}$. Hence, as these zero valued functions are blended with the local extensions in the first layer of partitions, i.e. the extension partitions, the global extension will be forced to zero over the overlapping region. Therefore zero partitions should be distributed such that
$f^{e}$ has a controlled transition to zero, see figure \ref{fig:nemhext}. Such a global extension $\mhrhse$ of $\mhrhs$ is given by
\begin{equation}
  \label{eq:functionextension}
      \mhrhse(\x) = \begin{cases}\mhrhs(\x),\quad \x\in\X_{\domain},\\\sum\limits_{i = 1}^{N_{\partition}+N_{\partition}^{0}}w^{0}_{i}(\x)\mhrhse_{i}(\x) = \sum\limits_{i = 1}^{N_{\partition}}w^{0}_{i}(\x)\mhrhse_{i}(\x),\quad \x \in \bigcup\limits_{i=1}^{N_{\partition}}\X_{E,i},\\
      0, \quad \text{otherwise}.
      \end{cases}
\end{equation}
Here $\X_{\domain} = \X\cap\domain$ and $\X_{E,i}$ is $\Xext$ for partition $\partitioni$. Thus the only difference between how the extensions \eqref{eq:fe_notcompact} and \eqref{eq:functionextension} are computed is how the weight functions are constructed.

The global extension inherits the regularity of $\rbf^{q}$ where it is extended to zero. However, there is no guarantee that $f^{e}$, or $\tilde{f}^{e}$, is of regularity $q$ over the boundary of $\Omega$. Inside we have $f^{e}= f$, while outside $f^{e}$ is equal to an extension based on a weighted sum of local extensions. However, in practice we observe that the radial basis functions $\rbf^{q}$ sets the regularity of the approximation.

%Moreover, the $(q+1)$th derivative of $\rbf^{q}$ is of bounded variation
As $\rbf^{q}$ we use one of the compactly supported Wu functions, which are tabulated after their regularity $q$, see table \ref{tab:Wu} or \cite{Fasshauer:2007:MAM:1506263}. There are other options, but the Wu functions have compact support and are simple to implement. Note that they have lower regularity at the origin, e.g. the Wu function listed as $C^{4}$ is only $C^{2}$ at that point. The partitions centers $\{\rbfcenteri\}_{i = 1}^{N_{\partition}}$ are set to be nodes on the regular grid that are the closest to be boundary, yet still in $\X_{\domain}$. Thus evaluation of weight functions at the origin is omitted and higher regularity is maintained. With this, we have described how local extensions are combined into a global one. For instructions on how to set parameters to achieve high accuracy and a computational complexity linear in the number of grid points $\Nu$, as well as proofs for the regularity of the extension, see \cite{FRYKLUNDPUX,puxapprox}.

When the domains are time-dependent, both definition \eqref{eq:fe_notcompact} and \eqref{eq:functionextension} are used for the global function extension. For each domain at each substep in $[\tn,\tnp]$ we construct an extension; first \eqref{eq:fe_notcompact} is created for each domain which the modified Helmholtz equation is not currently solved on, i.e. create an extension without enforcing compact support with zero partitions. By doing so, missing data in the slices are created, which can be combined to create a right-hand side for the modified Helmholtz equation, see section \ref{sss:imexsdctimedep}. The right-hand side for the modified Helmholtz equation is in turn extended with \eqref{eq:functionextension} such that it has a prescribed regularity $q$ in the plane.

%These extensions enter in the right hand $f$ for the modified Helmholtz equation, and then extended with PUX using definition \eqref{eq:functionextension}, i.e.  , from the domain we solve on.

%Now data may not be available at all nodes in $\domain$. Thus we redefine $\X_{\domain}$ to be the set of nodes in $\X$ that are in $\domain$ where data exists, i.e. the slice where data does not exist is excluded. For a given partition $\partitioni$ the set $\X_{E,i}$ now also contains the slice. The partitions are distributed as before, thus zero partitions do not intersect the slice, and extrapolated data in the slice is not suppressed to zero.

Given two time stamps $t_j$ and $t_{m+1}$, using the notation from section \eqref{sss:imexsdctimedep}, available data $\mhrhs_{j}$  in $\domain_{j}$ are extended into $\slice{t_{j}}{t_{m+1}}\cup \suppbox$, meaning that for each time-step with SDC of order $K$, there are $K-1 + K(K-1)^2$ calls for PUX: there are $K-1$ extensions for the initial approximations $\{\ouakm{0}{m}\}_{m = 1}^{K-1}$. Then there are $K-1$ applications of the correction equation, each with $K$ values $\mhrhs_j$. The total cost is efficiently reduced by precomputing a single interpolation matrix $\rbfmat$.

\begin{table}
\centering
{\def\arraystretch{1.3}
\begin{tabular}{cc}
\toprule
\textbf{Regularity} & $\rbf^{q}(r)$ \\
\midrule
$\rbf^{1}\in C^{1}_{0}$ & $(1- r)_{+}^{2}(2+ r)$ \\
$\rbf^{2}\in C^{2}_{0}$ & $(1- r)_{+}^{3}(8+9 r+3r^{2})$ \\
$\rbf^{3}\in C^{3}_{0}$ & $(1-r)_{+}^{4}(4+16r+12r^{2}+3r^{3})$ \\
$\rbf^{4}\in C^{4}_{0}$ & $(1-r)_{+}^{5}(8+40r+48r^{2}+25r^{3}+5r^{4})$ \\
$\rbf^{5}\in C^{5}_{0}$ & $(1-r)_{+}^{6}(6+36r+82r^{2}+72r^{3}+30r^{4}+5r^{5})$ \\
\bottomrule
\end{tabular}
}
\caption{Wu functions $\psi^{q}\in C^{q}_{0}$, with compact support in $r\in (0,1)$ \cite{Fasshauer:2007:MAM:1506263}. Here $(\cdot)_{+} = \max(0,\cdot)$. The listed regularity excludes evaluation at the origin.}
\label{tab:Wu}
\end{table}

%%% Local Variables:
%%% mode: latex
%%% TeX-master: "manuscript.tex"
%%% End:

%% file: numericalexamples.tex
\subsection{Implementation}
\label{ss:resultsimpl}
To implement the software required to solve the advection-diffusion equation with the methods proposed in this paper we used MATLAB, except for the spectral Ewald method which is written in C\texttt{++}. There are also external libraries: for NUFFT we used the library from \cite{nufftlink}, based on \cite{LeslieGreengard2004AtNF}. The incomplete modified Bessel functions are not implemented in standard numerical libraries, but can be computed with the methods proposed in \cite{HarrisFrankE2009MfiB}; note that the precision is limited to ten digits. For the RBF-QR algorithm, used in PUX, we use the code from \cite{RBFQRlink}, as implemented in \cite{MR2801193}. To label points as inside or outside $\domain$ we use a spectral Ewald method for the Stokes potential \cite{LINDBO20118744}, as well as compute the velocity field in numerical example in section \ref{ss:resultsperiodic} and section \ref{ss:resultsdrop}.

Let $\mhu_{*}(\x_j)$ and $\adu_{*}(t,\x_{j})$ denote approximate solutions to  the modified Helmholtz equation \eqref{eq:mh}, \eqref{eq:mhbc}, and the advection-diffusion equation \eqref{eq:adu}--\eqref{eq:adbdry}, respectively, at the point $\x_j$ at the instance $t$, and let $\mhu$ and $\adu$ denote the corresponding solutions or a computed reference solutions. Let $e = \mhu - \mhu_{*}$ or $e = \adu - \adu_{*}$, then the $\ell_{2}$ error and the $\ell_{\infty}$ error are defined as
\begin{equation}
  \norm{e}_{2}= \ \sqrt{\frac{1}{N_{\domain}}\sum\limits_{j = 1}^{N_{\domain}}\abs{e(\x_{j})}^2}, \quad \quad\norm{e}_{\infty}= \max\limits_{j = 1,\ldots,N_{\domain}}\abs{e(\x_j)},
\end{equation}
at time $t$, over the points $\{\x_{j}\}_{j = 1}^{N_{\domain}}$ which are the points in $\X$ that are in $\domain(t)$.

\subsection{The modified Helmholtz equation}
\label{ss:resultsmodhelm}
To validate our proposed method for solving the modified Helmholtz equation \eqref{eq:mh}, \eqref{eq:mhbc} we consider the function
\begin{equation}
  \label{eq:example1u}
  \mhu(\x) = \sin(x_{1})\sin(x_{2})\,\exp{-(x_1^2 + x_2^2)/10}
\end{equation}
on the domain given by the parametrization
\begin{equation}
  \label{eq:starfish}
  \gamma(t) = (1 + 0.3\cos(5t))\,\exp{it} - 0.1045 + i5/439,\quad t\in[0,2\pi),
\end{equation}
with $a = 0.3$, $b = 5$,  and $c =  - 0.1045 + i5/439$.

If the Neumann boundary conditions and the forcing function are given by $\mhbc = \partial \mhu/\partial \hat{n}$ and $\mhrhs = \mhaa\mhu-\Delta\mhu$, respectively, then the exact solution is given by \eqref{eq:example1u}. 
\begin{figure}[htbp]
% Generated by: est_plots.jl
  \centering
  \includegraphics[width=0.49\textwidth,trim={5cm 9cm 5cm 7cm},clip]{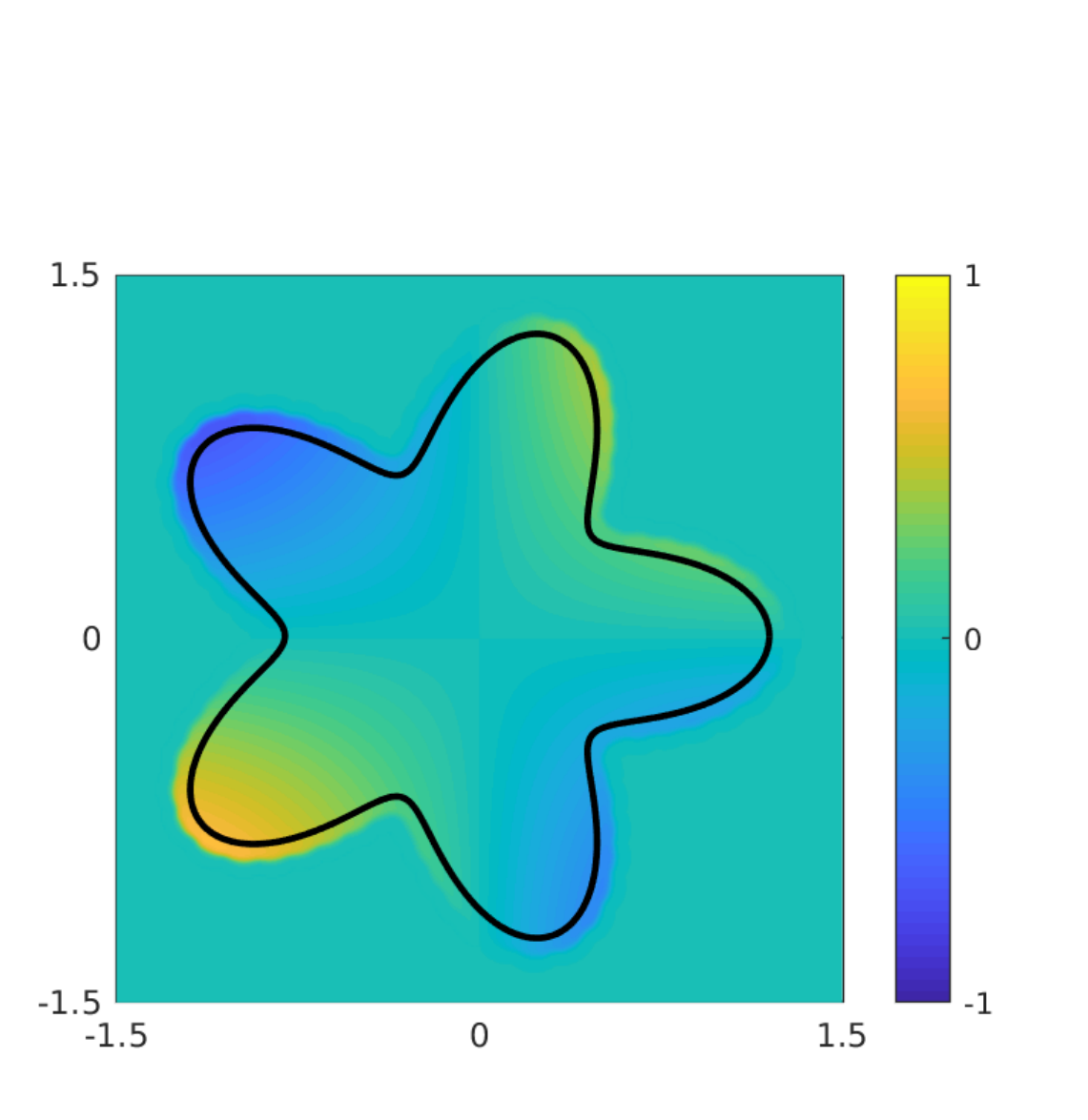}
  \hfill
    \includegraphics[width=0.48\textwidth,trim={5cm 9cm 5cm 7cm},clip]{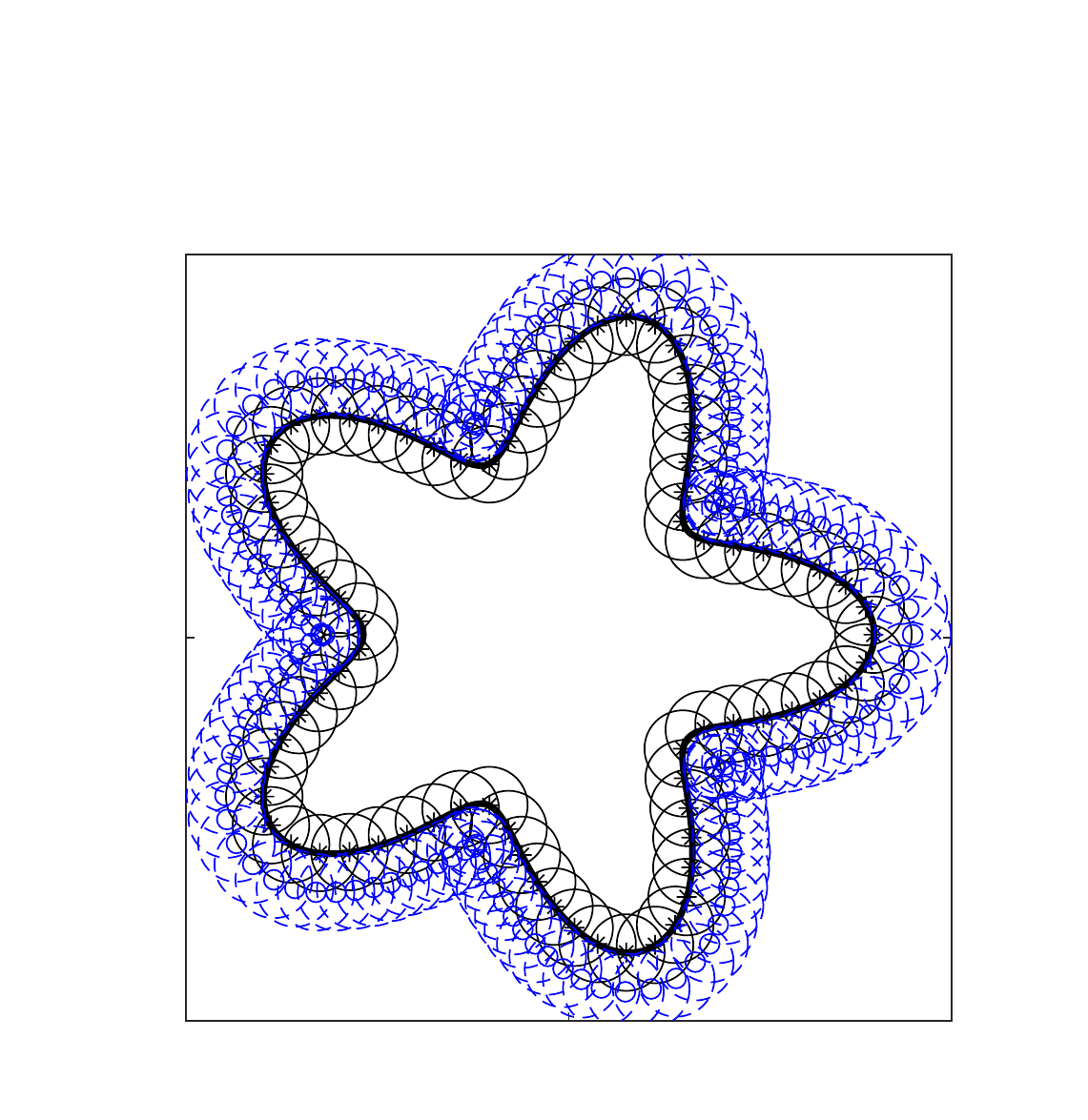}
    \caption{Left: The extension of function $\mhrhs = \mhaa\mhu-\Delta\mhu$ with $\mhu$ given by \eqref{eq:example1u}, with $\mhaa = 10$, constructed with PUX. Right:  The discs with black boundary are the extension partitions, while the zero partitions have a blue dashed boundary.}
  \label{fig:nemhext}
\end{figure}
For PUX we set the partition radius $R = 0.15$, which gives $77$ extension partitions, and $\varepsilon = 2$; we keep this value for $\varepsilon$ for all subsequent numerical experiments. An extension of $\mhrhs$ constructed with PUX, and a distribution of partitions along the boundary of the star-shaped domain, given by the parametrization \eqref{eq:starfish}, is shown in figure \ref{fig:nemhext}.

\begin{figure}[htbp]
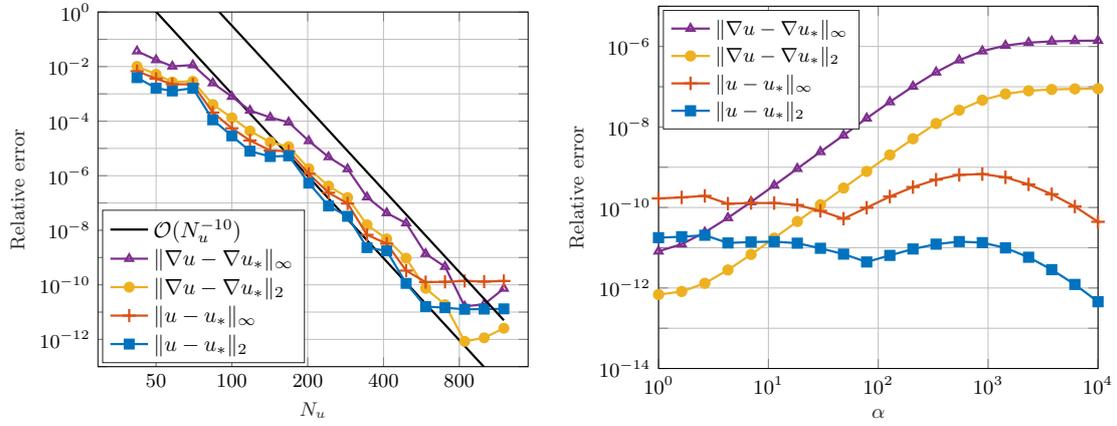

% Generated by: est_plots.jl
  \centering
  \includestandalone[width=0.47\textwidth]{mh_conv_a10}
  \hfill
    \includestandalone[width=0.50\textwidth]{mh_conv_nxe800}
  \caption{Errors in numerical solution for the modified Helmholtz equation with solution  \eqref{eq:example1u} on a starfish shaped domain. Left: Varying $\Nu$ for fixed $\mhaa = 10$. Right: Varying $\mha$ for fixed $\Nu = 800$.}
  \label{fig:nemherr}
\end{figure}

The boundary is discretized with $\Npan = 200$ panels, each with $\Nq = 16$ Gauss-Legendre nodes. For the first test we set $\mha = 10$ and solve with $\Nu$ ranging from $40$ to $1200$. The results are shown in figure \ref{fig:nemherr}. All errors decay as $\ordo{\Nu^{-10}}$, as expected \cite{FryklundFredrik2020Aien,FRYKLUNDPUX}. Also, we see that the errors level out at $10^{-10}$ and $10^{-11}$. This is due to the spectral Ewald method, where the involved algorithms for evaluating the incomplete modified Bessel functions are limited to about ten digits of accuracy \cite{HarrisFrankE2009MfiB}.

The next test is to fix $\Nu = 800$ and study the errors' dependence on $\mha$. The errors for $\mha$ ranging from $1$ to $10^{4}$ are shown in figure \ref{fig:nemherr}. For the relative $\ell_{\infty}$ error in $\mhu$ about one digit may be lost for $\mha \sim 10^{3}$, while for $\nabla\mhu$ both errors clearly deteriorate, as several orders of magnitudes in accuracy are lost. Recall that $\mha = \ordo{(\diffconst\dt)^{-1/2}}$, thus a high resolution in time results in more difficult problems in space. However, not until $\mha\approx 10^{2}$ the $\ell_{\infty}$ error is larger than $10^{-8}$, for this particular problem, which corresponds to time-steps of magnitude $10^{-4}$; the smallest value we consider for $\diffconst \dt$ is of the magnitude $10^{-4.5}$. From what we have observed this error does not limit the accuracy of  $\adu_{*}$ when solving the advection-diffusion equation.

\subsection{Convergence for rigid body in rotational flow}
\label{ss:rotconv}
\begin{figure}[htbp]
% Generated by: est_plots.jl
  \centering
  \includegraphics[width=0.49\textwidth,trim={0cm 0cm 0cm 0cm},clip]{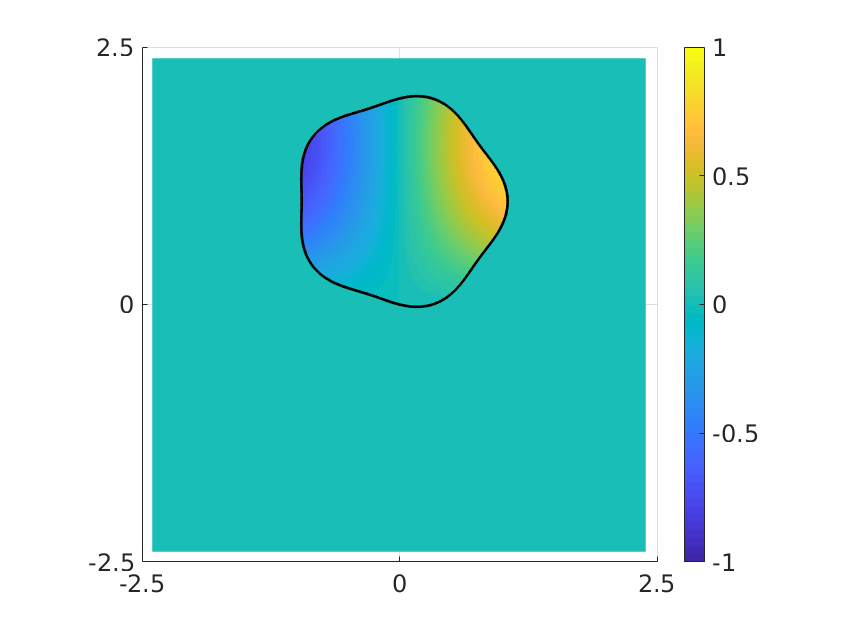}
  \hfill
    \includegraphics[width=0.48\textwidth,trim={0cm 0cm 0cm 0cm},clip]{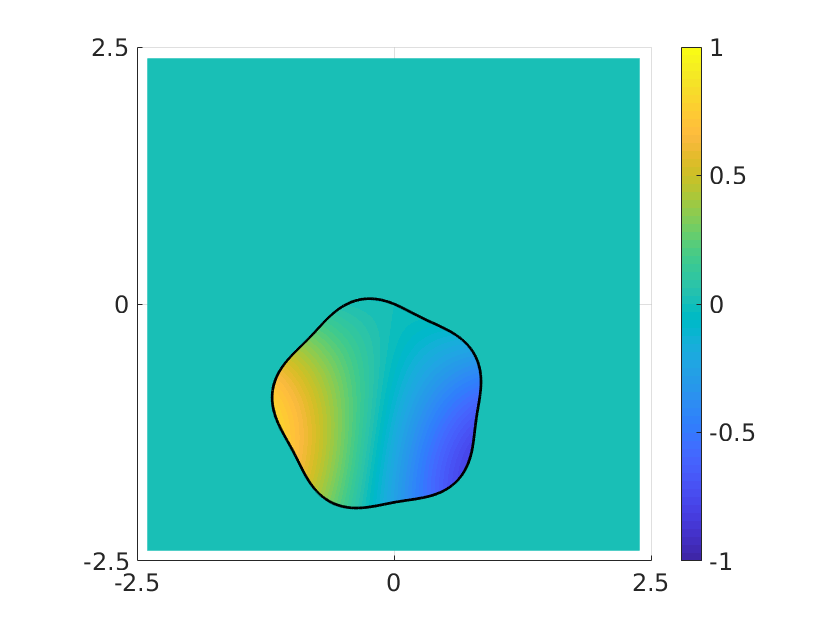}
    \caption{The function \eqref{eq:adurot} at the initial time $t = 0$ (left), and at the terminal time $t= 3$ (right), on a domain in a rotational flow \eqref{eq:rotfield}.}
  \label{fig:rotstarfish}
\end{figure}

We now solve the advection-diffusion equation \eqref{eq:adu}--\eqref{eq:adbdry} on a rigid body in a rotational flow, in order to show third and fourth order convergence in time. Consider
\begin{equation}
  \label{eq:adurot}
  \adu(t,\x) = \sin(\,\norm{\x}\cos(\tan^{-1}(x_{2}/x_{1}) -t))\sin(\,\norm{\x}\sin(\tan^{-1}(x_{2}/x_{1})-t)),
\end{equation}
and a velocity field $\velf$ given analytically by
\begin{equation}
  \label{eq:rotfield}
    \velf(\x) = w(-x_{2}, x_{1}),
\end{equation}
with $w=1$, which is a positively oriented rotational field around the origin. Let \eqref{eq:adurot} define the initial data $\adic$ at $t = 0$,  the Neumann boundary condition $\adbc$, and the source term $\force(t,\x) = -2\adu(t,\x)$, then \eqref{eq:adurot} is an analytical solution. The domain $\domain$ is given by the parametrization \eqref{eq:starfish}, with $a = 0.05$, $b = 5$, and $c = 0+i$, see figure \ref{fig:rotstarfish} for a visualization.

The boundary is discretized with $\Npan = 100$ panels, each with $\Nq = 16$ Gauss-Legendre nodes. The periodic box $\suppbox=[-\tfrac{L}{2},\tfrac{L}{2}]$, with $L = 4.8$, is discretized in each spatial dimension with $\Nu = 500$ uniformly distributed points; with these parameters the spatial error will not dominate the temporal error in the range of errors we investigate. The partition radius $R$ in PUX is set to $0.15$, giving $77$ extension partitions.

\begin{figure}[htbp]
%  \centering
\includestandalone[width=0.49\textwidth,clip]{errorconv}
  \hfill
\includegraphics[width=0.49\textwidth,trim={3cm 0cm 0cm 0cm},clip]{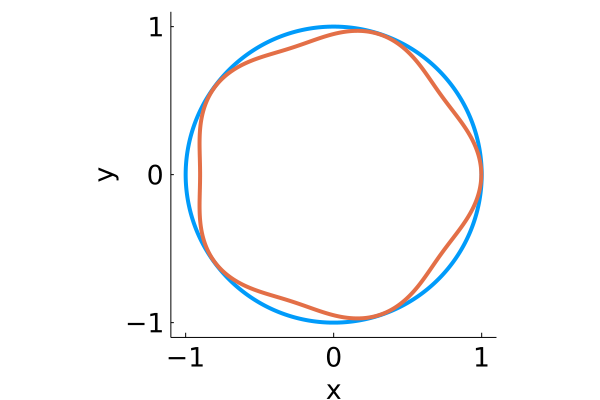}
  \caption{Left: Convergence plot of the absolute $\ell_{2}$ error at terminal time $t = 3$, as compared with
an analytical solution, for the numerical example described in section \ref{ss:rotconv}. Right: The domains used in section \ref{ss:cflcond} to study stability.}
  \label{fig:rotconv}
\end{figure}

We initialize the solver at $t = 0$, and measure the absolute $\ell_{2}$ error at the terminal time $t = 3$, see figure \ref{fig:rotstarfish} for a visualization of the solution at these instances in time. To show convergence we successively reduce the magnitude of time-step $\dt$. The results, using a third and fourth order SDC method, are shown in figure \ref{fig:rotconv}.  We observe the expected third order convergence. We also observe the expected fourth order convergence, but only between two data points. Larger time-steps, without changing $\Nu$, violates a CFL-type condition, as we observe that the solution becomes unstable. At the other end the error levels out at $10^{-10}$ as expected, due to limitations in accuracy in computing incomplete modified Bessel functions in the spectral Ewald method.

%%%%%%%%%%%%%%%%%%%%%%%%%%%%%%%%%%%%%%%%%%%%%%%%%%%%%%%%%%%%%%%%%%%%%%%%%%%%%%%%%
\subsection{Conservation of mass}
\label{ss:convmass}
We study the conservative properties of our solver on a time-dependent geometry $\domain(t)$ by computing the mass error  $e_{M}(t) = \abs{M(t) - M(\tO)}$ at time $t$, where
\begin{equation}
  \label{eq:masserr}
  M(t) = \int_{\domain(t)} \abs{\adu_{*}(t,\x)} \,\D \x.
\end{equation}
Since $\domain(t)$ cuts through the uniform grid, it is however difficult to numerically evaluate this integral in the general case, without introducing errors that obscure the errors in the mass conservation. Here, we do this for the special case where $\adu_{*}$ is the numerical solution to the advection-diffusion equation \eqref{eq:adu}--\eqref{eq:adbdry} with
\begin{equation}
  \label{eq:initdropsimu}
  \adu(t,\x) = \cos\left(k\sqrt{(x_{1}-c_{1}(t))^{2} + (x_{2}-c_{2}(t))^{2}}\right) + 1,
\end{equation}
where $k=2\pi$. Then the integral \eqref{eq:masserr} is identically equal to $2\pi$ for $\adu_{*} = \adu$ for all $t$. The rotational field is given by \eqref{eq:rotfield} with $w=1$, $\diffconst = 10^{-2}$ and the domain $\domain(t)$ is the unit circle centered at $c(t) = (-\sin(t), \cos(t))$.

The integral  \eqref{eq:masserr} is approximated by transforming to polar coordinates and then applying Gauss-Legendre quadrature in the radial direction and the trapezoidal rule in the angular direction, each with $150$ nodes. We run the simulation from $\tO = 0$ to the terminal time $t=3$.  The boundary is discretized with $\Npan = 100$ panels, each with $\Nq = 16$ Gauss-Legendre nodes.

We compute the mass error at $t = 3$ for a range of values for $N_{\dt} = t/\dt$ and $\Nu$, using SDC of order $3$ and order $4$. We also compute the max error $e_{\infty}$ at the terminal time by comparing $\adu_{*}$ with a reference solution using $\Nu = N_{\dt} = 800$ and SDC order $4$. The results are shown in table \ref{tab:masserr}. Clearly the mass error is at the level of the max error, or below. Furthermore, $\Nu = N_{\dt} = 100$ is clearly insufficient for SDC of order $3$ to yield useful results.

\begin{table}
\centering
{\def\arraystretch{1.3}
\begin{tabular}{c|c|c}
\toprule
\textbf{$e_{M}(t=3)$}  & \textbf{$\norm{e}_{\infty}$} & \textbf{$N_{\dt}$} = \textbf{$\Nu$}\\
\midrule
$3.7 \times 10^{-7}$ & $1.3 \times 10^{-5}$ & $400$\\
$1.4 \times 10^{-6}$ & $1.3 \times 10^{-4}$ & $200$\\
$9.3 \times 10^{-3}$ & $1.7 \times 10^{0}$ & $100$\\
%$3.9 \times 10^{-2}$ & $2.8 \times 10^{-3}$ & $50$\\  
\bottomrule
\end{tabular}
\quad
\begin{tabular}{c|c|c}
\toprule
\textbf{$e_{M}(t=3)$}  & \textbf{$\norm{e}_{\infty}$} & \textbf{$N_{\dt}$} = \textbf{$\Nu$}\\
\midrule
$1.3 \times 10^{-7}$ & $4.6 \times 10^{-7}$ & $400$\\
$2.0 \times 10^{-6}$ & $7.1 \times 10^{-5}$ & $200$\\
$6.3 \times 10^{-3}$ & $1.8 \times 10^{-2}$ & $100$\\
%$4.1 \times 10^{-1}$ & $9.9 \times 10^{-1}$ & $50$\\  
\bottomrule
\end{tabular}
}
\caption{Mass error and max error for the problem in section \ref{ss:convmass}. Left: SDC $3$. Right: SDC $4$.}
\label{tab:masserr}
\end{table}
%\begin{figure}[htbp]
% Generated by: est_plots.jl
%  \centering
%  \includegraphics[width=0.49\textwidth,trim={0cm 0cm 0cm 0cm},clip]{fig/cfl/circandstar.png}
%  \label{fig:cfldomain}
%\end{figure}

%%%%%%%%%%%%%%%%%%%%%%%%%%%%%%%%%%%%%%%%%%%%%%%%%%%%%%%%%%%%%%%%%%%%%%%%%%%%%%%%%
\subsection{Stability}
\label{ss:cflcond}
Since we are treating the diffusion term implicitly and the advection and source term explicitly, we expect to have a first order CFL condition of the form \eqref{eq:CFL}. For time-dependent geometries, we extend data into slices $S$. Here, we want to investigate if this extension affects the CFL condition. To do so, we solve the advection-diffusion equation on two different domains. Consider the unit circle centered at the origin with a star shaped domain inscribed in it. The star shaped domain is given by the parametrization \eqref{eq:starfish} with the same parameters as in section \ref{ss:rotconv}, except that the radius is set to one and it is centered at the origin, see figure \ref{fig:rotconv}. The rotational velocity field \eqref{eq:rotfield} is imposed, thus the slices $S$ are empty for the circle, but non-empty for the star shaped domain. 
We solve the advection-diffusion equation given in section \ref{ss:rotconv} for a variety of values of $w$, $\dt$, $\dx$ and $\diffconst$ to study if the solution becomes unsteady for the star shaped domain, but not for the circle. We set $L = 4.8$ and run the simulation to time $t = 3$. The boundaries are discretized with $100$ Gauss-Legendre panels, each with $16$ nodes.
 The results are summarized in table \ref{tab:cfl}, which shows that using PUX to accommodate for missing data, due to time-dependent geometry, does not impose a stricter CFL-type condition than the overall method. We leave further investigation of the CFL-type condition of the advection-diffusion equation solver to the future.

%Extending data into slices $S$, arising due to time-dependent geometry, might impose a CFL-type condition of the form \eqref{eq:CFL}. To investigate this, we solve the advection-diffusion equation on two different domains. Consider the unit circle centered at the origin with a star shaped domain inscribed in it. The star shaped domain is given by the parametrization \eqref{eq:starfish} with the same parameters as in section \ref{ss:rotconv}, see figure \ref{fig:rotstarfish}, except that the radius is set to one and it is centered at the origin. The rotational velocity field \eqref{eq:rotfield} is imposed, thus the slices $S$ are empty for the circle, but non-empty for the star shaped domain. 
%We solve the advection-diffusion equation given in section \ref{ss:rotconv} for a variety of values of $w$, $\dt$, $\dx$ and $\diffconst$ to study if the solution becomes unsteady for the star shaped domain, but not for the circle. We set $L = 4.8$ and run the simulation to time $t = 3$. The boundaries are discretized with $100$ Gauss-Legendre panels, each with $16$ nodes.
% The results are summarized in table \ref{tab:cfl}, which shows that using PUX to accommodate for missing data, due to time-dependent geometry, does not impose a stricter CFL-type condition than the overall method. We leave further investigation of the CFL-type condition of the advection-diffusion equation solver to the future.

\begin{table}
\centering
{\def\arraystretch{1.3}
\begin{tabular}{ccccccc}
\toprule
\textbf{Stable-circle} & \textbf{Stable-Star} & $w$ & $\dt$ & $\dx$ & $\max_{\x\in\bar{\domain}}\norm{\velf(\x)}\dt/\dx$& $\diffconst$\\
\midrule
\textcolor{green}{\checkmark} & \textcolor{green}{\checkmark} & $2$ & $6\times10^{-2}$ & $9.6\times10^{-2}$ & $1.25$ & $10^{-3}$\\
\textcolor{green}{\checkmark} & \textcolor{green}{\checkmark} & $4$ & $6\times10^{-3}$ & $2.4\times10^{-2}$ & $1$ & $10^{-3}$\\
\textcolor{red}{$\times$} & \textcolor{red}{$\times$} & $4$ & $3\times10^{-3}$ & $9.6\times10^{-3}$ & $1.25$ & $10^{-3}$\\
\textcolor{green}{\checkmark} & \textcolor{green}{\checkmark} & $2$ & $3\times10^{-3}$ & $9.6\times10^{-3}$ & $0.625$ & $10^{-3}$\\
\textcolor{red}{$\times$} & \textcolor{red}{$\times$} & $4$ & $6\times10^{-2}$ & $9.6\times10^{-2}$ & $2.5$ & $10^{-3}$\\
\textcolor{green}{\checkmark} & \textcolor{green}{\checkmark} & $4$ & $3\times10^{-3}$ & $9.6\times10^{-3}$ & $1.25$ & $10^{-2}$\\
\textcolor{green}{\checkmark} & \textcolor{green}{\checkmark} & $2$ & $6\times10^{-3}$ & $9.6\times10^{-2}$ & $0.125$ & $10^{-3}$\\
\textcolor{green}{\checkmark} & \textcolor{green}{\checkmark} & $2$ & $6\times10^{-3}$ & $9.6\times10^{-2}$ & $0.125$ & $10^{-3}$\\
\textcolor{green}{\checkmark} & \textcolor{green}{\checkmark} & $2$ & $6\times10^{-3}$ & $9.6\times10^{-2}$ & $0.125$ & $10^{-3}/2$\\
\textcolor{green}{\checkmark} & \textcolor{green}{\checkmark} & $2$ & $3\times10^{-3}$ & $9.6\times10^{-3}$ & $0.0625$ & $10^{-3}/4$\\
\textcolor{green}{\checkmark} & \textcolor{green}{\checkmark} & $4$ & $6\times10^{-3}$ & $9.6\times10^{-2}$ & $0.25$ & $10^{-3}$\\
\textcolor{red}{$\times$} & \textcolor{red}{$\times$} & $4$ & $6\times10^{-3}$ & $1.9\times10^{-2}$ & $1.25$ & $10^{-3}$\\
\textcolor{green}{\checkmark} & \textcolor{green}{\checkmark} & $4$ & $6\times10^{-3}$ & $2.4\times10^{-2}$ & $1$ & $10^{-3}$\\
\bottomrule
\end{tabular}
}
\caption{Results for study of stability in section \ref{ss:cflcond}. In the table \textcolor{green}{\checkmark} indicates a stable solution, while \textcolor{red}{$\times$} indicates an unstable solution.}
\label{tab:cfl}
\end{table}
%%%%%%%%%%%%%%%%%%%%%%%%%%%%%%%%%%%%%%%%%%%%%%%%%%%%%%%%%%%%%%%%%%%%%%%%%%%%%%%%%

\subsection{Adaptivity for rigid body in rotational flow}
\label{ss:rotadap}
We now solve the advection-diffusion equation in the setting of section \ref{ss:rotconv} with an adaptive time-stepper, as introduced in section \ref{sss:adaptivity}. We run the code for the tolerances $10^{-5}$, $10^{-6}$, $10^{-7}$ and $10^{-8}$ for time-step orders $K = 3$ and $K = 4$; for the tolerance $10^{-8}$ we set $\Nu = 500$, and $\Nu = 400$ for the other tolerances. The other parameters are set as in the previous section.

\begin{figure}[htbp]
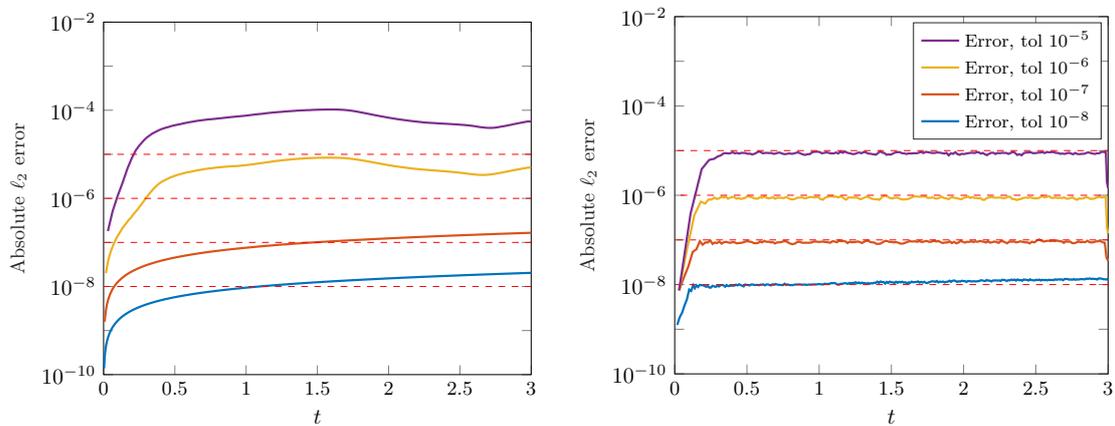

  \centering
      \includestandalone[width=0.49\textwidth]{adaptivesdc3}
  \hfill
  \includestandalone[width=0.49\textwidth,clip]{adaptivesdc4}
  \caption{The absolute $\ell_{2}$ errors plotted over time for solving the advection-diffusion equation for a rigid body in a rotational flow. The dashed red lines are different error tolerances. Left: Results for $K = 3$. Right: Results for $K = 4$.}
  \label{fig:adaptiveerror}
\end{figure}

 We measure the error by comparing with the analytical solution \eqref{eq:adurot} at each time-step, from $t_{0} = 0$, to the terminal time $t = 3$. The absolute $\ell_{2}$ errors over time are plotted in figure \ref{fig:adaptiveerror}. We see that SDC for both orders $K$ perform well, but $K=4$ gives results closer to the set tolerance. With a $K$th order method in total $N_{\text{mH}} = N_{\dt}K(K-1)$ modified Helmholtz equations are solved for $N_{\dt}$ time-steps, including the time-steps that were rejected due to not satisfying the set tolerance. A compilation of these values is presented in table \ref{tab:rotadap}. Clearly, $K = 3$ is preferable in terms of reducing $N_{\text{mH}}$ for the tested tolerances above $10^{-8}$. However, the errors are at some points about one order larger than the set tolerances, see figure \ref{fig:adaptiveerror}, thus smaller time-steps would be required to satisfy the tolerance, and consequently increasing $N_{\text{mH}}$. A safety factor in \eqref{eq:adaptiveScheme} less than $0.9$ is recommended if using SDC of order $K=3$ for less strict tolerances.

From studying convergence and adaptivity for a rigid body motion in section \ref{ss:rotconv} and section \ref{ss:rotadap} we conclude that, in terms of minimizing $N_{\mathrm{mH}}$, i.e. the total number of modified Helmholtz equations to solve, a fourth order SDC method is preferable over a third order SDC method only when considering tolerances stricter than $10^{-7}$.

\begin{table}
\centering
{\def\arraystretch{1.3}

\begin{tabular}{c|c|c|c|c}
\toprule
\textbf{Order $(K)$} & \textbf{Tolerance}  & \textbf{Maximum $\ell_{2}$ error} & \textbf{$N_{\dt}$} & \textbf{$N_{\mathrm{mH}}$}\\
\midrule
$4$ & $10^{-8}$ &$1.3 \times 10^{-8}$& $360$ & $4320$\\  
$3$ & $10^{-8}$ & $2.0\times 10^{-8}$  & $883$ & $5322$\\
$4$ & $10^{-7}$ &$1.0 \times 10^{-7}$& $212$ & $2544$\\
$3$ & $10^{-7}$ & $1.7\times 10^{-7}$ & $411$ & $2466$\\
$4$ & $10^{-6}$ &$9.7 \times 10^{-7}$& $209$ & $2508$\\  
$3$ & $10^{-6}$ & $8.4\times 10^{-6}$  & $223$ & $1338$\\
$4$ & $10^{-5}$ &$9.7 \times 10^{-6}$& $203$ & $2436$\\
$3$ & $10^{-5}$ & $1.0 \times 10^{-4} $  & $223$ & $1278$\\
\bottomrule
\end{tabular}
}
\caption{Here we report data from solving the advection-diffusion equation for a rigid body in a rotational flow, see section \ref{ss:rotadap}, with an adaptive time-stepper of orders $K$, and different tolerances. The reported error is the maximum over the entire simulation. In the table $N_{\dt}$ is the total number of time-steps, including the rejected time-steps, and $N_{\text{mH}} = N_{\dt}K(K-1)$ is the total number of modified Helmholtz equations solved.}
\label{tab:rotadap}
\end{table}

\iffalse
\begin{table}
\centering
{\def\arraystretch{1.3}

\begin{tabular}{c|c|c|c}
\toprule
\textbf{Order $(K)$} & \textbf{Tolerance}  & \textbf{$N_{\dt}$} & \textbf{$N_{\mathrm{mH}}$}\\
\midrule
  
$3$ & $10^{-8}$ & $883$ & $5322$\\
$3$ & $10^{-7}$ & $411$ & $2466$\\
$3$ & $10^{-6}$ & $223$ & $1338$\\
$3$ & $10^{-5}$ & $223$ & $1278$\\
\bottomrule
\end{tabular}
\quad
\begin{tabular}{c|c|c|c}
\toprule
\textbf{Order $(K)$} & \textbf{Tolerance}  & \textbf{$N_{\dt}$} & \textbf{$N_{\mathrm{mH}}$}\\
\midrule
$4$ & $10^{-8}$ & $360$ & $4320$\\
$4$ & $10^{-7}$ & $212$ & $2544$\\
$4$ & $10^{-6}$ & $209$ & $2508$\\
$4$ & $10^{-5}$ & $203$ & $2436$\\
\bottomrule
\end{tabular}
}
\caption{Here we report data from solving the advection-diffusion equation for a rigid body in a rotational flow, see section \ref{ss:rotadap}, with an adaptive time-stepper of orders $K$, and different tolerances. In the table $N_{\dt}$ is the total number of time-steps, including the rejected time-steps, and $N_{\text{mH}} = N_{\dt}K(K-1)$ is the total number of modified Helmholtz equations solved.}
\label{tab:rotadap}
\end{table}
\fi
\subsection{Periodic channel}
\label{ss:resultsperiodic}
\begin{figure}[htbp]
% Generated by: est_plots.jl
  \centering
    \includegraphics[width=0.49\textwidth,trim={0cm 0.3cm 0.5cm 0.3cm},clip]{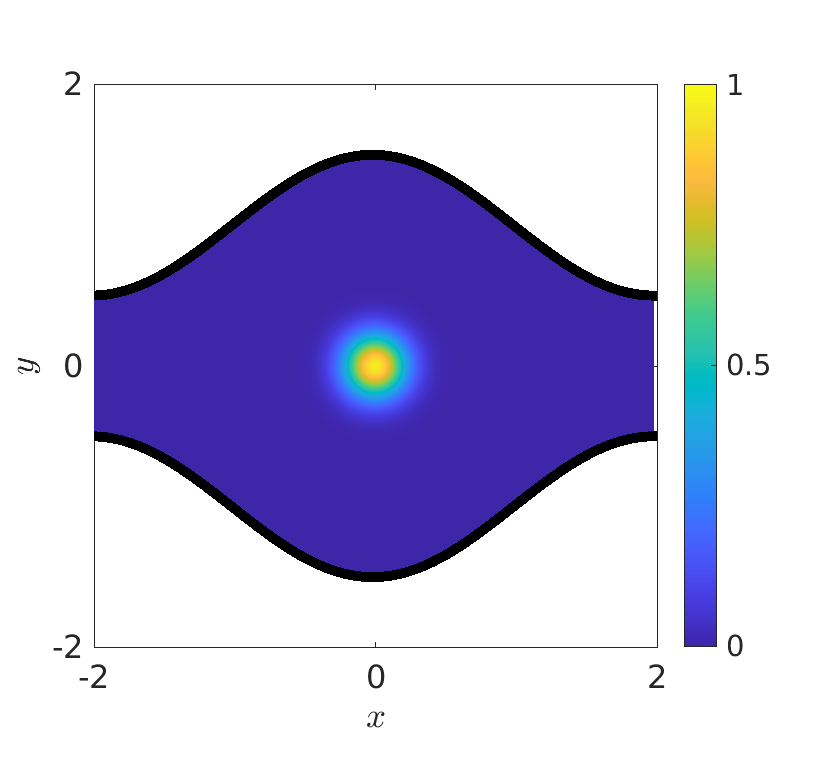}
    \hfill
      \includegraphics[width=0.455\textwidth,trim={0.3cm 0.1cm 0.3cm 0.5cm},clip]{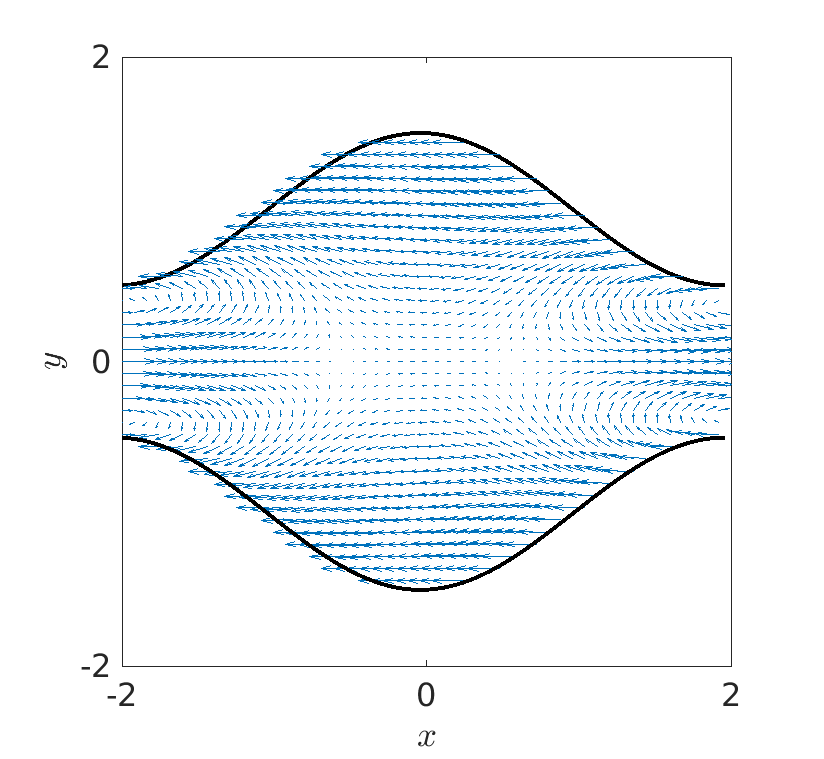}
  \caption{Left: The initial data \eqref{eq:icchannel} evaluated in a channel with sinusodial walls. Right: The velocity field $\velf$.}
  \label{fig:channelinit}
\end{figure}

We now consider the advection-diffusion equation with homogeneous Neumann boundary conditions in a singly periodic domain, with an adaptive time-stepper as described in section \ref{sss:adaptivity}. The domain $\domain$ is a channel with sinusodial walls in a box $\suppbox = [-L/2,L/2]$ where $L = 4$, that is periodic in the $x_{1}$-direction. The initial data $\adic$, with $\tO = 0$, is
\begin{equation}
  \label{eq:icchannel}
  \adic(\x) = \exp{-20x_{1}^{2}}\exp{-20x_{2}^{2}},
\end{equation}
which satisfies the imposed homogeneous Neumann boundary conditions, see figure \ref{fig:channelinit} for a visualization of the initial data and the geometry. The source term $\force$ is set to be identically equal to zero for all $t > 0$, and $\diffconst = 0.1$. The walls are translated with a constant velocity of magnitude $4$, from the right to the left. The velocity $\velf$ is found by solving the Stokes equations in the domain with an imposed pressure gradient of $(24,0)$, i.e. from the left to the right, and with no-slip boundary conditions. The Stokes equations are solved using a boundary integral equation \cite{BystrickyLukas2020Aaie}, using the same discretization of the walls as for the modified Helmholtz equation. Using the periodicity of the problem, the velocity $\velf\fp{\x}$ can be computed at every point in $\X\cap\domain(t)$ for any value for $t$ using post-processing. See figure \ref{fig:channelinit} for a visualization of the flow field.% For the CFL-type condition \eqref{eq:CFL} we have $\max_{\x\in\bar{\domain}}\norm{\velf(\x)} \approx 4.21$ for all $t$.
%The velocity field $\velf$ is a Stokes flow with a pressure gradient equal to $(24,0)$. The Stokes equations are solved once for the density in a boundary integral formulation, using the same discretization as for the modified Helmholtz equation. Then

% Old no trim

\begin{figure}[htbp]
% Generated by: est_plots.jl
  \centering
    \includegraphics[width=0.49\textwidth,trim={0cm 0cm 0cm 0cm},clip]{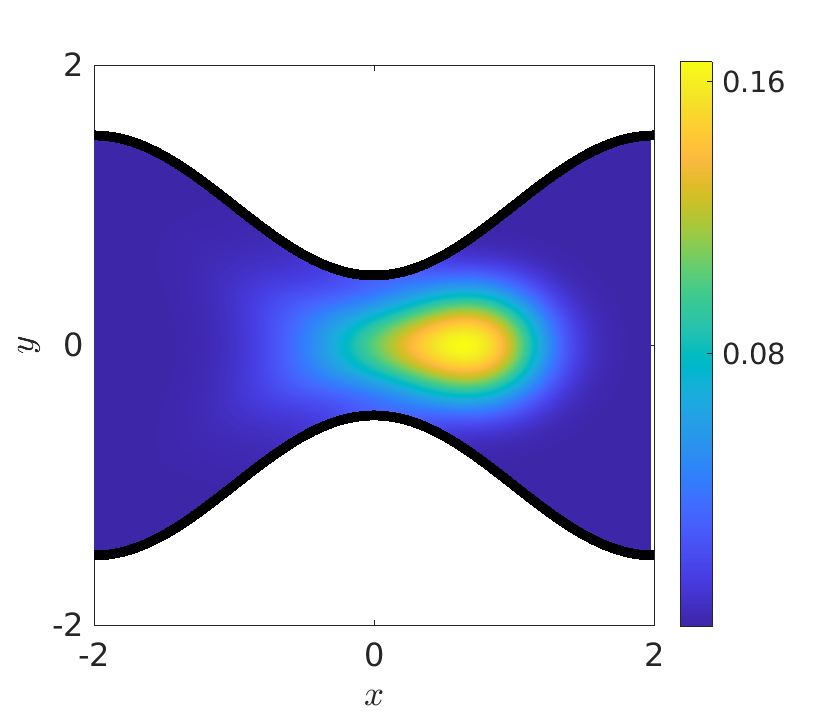}
    \hfill
      \includegraphics[width=0.49\textwidth,trim={0cm 0cm 0cm 0cm},clip]{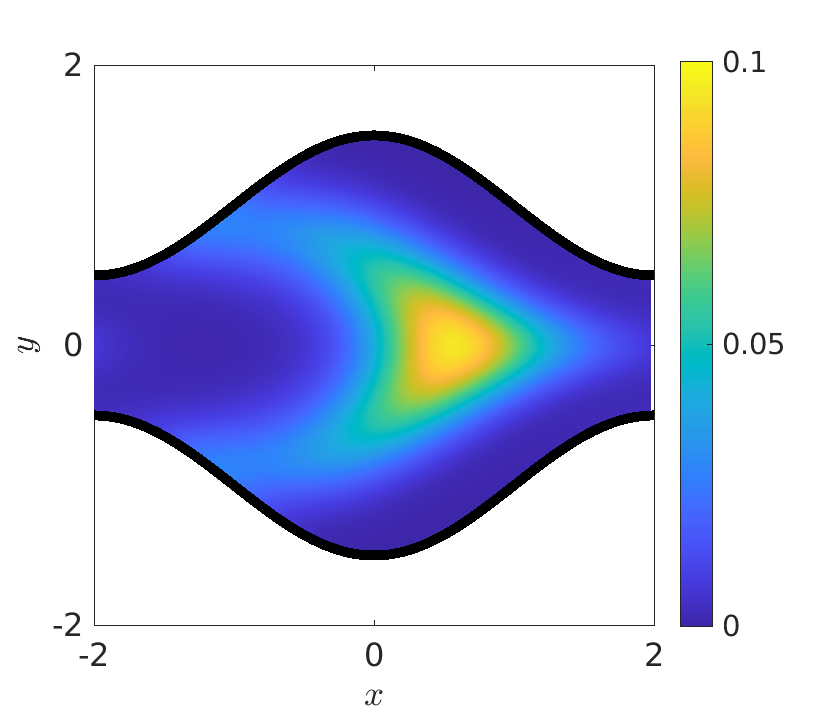}
  \caption{The numerical solution to the advection-diffusion equation, for the problem introduced in section \ref{ss:resultsperiodic}, evaluated at $t = 0.5$ (left) and the terminal time $t = 1$ (right).}
  \label{fig:channeltimestamps}
\end{figure}

We set $\Nu = 400$, and $\Npan = 100$ panels on each wall with $\Nq = 16$, and partition radius $R = 0.3$. A reference solution is computed with $\Nu = 400$ and a tolerance $10^{-10}$ for the estimated absolute $\ell_{2}$ error, with $K = 4$. The absolute $\ell_{2}$ error is measured at $t = 1$ by comparing the numerical solution with the reference solution. The results are presented in table \ref{tab:channel}. Two snapshots of the numerical solution at $t = 0.5$ and $t = 1$ are shown in figure \ref{fig:channeltimestamps}.

\begin{table}
\centering
{\def\arraystretch{1.3}
\begin{tabular}{c|c|c|c|c}
\toprule
\textbf{Order $(K)$} & \textbf{Tolerance} & \textbf{$\ell_{2}$ error} & \textbf{$N_{\dt}$} & \textbf{$N_{\mathrm{mH}}$}\\
\midrule
$4$ & $10^{-8}$ & $2.1\times 10^{-8}$ & $178$ & $2136$\\
$3$ & $10^{-8}$ & $2.7\times 10^{-8}$ & $561$ & $3366$\\
$4$ & $10^{-7}$ & $4.5\times 10^{-8}$ & $156$ & $1872$\\
$3$ & $10^{-7}$ & $3.0\times 10^{-7}$ & $255$ & $1530$\\
$4$ & $10^{-6}$ & $2.8\times 10^{-7}$ & $143$ & $1716$\\
$3$ & $10^{-6}$ & $1.3\times 10^{-5}$ & $157$ & $942$\\
$4$ & $10^{-5}$ & $3.4\times 10^{-6}$ & $136$ & $1632$\\
$3$ & $10^{-5}$ & $1.3\times 10^{-4}$ & $129$ & $774$\\
\bottomrule
\end{tabular}
}
\caption{Here we report data from solving the advection-diffusion equation in a periodic channel, see section \ref{ss:resultsperiodic}, with an adaptive time-stepper for different orders $K$ and tolerances. In the table the reported absolute $\ell_{2}$ error is measured at the terminal time $t = 1$, $N_{\dt}$ is the number of time-steps, including the rejected time-steps, and $N_{\mathrm{mh}} = N_{\dt}K(K-1)$ is the total number of modified Helmholtz equations solved.}
\label{tab:channel}
\end{table}

For $K = 4$ and tolerance $10^{-8}$ the absolute $\ell_{2}$ error is slightly above the set error tolerance; for the other tested less strict tolerances the $\ell_{2}$ error is well below the set tolerance. For $K = 3$ we get similar results for the tolerance $10^{-8}$, but by solving $1230$ more modified Helmholtz equations. As for the numerical example in section \ref{ss:rotadap}, the error is one order larger than the set tolerance for tolerances $10^{-5}$ and $10^{-6}$. The results are consistent, as we observe again that $K = 3$ is preferable, in terms of efficiency, over $K=4$ for less strict tolerances.

%Same settings for old

\subsection{Deforming drop}
\label{ss:resultsdrop}
We now consider the advection-diffusion equation on a time-dependent domain $\domain$ which deforms under the effect of a given velocity field, resembling a deforming drop. This is done for two different initial values and velocity fields: in the first setting we show convergence in time, and in the second we demonstrate robustness.

Consider the initial data \eqref{eq:icchannel} at $t = 0$ in a circle of radius $1.5$ centered at the origin, where $\partial \adic /\partial \normalscal = 0$ on the boundary. As in the previous example we impose homogeneous Neumann boundary conditions and set the source term $\force$ identically equal to zero. The diffusion coefficient $\diffconst$ is set to be equal to $0.1$.

\begin{figure}[htbp]
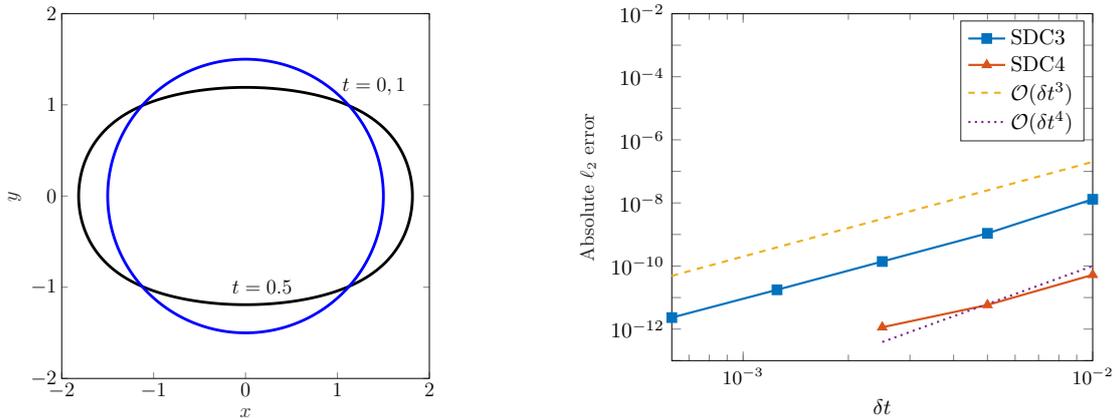

  \centering
\includestandalone[width=0.39\textwidth,clip]{circletimestamps}
  \hfill
      \includestandalone[width=0.49\textwidth]{errorconv_new}
  \caption{Left: The domain at initial time $t=0$ and $t = 1$ (blue), and domain at time $t = 0.5$ (black). Right: Convergence plot for relative max error for solving the advection-diffusion equation at terminal time $t = 1$, as compared with a reference solution. The magnitude of the solution is $\ordo{1}$.}
  \label{fig:deformdroperror}
\end{figure}

The time-dependent velocity field is given by
\begin{equation}
  \label{eq:extfield}
    \velf(\x) = (\cos(\pi t)\cos(x_{2})\sin(x_{1}),-\cos(\pi t)\sin(x_{2})\cos(x_{1})),
\end{equation}
which is an extensional flow, meaning the circle at time $t=0$ will become elongated and then go back to its initial state, see figure \ref{fig:deformdroperror}. After $t=0$ we have no parametrization to describe the boundary, thus we move the discretization of it according to the velocity field.

A reference solution is computed at time $t=1$ on a box with $L = 5$, $\Nu = 500$, and $\dt = 3.1250\times 10^{-4}$. The boundary is discretized with $\Npan = 100$ panels, each with $\Nq = 16$ Gauss-Legendre nodes. The partition radius in PUX is set to $0.4$, and the partitions are redistributed along the boundary as it changes over time.

The time-step $\dt$ is successively increased with a factor of two and we solve the advection-diffusion equation with SDC of order $3$ and $4$. The absolute $\ell_{2}$ error is computed at $t = 1$ by comparing with the reference solution. The results are shown in figure \ref{fig:deformdroperror}, where we clearly achieve third order convergence. Fourth order convergence is almost obtained; some other error or errors are interfering with the temporal convergence. For larger time-steps a CFL-type condition is violated, thus the convergence line cannot be continued. For smaller time-steps the error line would level out, again due to limitations in computing the incomplete modified Helmholtz equations. As seen before the third order SDC is more efficient for less strict tolerances.

We now study the robustness of the method simulating a drop with more complicated initial data and velocity field. Thus let
\begin{equation}
  \label{eq:initdropsimu}
  \adic(\x) = \cos(k\sqrt{(x_{1}-c_{1})^{2} + (x_{2}-c_{2})^{2}}) + 1,
\end{equation}
where $k = 2\pi/1.5$, $\bm{c} = (0,0.5)$, in a circle with radius $1.5$ centered at $\bm{c}$. The imposed velocity field is a superposition of the a rotational field \eqref{eq:rotfield} with $w=1$ and the extensional field \eqref{eq:extfield}.

We run the code for the time interval $[0,5]$ for homogeneous Neumann boundary conditions and a zero valued source term $\force$, with $\Nu = 300$, with the $L$ and boundary discretization and partition radius as above. The time-step is $0.01$.

\begin{figure}[htbp]
% Generated by: est_plots.jl
  \centering
  \includegraphics[width=0.49\textwidth,trim={3cm 7cm 3cm 7cm},clip]{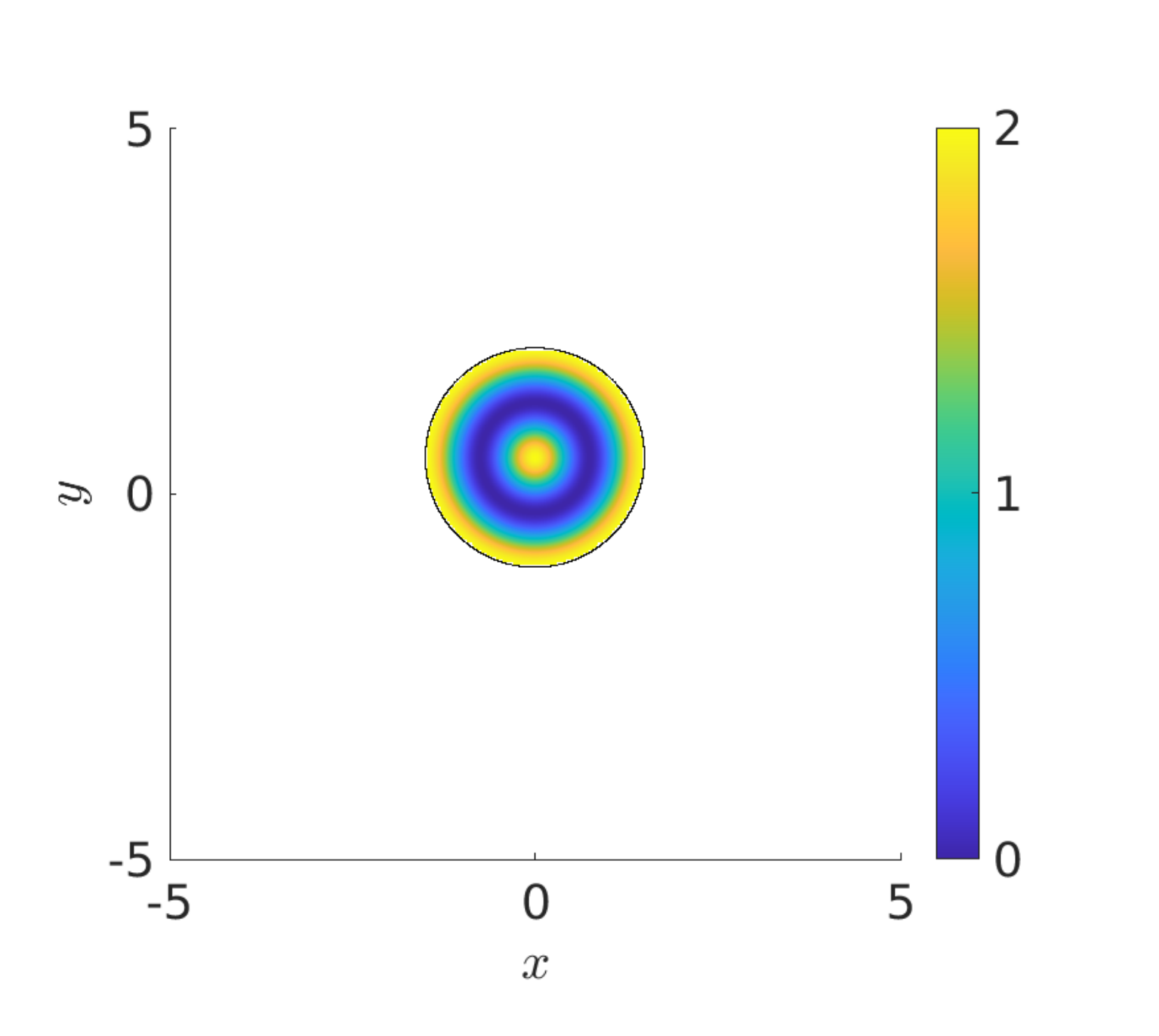}
    \hfill
      \includegraphics[width=0.49\textwidth,trim={3cm 7cm 3cm 7cm},clip]{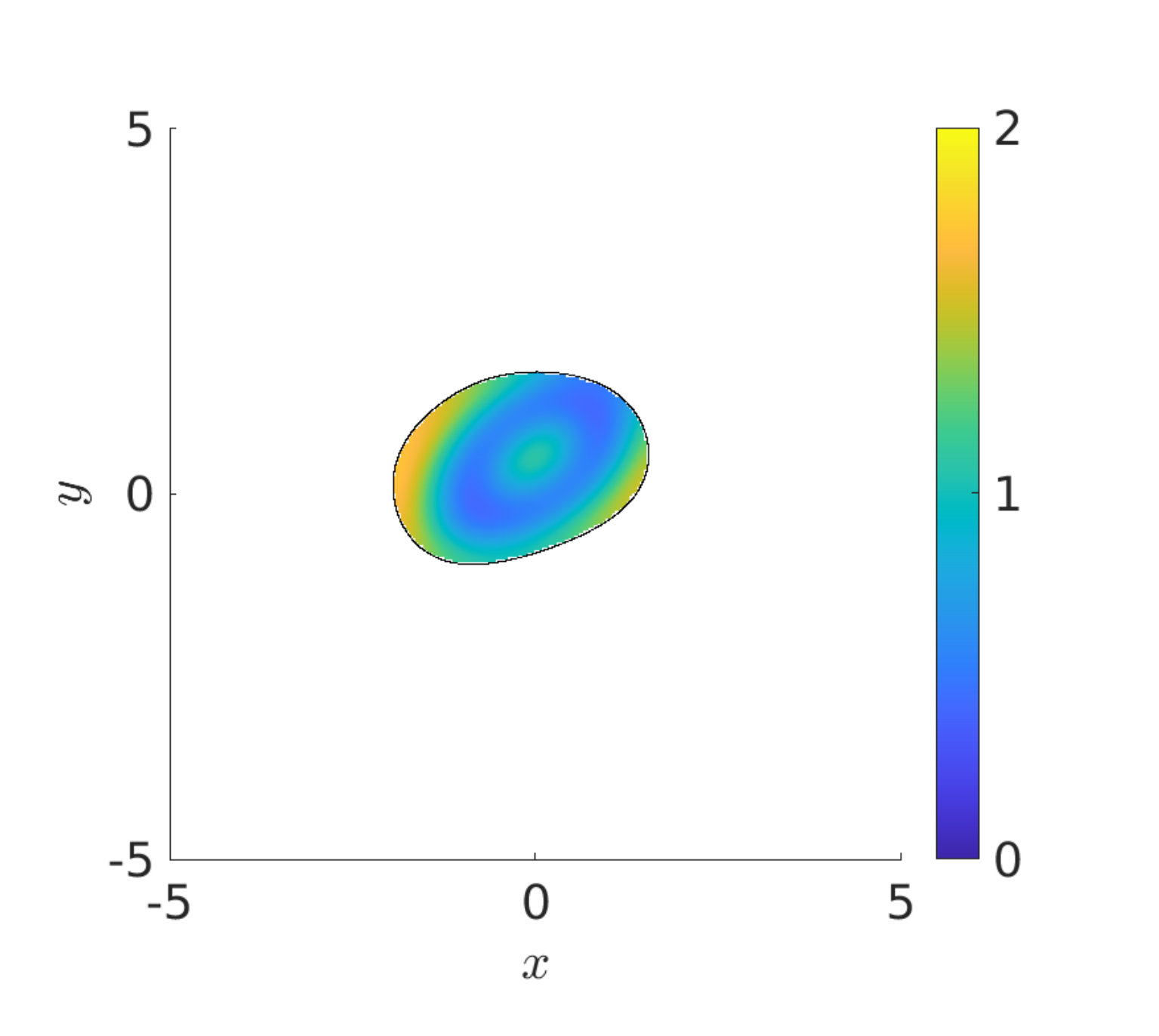}
%  \caption{Snapshots of simulation based on numerical solution of the advection-diffusion equation with initial data \eqref{eq:initdropsimu} at different times. Left: initial data at $t = 0$. Right: Numerical solution at $t = 0.5$.}
 % \label{fig:deformdrop1}
%\end{figure}

%\begin{figure}[htbp]
% Generated by: est_plots.jl
  \centering
  \includegraphics[width=0.49\textwidth,trim={3cm 7cm 3cm 7cm},clip]{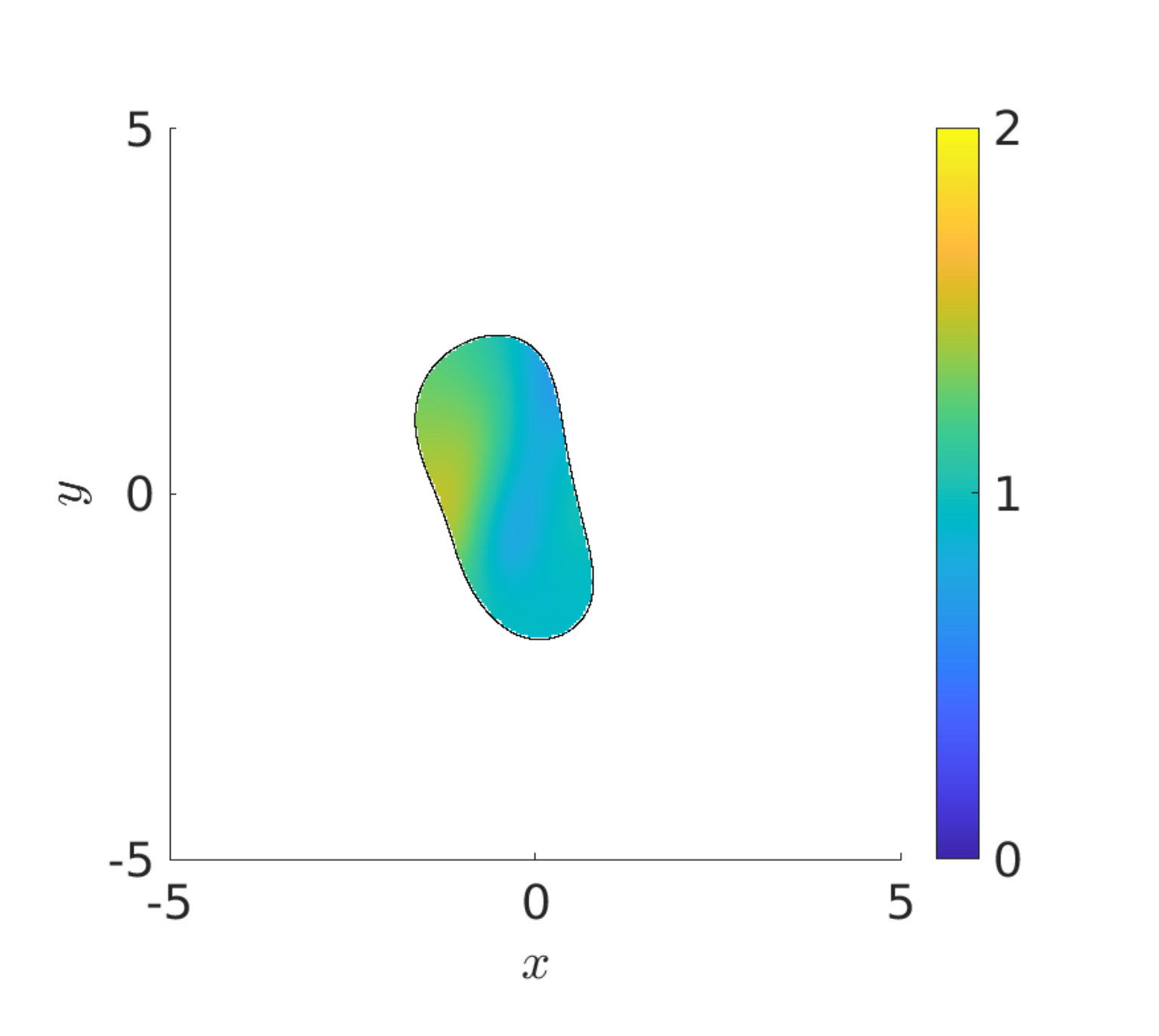}
    \hfill
  \includegraphics[width=0.49\textwidth,trim={3cm 7cm 3cm 7cm},clip]{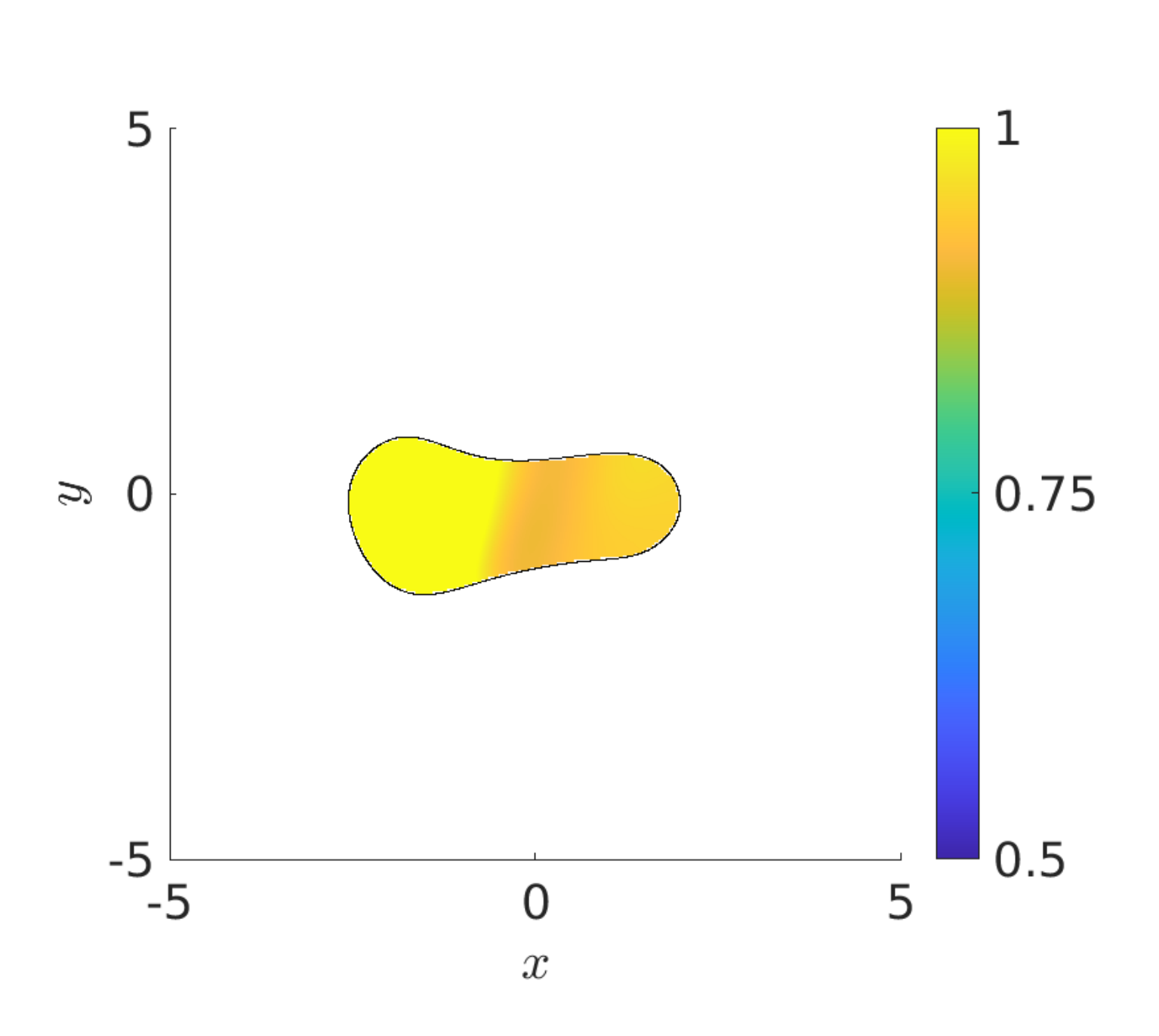}
%  \caption{Snapshots of simulation based on numerical solution of the advection-diffusion equation with initial data \eqref{eq:initdropsimu} at different times. Notice the different ranges for the colorbars. Left: Numerical solution at $t = 1.5$. Right: Numerical solution at $t = 2.5$.}
 % \label{fig:deformdrop2}
%\end{figure}

%\begin{figure}[htbp]
% Generated by: est_plots.jl
  \centering
  \includegraphics[width=0.49\textwidth,trim={3cm 7cm 3cm 7cm},clip]{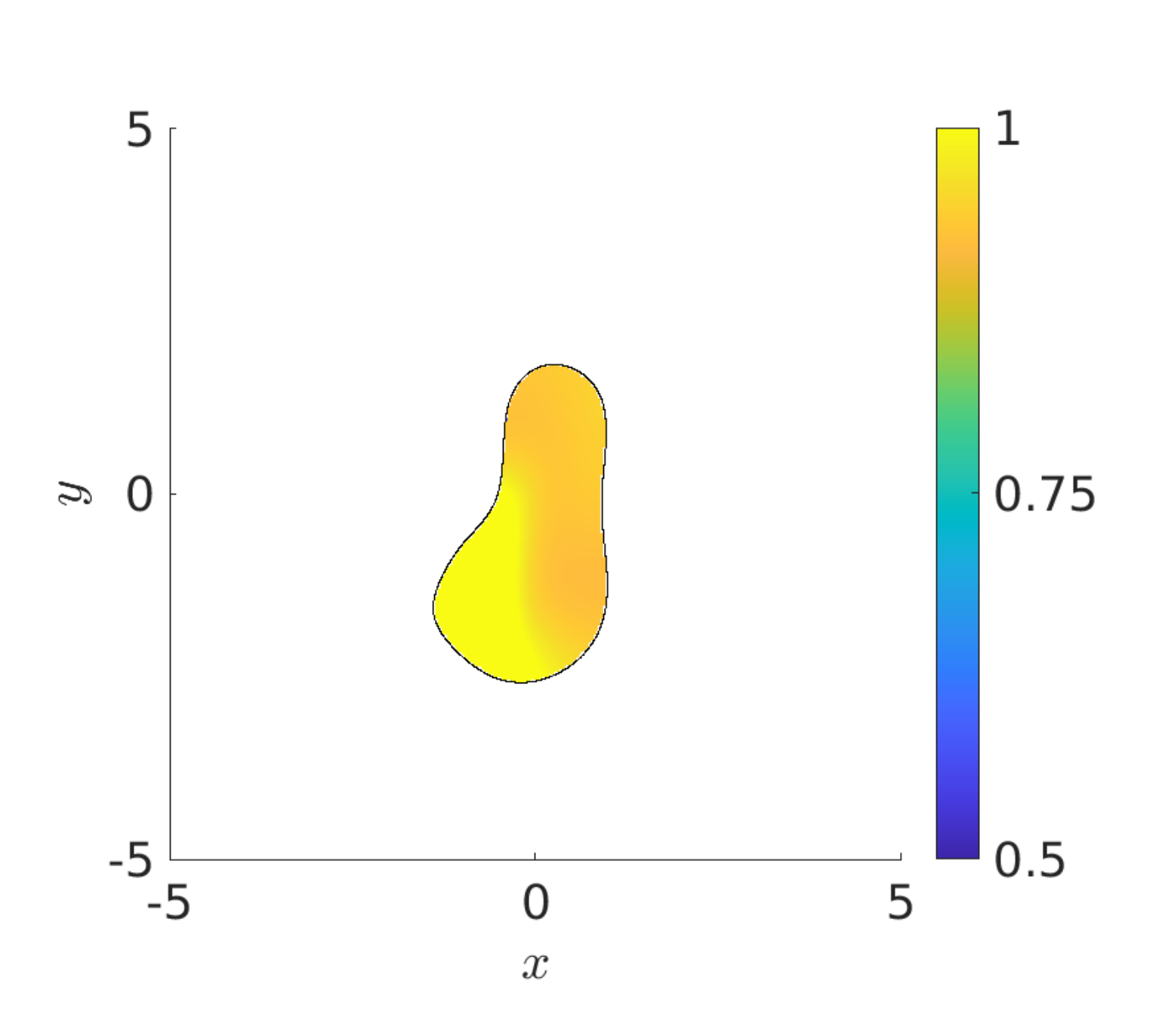}
    \hfill
      \includegraphics[width=0.49\textwidth,trim={3cm 7cm 3cm 7cm},clip]{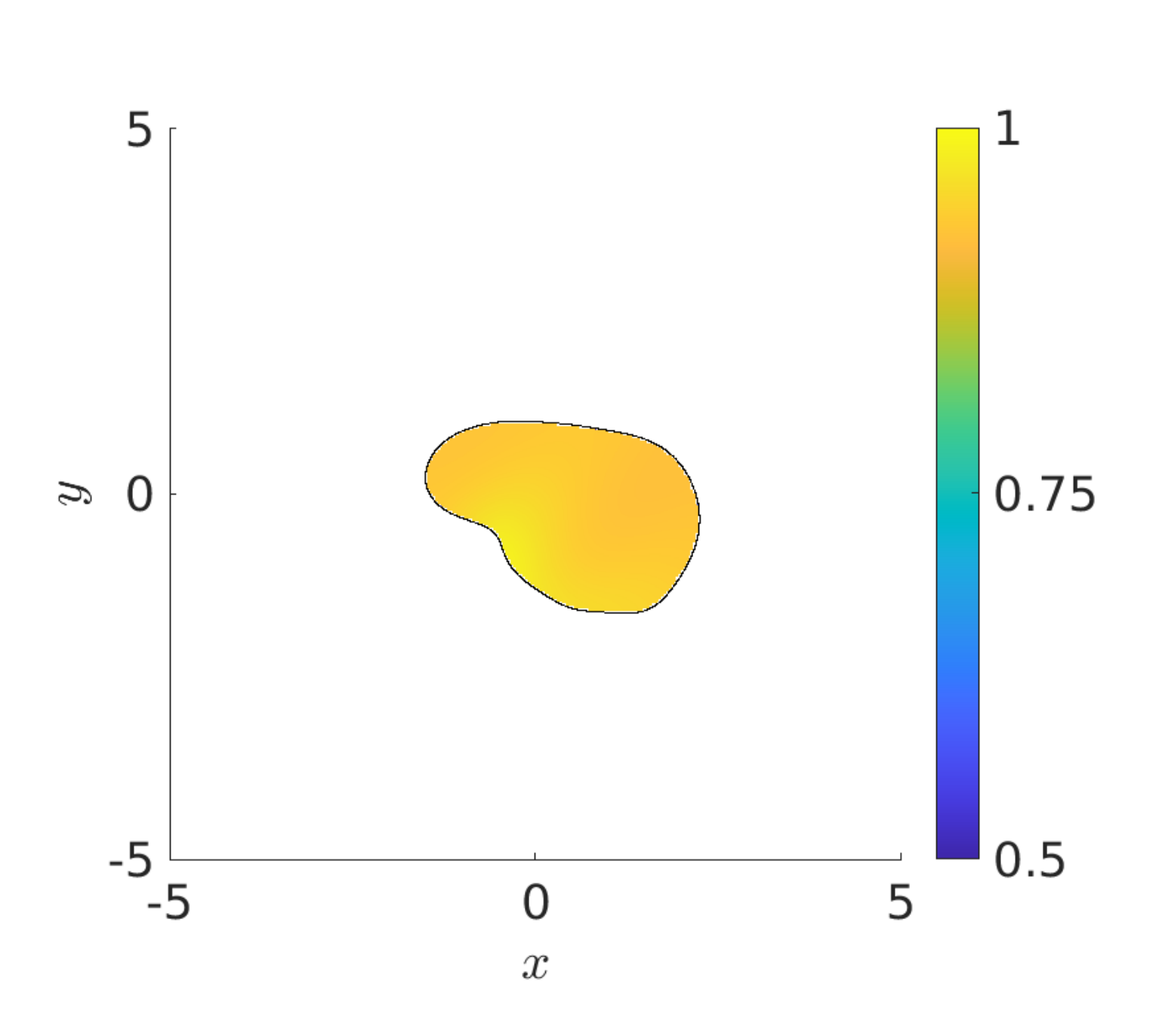}
  \caption{Snapshots of simulation based on numerical solution of the advection-diffusion equation with initial data \eqref{eq:initdropsimu} at different times. Notice the different ranges of the color bars. From left to right, top to bottom: $t = 0,\,0.5,\,1.5,\,2.5,\,3.5,\,4.5$.}
  \label{fig:deformdrop}
\end{figure}

The results are shown in figure \ref{fig:deformdrop}. Clearly the method is robust and can handle deformations into entirely new shapes. As expected the solution goes towards an equilibrium state due to diffusion.

%%% Local Variables:
%%% mode: latex
%%% TeX-master: "manuscript.tex"
%%% End:

%% file: mh_conv_a10.tex
\definecolor{mycolor1}{rgb}{0.00000,0.44700,0.74100}%
\definecolor{mycolor2}{rgb}{0.85000,0.32500,0.09800}%
\definecolor{mycolor3}{rgb}{0.92900,0.69400,0.12500}%
\definecolor{mycolor4}{rgb}{0.49400,0.18400,0.55600}%
\definecolor{mycolor5}{rgb}{0.46600,0.67400,0.18800}%
\definecolor{mycolor6}{rgb}{0.00000,0.74902,0.74902}%
\begin{tikzpicture}

\begin{loglogaxis}[%
xticklabel style={/pgf/number format/fixed},
xmin=0,
xmax=1500,
xlabel style={font=\small\color{white!15!black}},
xlabel={$N_{u}$},
xtick={50,100,200,400,800},
xticklabels={$50$,$100$,$200$,$400$,$800$},
xticklabel style = {font=\small},
xminorticks=false,
ymin=1e-13,
ymax=1,
yminorticks=true,
ylabel style={font=\small\color{white!15!black}},
ylabel={Relative error},
ytick={1,1e-1,1e-2,1e-3,1e-4,1e-5,1e-6,1e-7,1e-8,1e-9,1e-10,1e-11,1e-12,1e-13,1e-14},
yticklabels={$10^{0}$,$$,$10^{-2}$,$$,$10^{-4}$,$$,$10^{-6}$,$$,$10^{-8}$,$$,$10^{-10}$,$$,$10^{-12}$,$$,$10^{-14}$},
yticklabel style = {font=\small},
axis background/.style={fill=white},
xmajorgrids,
xminorgrids,
ymajorgrids,
legend style={font=\normalsize,at={(0.01,0.01)}, anchor=south west, legend cell align=left, align=left, draw=white!15!black},
clip mode=individual,
]
\addplot [color=black, line width=1.0pt, mark options={solid, black}]
  table[row sep=crcr]{%
42	5.85473302403221\\
50	1.024\\
58	0.232124009929528\\
70	0.0354013317464144\\
84	0.00571751271878145\\
100	0.001\\
118	0.000191064466913606\\
142	2.99996754374043e-05\\
168	5.58350851443501e-06\\
202	8.84069291692366e-07\\
242	1.45159792992328e-07\\
288	2.54727083052625e-08\\
344	4.3093877648689e-09\\
412	7.09623255631185e-10\\
492	1.20323088858919e-10\\
588	2.02407564521926e-11\\
702	3.44055796227802e-12\\
840	5.71751271878145e-13\\
1004	9.60866101010267e-14\\
1202	1.58838334029522e-14\\
};
\addlegendentry{$\mathcal{O}(\Nu^{-10})$}

\addplot [color=mycolor4, line width=1.0pt, mark=triangle, mark options={solid, mycolor4}, mark size=2.0pt]
  table[row sep=crcr]{%
42	0.0368031766567311\\
50	0.017520502478261\\
58	0.0102283653288061\\
70	0.0113667996780811\\
84	0.0024803455120461\\
100	0.00079179528778441\\
118	0.000245822920885157\\
142	0.000136865190792863\\
168	9.06470581020844e-05\\
202	1.88363089817118e-05\\
242	4.85530997265533e-06\\
288	1.76049247521564e-06\\
344	1.63903358230871e-07\\
412	4.30128462389829e-08\\
492	1.82147146004504e-08\\
588	1.37156297689305e-09\\
702	4.63363828890922e-10\\
840	1.62661005666507e-11\\
1004	1.93598347002956e-11\\
1202	7.32977822007216e-11\\
};
\addlegendentry{$\|\nabla\mhu-\nabla\mhu_{*}\|_{\infty}$}

\addplot [color=mycolor3, line width=1.0pt, mark=*, mark options={solid, mycolor3}, mark size=2.0pt]
  table[row sep=crcr]{%
42	0.0101649586986913\\
50	0.00522837045556098\\
58	0.00273130057824989\\
70	0.00289395299002193\\
84	0.000408237109580724\\
100	0.000132551068336037\\
118	4.36354126827915e-05\\
142	1.69987599983588e-05\\
168	1.15548985159513e-05\\
202	1.81155455342188e-06\\
242	4.19570768532061e-07\\
288	1.57809273241574e-07\\
344	1.58256586152373e-08\\
412	4.82839653029613e-09\\
492	9.44861665971199e-10\\
588	7.43880604371609e-11\\
702	1.87494209435432e-11\\
840	8.48226253812222e-13\\
1004	1.14006351949698e-12\\
1202	2.56060101746613e-12\\
};
\addlegendentry{$\|\nabla\mhu-\nabla\mhu_{*}\|_{2}$}

\addplot [color=mycolor2, line width=1.0pt, mark=+, mark options={solid, mycolor2}, mark size=3.0pt]
  table[row sep=crcr]{%
42	0.0067993361848937\\
50	0.00359381608992077\\
58	0.00224325177411346\\
70	0.00224591716810221\\
84	0.000206866339120952\\
100	5.46715731193013e-05\\
118	1.95677650373994e-05\\
142	8.38834874851399e-06\\
168	7.86915946614093e-06\\
202	1.17678831737978e-06\\
242	2.3768000880715e-07\\
288	9.60803705839045e-08\\
344	6.80364661736835e-09\\
412	3.23963216541968e-09\\
492	3.31592158134102e-10\\
588	1.25945152609982e-10\\
702	1.28388525443696e-10\\
840	1.38491722986268e-10\\
1004	1.31273558235686e-10\\
1202	1.38066110994753e-10\\
};
\addlegendentry{$\|\mhu-\mhu_{*}\|_{\infty}$}

\addplot [color=mycolor1, line width=1.0pt, mark=square*, mark options={solid, mycolor1}, mark size=2.0pt]
  table[row sep=crcr]{%
42	0.00393103910674211\\
50	0.00162394858282795\\
58	0.00127722816843265\\
70	0.00161512603201258\\
84	0.000111477315810139\\
100	2.88327281264096e-05\\
118	7.94394214575315e-06\\
142	5.03601872320747e-06\\
168	5.33885877112663e-06\\
202	5.34558200150361e-07\\
242	7.7672817627203e-08\\
288	3.24665654396773e-08\\
344	2.30824600480951e-09\\
412	1.76169423629494e-09\\
492	1.12840803931226e-10\\
588	1.59568922668496e-11\\
702	1.45019351657634e-11\\
840	1.26569056042293e-11\\
1004	1.28716206620156e-11\\
1202	1.31516710899709e-11\\
};
\addlegendentry{$\|\mhu-\mhu_{*}\|_{2}$}

\addplot [color=black, line width=1.0pt, mark options={solid, black}]
  table[row sep=crcr]{%
42	1851.42914481471\\
50	323.817232401242\\
58	73.4040570988849\\
70	11.1948840521896\\
84	1.80803627423312\\
100	0.316227766016838\\
118	0.0604198895372876\\
142	0.00948673034480058\\
168	0.00176566042405578\\
202	0.000279567257115965\\
242	4.59035570534305e-05\\
288	8.05517764177172e-06\\
344	1.36274806578479e-06\\
412	2.24402576841845e-07\\
492	3.80495015901013e-08\\
588	6.40068919536775e-09\\
702	1.08799995826262e-09\\
840	1.80803627423312e-10\\
1004	3.03852540563786e-11\\
1202	5.0229091527992e-12\\
};

\end{loglogaxis}
\end{tikzpicture}%

%%% Local Variables:
%%% mode: latex
%%% TeX-master: t
%%% End:

%% file: mh_conv_nxe800.tex
\definecolor{mycolor1}{rgb}{0.00000,0.44700,0.74100}%
\definecolor{mycolor2}{rgb}{0.85000,0.32500,0.09800}%
\definecolor{mycolor3}{rgb}{0.92900,0.69400,0.12500}%
\definecolor{mycolor4}{rgb}{0.49400,0.18400,0.55600}%
\definecolor{mycolor5}{rgb}{0.46600,0.67400,0.18800}%
\definecolor{mycolor6}{rgb}{0.00000,0.74902,0.74902}%
\begin{tikzpicture}

\begin{loglogaxis}[%
xticklabel style={/pgf/number format/fixed},
xmin=1,
xmax=10000,
xminorticks=false,
xlabel style={font=\small\color{white!15!black}},
xlabel={$\mha$},
xtick={1,10,100,1000,10000},
xticklabels={$10^{0}$,$10^{1}$,$10^{2}$,$10^{3}$,$10^{4}$},
xticklabel style = {font=\small},
ymin=1e-14,
ymax=1e-5,
yminorticks=true,
ylabel style={font=\small\color{white!15!black}},
ylabel={Relative error},
ytick={1e-1,1e-2,1e-3,1e-4,1e-5,1e-6,1e-7,1e-8,1e-9,1e-10,1e-11,1e-12,1e-13,1e-14},
yticklabels={$$,$10^{-2}$,$$,$10^{-4}$,$$,$10^{-6}$,$$,$10^{-8}$,$$,$10^{-10}$,$$,$10^{-12}$,$$,$10^{-14}$},
yticklabel style = {font=\small},
axis background/.style={fill=white},
xmajorgrids,
%xminorgrids,
ymajorgrids,
legend style={font=\small,at={(0.01,0.66)}, anchor=south west, legend cell align=left, align=left, draw=white!15!black},
clip mode=individual,
]
\addplot [color=mycolor4, line width=1.0pt, mark=triangle, mark options={solid, mycolor4}, mark size=2.0pt]
  table[row sep=crcr]{%
10000	1.39715379960729e-06\\
6158.48211066027	1.38085425345185e-06\\
3792.69019073225	1.33968908144118e-06\\
2335.72146909012	1.24325385069567e-06\\
1438.44988828766	1.0532574015969e-06\\
885.866790410083	7.67191090418505e-07\\
545.559478116852	4.57577070667081e-07\\
335.981828628378	2.29164569439598e-07\\
206.913808111479	1.01582626917091e-07\\
127.427498570313	4.18054302708088e-08\\
78.4759970351462	1.64567066666449e-08\\
48.3293023857175	6.20396554816113e-09\\
29.7635144163132	2.41619228794732e-09\\
18.3298071083244	9.22043362634684e-10\\
11.2883789168469	3.53561025232532e-10\\
6.95192796177561	1.37459917951016e-10\\
4.2813323987194	5.52932647442095e-11\\
2.63665089873036	2.41235446635515e-11\\
1.62377673918872	1.22796913178978e-11\\
1	8.17203862345907e-12\\
};
\addlegendentry{$\|\nabla\mhu-\nabla\mhu_{*}\|_{\infty}$}

\addplot [color=mycolor3, line width=1.0pt, mark=*, mark options={solid, mycolor3}, mark size=2.0pt]
  table[row sep=crcr]{%
10000	8.95631369468293e-08\\
6158.48211066027	8.84192066915416e-08\\
3792.69019073225	8.56340027328454e-08\\
2335.72146909012	7.91966742951452e-08\\
1438.44988828766	6.62564614774238e-08\\
885.866790410083	4.64805392173769e-08\\
545.559478116852	2.61217696409936e-08\\
335.981828628378	1.21906441509268e-08\\
206.913808111479	5.0973097711969e-09\\
127.427498570313	2.02501124509972e-09\\
78.4759970351462	7.88228281875515e-10\\
48.3293023857175	3.03985076340898e-10\\
29.7635144163132	1.16919956511787e-10\\
18.3298071083244	4.49529336125985e-11\\
11.2883789168469	1.73693833114493e-11\\
6.95192796177561	6.82529813408723e-12\\
4.2813323987194	2.80735343894925e-12\\
2.63665089873036	1.30961558110718e-12\\
1.62377673918872	8.17787639346156e-13\\
1	6.86895869500738e-13\\
};
\addlegendentry{$\|\nabla\mhu-\nabla\mhu_{*}\|_{2}$}

\addplot [color=mycolor2, line width=1.0pt, mark=+, mark options={solid, mycolor2}, mark size=3.0pt]
  table[row sep=crcr]{%
10000	4.41077029835337e-11\\
6158.48211066027	1.05688846556171e-10\\
3792.69019073225	2.13591676640896e-10\\
2335.72146909012	3.74982904988165e-10\\
1438.44988828766	5.61584881358252e-10\\
885.866790410083	6.7624094295088e-10\\
545.559478116852	6.44882375708793e-10\\
335.981828628378	4.91396500005732e-10\\
206.913808111479	3.21383710591559e-10\\
127.427498570313	1.87128147827755e-10\\
78.4759970351462	9.87569044979006e-11\\
48.3293023857175	5.26807852328286e-11\\
29.7635144163132	8.24910271837869e-11\\
18.3298071083244	1.1505874428035e-10\\
11.2883789168469	1.28951539400846e-10\\
6.95192796177561	1.29234126835855e-10\\
4.2813323987194	1.24105949365862e-10\\
2.63665089873036	1.94131926679392e-10\\
1.62377673918872	1.78256436061988e-10\\
1	1.67774081688638e-10\\
};
\addlegendentry{$\|\mhu-\mhu_{*}\|_{\infty}$}

\addplot [color=mycolor1, line width=1.0pt, mark=square*, mark options={solid, mycolor1}, mark size=2.0pt]
  table[row sep=crcr]{%
10000	4.55577995648133e-13\\
6158.48211066027	1.22500703528925e-12\\
3792.69019073225	2.85711778154624e-12\\
2335.72146909012	5.80492328899908e-12\\
1438.44988828766	9.9203860791008e-12\\
885.866790410083	1.33792981812824e-11\\
545.559478116852	1.4213790923994e-11\\
335.981828628378	1.22766467514895e-11\\
206.913808111479	9.3190500695509e-12\\
127.427498570313	6.45556389320688e-12\\
78.4759970351462	4.43695943422927e-12\\
48.3293023857175	6.94126744316322e-12\\
29.7635144163132	9.70158485853709e-12\\
18.3298071083244	1.29784122299417e-11\\
11.2883789168469	1.4178353290888e-11\\
6.95192796177561	1.39024494365861e-11\\
4.2813323987194	1.3145506953988e-11\\
2.63665089873036	2.04451592715527e-11\\
1.62377673918872	1.87933416039565e-11\\
1	1.77849937566641e-11\\
};
\addlegendentry{$\|\mhu-\mhu_{*}\|_{2}$}

\end{loglogaxis}
\end{tikzpicture}%

%%% Local Variables:
%%% mode: latex
%%% TeX-master: t
%%% End:

%% file: errorconv.tex
\definecolor{mycolor1}{rgb}{0.00000,0.44700,0.74100}%
\definecolor{mycolor2}{rgb}{0.85000,0.32500,0.09800}%
\definecolor{mycolor3}{rgb}{0.92900,0.69400,0.12500}%
\definecolor{mycolor4}{rgb}{0.49400,0.18400,0.55600}%
\begin{tikzpicture}

\begin{loglogaxis}[%
xmin=0.001,
xmax=0.01,
xminorticks=true,
xticklabel style = {font=\small},
xlabel={$\delta t$},
ymode=log,
ymin=1e-13,
ymax=1e-4,
yminorticks=true,
ylabel style={font=\small\color{white!15!black}},
ylabel={Absolute $\ell_{2}$ error},
ytick={1,1e-1,1e-2,1e-3,1e-4,1e-5,1e-6,1e-7,1e-8,1e-9,1e-10,1e-11,1e-12,1e-13,1e-14},
yticklabels={$10^{0}$,$$,$10^{-2}$,$$,$10^{-4}$,$$,$10^{-6}$,$$,$10^{-8}$,$$,$10^{-10}$,$$,$10^{-12}$,$$,$10^{-14}$},
axis background/.style={fill=white},
legend style={legend cell align=left, align=left, align=left,at={(0.03,0.80)},anchor=west, draw=white!15!black}
]
\addplot [color=mycolor1, line width=1.0pt, mark=square*, mark options={solid, mycolor1}]
  table[row sep=crcr]{%
0.001875	4.22617104936602e-09\\
0.00375	2.65914899832274e-08\\
0.0075	2.01597591727525e-07\\
};
\addlegendentry{SDC$3$}

\addplot [color=mycolor2, line width=1.0pt, mark=triangle*, mark options={solid, mycolor2}]
  table[row sep=crcr]{%
0.001875	1.83737972609818e-10\\
0.00375	7.95971890167656e-10\\
0.0075	1.18791211744236e-08\\
};
\addlegendentry{SDC$4$}

\addplot [color=mycolor3, dashed, line width=1.0pt]
  table[row sep=crcr]{%
0.001875	1.318359375e-09\\
0.00375	1.0546875e-08\\
0.0075	8.4375e-08\\
};
\addlegendentry{$\mathcal{O}(\delta t^{3})$}

\addplot [color=mycolor4, dotted, line width=1.0pt]
  table[row sep=crcr]{%
0.001875	1.95886762373886e-11\\
0.00375	3.13418819798218e-10\\
0.0075	5.01470111677149e-09\\
};
\addlegendentry{$\mathcal{O}(\delta t^4)$}

\end{loglogaxis}

\end{tikzpicture}%

%%% Local Variables:
%%% mode: latex
%%% TeX-master: "../../manuscript"
%%% End:

%% file: errorconv_new.tex
\definecolor{mycolor1}{rgb}{0.00000,0.44700,0.74100}%
\definecolor{mycolor2}{rgb}{0.85000,0.32500,0.09800}%
\definecolor{mycolor3}{rgb}{0.92900,0.69400,0.12500}%
\definecolor{mycolor4}{rgb}{0.49400,0.18400,0.55600}%
\begin{tikzpicture}

\begin{loglogaxis}[%
xmin=0.000625,
xmax=0.01,
xminorticks=true,
xticklabel style = {font=\small},
xlabel={$\delta t$},
ymode=log,
ymin=1e-13,
ymax=1e-2,
yminorticks=true,
ylabel style={font=\small\color{white!15!black}},
ylabel={Absolute $\ell_{2}$ error},
ytick={1,1e-1,1e-2,1e-3,1e-4,1e-5,1e-6,1e-7,1e-8,1e-9,1e-10,1e-11,1e-12,1e-13,1e-14},
yticklabels={$10^{0}$,$$,$10^{-2}$,$$,$10^{-4}$,$$,$10^{-6}$,$$,$10^{-8}$,$$,$10^{-10}$,$$,$10^{-12}$,$$,$10^{-14}$},
axis background/.style={fill=white},
legend style={legend cell align=left, align=left, draw=white!15!black}
]
\addplot [color=mycolor1, line width=1.0pt, mark=square*, mark options={solid, mycolor1}]
  table[row sep=crcr]{%
0.000625	2.31948069224286e-12\\
0.00125	1.75978870041848e-11\\
0.0025	1.38789751724218e-10\\
0.005	1.08886408602556e-09\\
0.01	1.30471731194644e-08\\
};
\addlegendentry{SDC$3$}

\addplot [color=mycolor2, line width=1.0pt, mark=triangle*, mark options={solid, mycolor2}]
  table[row sep=crcr]{%
0.0025	1.14492683750893e-12\\
0.005	5.81695592130002e-12\\
0.01	5.25114306711836e-11\\
};
\addlegendentry{SDC$4$}

\addplot [color=mycolor3, dashed, line width=1.0pt]
  table[row sep=crcr]{%
0.000625	4.8828125e-11\\
0.00125	3.90625e-10\\
0.0025	3.125e-09\\
0.005	2.5e-08\\
0.01	2e-07\\
};
\addlegendentry{$\mathcal{O}(\delta t^{3})$}

\addplot [color=mycolor4, dotted, line width=1.0pt]
  table[row sep=crcr]{%
0.0025	3.90625e-13\\
0.005	6.25e-12\\
0.01	1e-10\\
};
\addlegendentry{$\mathcal{O}(\delta t^{4})$}

\end{loglogaxis}

\end{tikzpicture}%
%%% Local Variables:
%%% mode: latex
%%% TeX-master: "../../manuscript"
%%% End:

%% file: conclusions.tex
We have presented a boundary integral equation-based method for solving the isotropic advection-diffusion equation in two dimensions on time-dependent geometry. First, the PDE is discretized in time by applying an IMEX SDC method, then an approximate solution in time is advanced by solving a sequence of modified Helmholtz equations. By using PUX \cite{FRYKLUNDPUX} the necessary quantities are smoothly extended to accommodate for missing data due to time-dependent geometry. The data is represented on a stationary underlying grid, while a discretization of the boundary is updated in time.  We obtain a tenth-order accurate method in the resolution of the underlying stationary uniform grid when solving the modified Helmholtz equation. To reduce the computational complexity of computing the involved sums in the boundary integral method we present a spectral Ewald method, which is of complexity $\mathcal{O}(N\log N)$ for $N$ targets and sources. The time-stepping method is adaptive, providing errors down to $10^{-8}$, and we show up to fourth order convergence in time for time-dependent deformable geometry. Moreover, the single interpolation matrix used in PUX is precomputed once and reused for all time-steps. In the precomputation step the condition number of the extension problem is significantly reduced.

This paper is an intermediate step in developing a complete solver for the simulation of drop dynamics in Stokes flow, with surfactants dissolved in the bulk fluid and on the drops' surfaces. This would mean full coupling between the advection-diffusion equation and Stokes equations. 

One of the challenges that remains is to make PUX adaptive, as high resolution will be required in the gap between two closely interacting drops. Since PUX currently requires a uniform grid for efficiency, high grid resolution is imposed globally. Another future direction is to investigate alternatives to IMEX SDC methods, in order to reduce the number of modified Helmholtz equations to solve per time-step. Future work also includes closer investigation of the CFL-type conditions that are imposed by using elliptic marching with SDC.

%%% Local Variables:
%%% mode: latex
%%% TeX-master: "manuscript.tex"
%%% End:

%% file: appsupewald.tex
\subsection{Derivation of the Ewald decomposition for $\mhhu$}
We present the derivation of the Ewald decomposition for $\mhhu$ in the free-space case, seeking $\ewrealker$ and $\widehat{\ewfourierker}$ in \eqref{eq:ewalddecomp}, and $\ewselfker\fp{\ewparam}$ in \eqref{eq:ewalddecompself}.
%The derivation for the periodic setting is analogous and can be found in \cite{PalssonSara2019SaEs}. 

Consider the sum \eqref{eq:mhhslpdiscewald}, it can be seen as the solution to
\begin{equation}
  \label{eq:mhhslpdiscpde}
 (\alpha^{2}-\Delta) \frac{1}{2\pi}\sum\limits_{n = 1}^{\Nbdry}\kerdomain\fp{\y_{n},\x}\ewfunn = \sum\limits_{n = 1}^{\Nbdry}\diracfn(\x), \quad \x\in\domain,
\end{equation}
where $\diracfn(\x)=\delta\fp{\x-\y_{n}}\ewfunn$, with $\delta$ being the Dirac delta function. To produce an Ewald decomposition of \eqref{eq:mhhslpdiscpde} a screening function is required. Note that for $\mha = 0$ the modified Helmholtz equation is reduced to the Poisson equation. Inspired by the screening function for the Poisson equation \cite{EwaldPP1921DBou} we present $\ewscreen$, see \eqref{eq:ewaldscreening}, and split $\diracfn$ accordingly
\begin{equation}
  \label{eq:ewalddecompkernel}
  \diracfn\fp{\x} = \underbrace{ \diracfn\fp{\x} - \fp{ \diracfn\fp{\cdot}\ast\ewscreen\fp{\cdot,\ewparam}}\fp{\x}}_{\diracfnr\fp{\x}}
  + \underbrace{\fp{\diracfn\fp{\cdot}\ast\ewscreen\fp{\cdot,\ewparam}}\fp{\x}}_{\diracfnf\fp{\x}}.
\end{equation}
By the linearity of the operator $\mhaa - \Delta$ the quantity $\mhhu$ can be decomposed as $\mhhu(\x) = \sum_{n}\ewur_{n}\fp{\x}+\ewuf_{n}\fp{\x}$, where
\begin{align}
  (\mhaa - \Delta)\ewur_{n}\fp{\x} &= \diracfnr\fp{\x},\\
  \label{eq:ewewufnpde}
  (\mhaa - \Delta)\ewuf_{n}\fp{\x} &= \diracfnf\fp{\x}.
\end{align}
We start out with finding a closed form expression for $\ewuf$ in the frequency domain. To this end, write
\begin{equation}
(\mhaa + \Delta) \ewuf\fp{\x} = \frac{1}{4\pi^{2}}\int\limits_{\mathbb{R}^{2}}(\mhaa + \norm{\fovarvec}^{2})\widehat{\ewuf}e^{i\fovarvec\cdot\x}\,\D\fovarvec,
\end{equation}
which is equal to
\begin{equation}
\diracfnf\fp{\x} = \frac{1}{4\pi^{2}}\int\limits_{\mathbb{R}^{2}}\widehat{\diracfnf\fp{\fovarvec}}e^{i\fovarvec\cdot\x}\,\D\fovarvec = \frac{\ewfunn}{4\pi^{2}}\int\limits_{\mathbb{R}^{2}}\widehat{\ewscreen\fp{\fovarvec}}e^{i\fovarvec\cdot(\x-\y_{n})}\,\D\fovarvec.
\end{equation}
By the orthogonality condition $(\mhaa + \norm{\fovarvec}^{2})\widehat{\ewuf} = \ewfunn \widehat{\ewscreen\fp{\fovarvec}}e^{-i\fovarvec\cdot\y_{n}}$ for all $\fovarvec\in\mathbb{R}^{2}$, thus

\begin{equation}
\ewuf_{n}\fp{\x} = \frac{\ewfunn}{4\pi^{2}}\int\limits_{\mathbb{R}^{2}}\frac{1}{\mhaa + \norm{\fovarvec}^{2}}\widehat{\ewscreen\fp{\fovarvec}}e^{i\fovarvec\cdot(\x-\y_{n})}\,\D\fovarvec.
\end{equation}
If we define $\ewuf = \sum_{n}\ewuf_{n}$, then by comparing to \eqref{eq:ewalddecomp} it is clear that

\begin{equation}
  \widehat{\ewfourierker}(\fovarvec,\ewparam) = \frac{2\pi}{\mhaa+\norm{\fovarvec}^{2}}e^{-(\mhaa + \norm{\fovarvec}^{2})/4\ewparam^{2}},\quad \fovarvec\in\mathbb{R}^{2},
\end{equation}
with $\besselK{0}\fp{\cdot,\cdot}$ defined as \eqref{eq:incompbesselK}.

Deriving a closed form for $\ewrealker$ is not as straightforward, as such has not been found by direct computation of $\kerdomain\fp{\x,\y_{n}} - \fp{\kerdomain\fp{\cdot,\y_{n}}\ast\ewscreen\fp{\cdot,\ewparam}}\fp{\x}$. Instead, it is expressed in Fourier space as
\begin{equation}
  \label{eq:ewfouriertransrealker}
  \widehat{\ewrealker}\fp{\fovarvec,\ewparam} =   \widehat{\kerdomain}\fp{\fovarvec} - \widehat{\kerdomain}\fp{\fovarvec}\widehat{\ewscreen}(\fovarvec,\ewparam) = \frac{2\pi}{\mhaa + \norm{\fovarvec}^{2}}\left( 1 - e^{-(\mhaa + \norm{\fovarvec}^{2})/4\ewparam}\right).
\end{equation}
Introduce polar coordinates so that $\fovarvec = \kappa[\cos\fp{\theta},\sin\fp{\theta}]$ and $\mathbf{r} = r[\cos\fp{\beta},\sin\fp{\beta}]$. The inverse Fourier transform of \eqref{eq:ewfouriertransrealker} can then be written as
\begin{align}
  \ewrealker(\mathbf{r},\ewparam) &= \frac{1}{4\pi^{2}}\int\limits_{\mathbb{R}^{2}}\widehat{\ewrealker}\fp{\fovarvec,\ewparam}e^{i\fovarvec\cdot\mathbf{r}}\,\D\fovarvec = \frac{1}{2\pi}\int\limits_{\kappa = 0}^{\infty}\int\limits_{\theta = 0}^{2\pi} \frac{1 - e^{-(\mhaa + \kappa^{2})/4\ewparam}}{\mhaa + \kappa^{2}}e^{i\kappa r \cos(\theta - \beta)}\kappa \,\D\theta \D\kappa\\
  &= \underbrace{\int\limits_{0}^{\infty}\frac{\kappa}{\mhaa + \kappa^{2}}\op{J}\fp{\kappa r}\,\D\kappa}_{= I_{1}} - \underbrace{\int\limits_{0}^{\infty}\frac{\kappa}{\mhaa + \kappa^{2}}\op{J}\fp{\kappa r} e^{-(\mhaa + \kappa^{2})/4\ewparam}\,\D\kappa}_{=I_{2}}.
\end{align}
The first integral $I_{1}$ is identically equal to $\besselK{0}\fp{\mha r}$. To evaluate $I_{2}$ we proceed as in \cite{TornbergAnnaKarin2016TEsf} and make the change of variables $\lambda = 1/4\ewparam^{2}$ and differentiate $I_{2}$ with respect to $\lambda$. By the dominated convergence theorem we may interchange integration and differentiation,
\begin{equation}
  \pd{I_{2}}{\lambda} = e^{-\mhaa\lambda}\int\limits_{0}^{\infty}(\mhaa + \kappa^{2})\frac{\kappa}{\mhaa + \kappa^{2}}\op{J}\fp{\kappa r} e^{-\kappa^{2}\lambda}\,\D\kappa =- e^{-\mhaa\lambda}\frac{e^{-r^{2}/4\lambda}}{2\lambda}.
\end{equation}

Before proceeding we study the limit values of $\partial I_{2}/\partial\lambda$ as $\lambda$ goes to zero and infinity, respectively. Let $\rho = 1/4\lambda$, then
\begin{equation}
\lim\limits_{\lambda\rightarrow 0_{+}} -e^{-\mhaa\lambda}\frac{e^{-r^{2}/4\lambda}}{2\lambda} =  \lim\limits_{\rho\rightarrow \infty_{+}}-2e^{-\mhaa/4\rho}e^{-r^{2}\rho}\rho = 0,
\end{equation}
and
\begin{equation}
\lim\limits_{\lambda\rightarrow \infty} -e^{-\mhaa\lambda}\frac{e^{-r^{2}/4\lambda}}{2\lambda} = 0.
\end{equation}
We now integrate $\partial I_{2}/\partial\lambda$ with respect to $\lambda$,
\begin{equation}
  \begin{aligned}
  I_{2} &= -\int\limits_{0}^{\lambda} e^{-\mhaa\rho}\frac{e^{-r^{2}/4\rho}}{2\rho}\,\D\rho = \frac{1}{2}\int\limits_{\lambda}^{\infty} e^{-\mhaa\rho}\frac{e^{-r^{2}/4\rho}}{\rho}\,\D\rho = \frac{1}{2}\int_{1}^{\infty}\frac{e^{-\mhaa t/4\ewparam^{2}}e^{-r^{2}\ewparam^{2}/t}}{t}\,\D t \\&= \frac{1}{2}\besselK{0}\left(\frac{\mhaa}{4\ewparam^{2}},r^{2}\ewparam^{2}\right),
  \end{aligned}
\end{equation}
by definition \eqref{eq:incompbesselK}, and reinserting $\lambda = 1/4\ewparam^{2}$. We now have

\begin{equation}
\ewrealker(\y-\x,\ewparam) = I_{1} - I_{2} =
\besselK{0}\fp{\mha\norm{\y-\x}}-\frac{1}{2}\besselK{0}\fp{\frac{\mhaa}{4\ewparam},\norm{\y-\x}^{2}\ewparam^{2}}=
\frac{1}{2}\besselK{0}\fp{\norm{\y-\x}^{2}\ewparam^{2},\frac{\mhaa}{4\ewparam}},
\end{equation}
by using the relation
\begin{equation}
  \label{eq:relationincompmodbess}
  \besselK{0}(\rho_{1},\rho_{2}) =   2\besselK{0}(2\sqrt{\rho_{1}\rho_{2}})-  \besselK{0}(\rho_{2},\rho_{1}) .
\end{equation}
The self-interaction term is obtained by passing the limit
\begin{equation}
\ewselfker\fp{\ewparam} = \lim\limits_{\norm{\mathbf{r}}\rightarrow 0} \frac{1}{2\pi}(\ewrealker(\mathbf{r},\ewparam) - \kerdomain(0,\mathbf{r})) = \lim\limits_{\norm{\mathbf{r}}\rightarrow 0} \frac{1}{4\pi}\besselK{0}\left(\frac{\mhaa}{4\ewparam^{2}},r^{2}\ewparam^{2}\right) = \frac{-1}{4\pi} E_{1}\left(\frac{\mhaa}{4\ewparam^{2}}\right),
\end{equation}
again using \eqref{eq:relationincompmodbess} together with $\besselK{\nu}(\rho) = E_{\nu+1}(\rho)$ \cite{HarrisFrankE2008IBgi}.
\subsection{Derivation of Ewald decomposition for $\nabla\mhhu$}
We present the derivation of the Ewald decomposition for $\mhhu$ in the free-space case, seeking $\ewrealkergrad$ and $\widehat{\ewfourierkergrad}$ in \eqref{eq:ewalddecompgrad}, and $\ewselfkergrad$. The derivation for the periodic setting is analogous and can be found in \cite{PalssonSara2019SaEs}. 

To obtain  $\ewrealkergrad$ and $\widehat{\ewfourierkergrad}$, simply take the gradient, with respect to $\x$, of \eqref{eq:ewalddecomp} and compare with the decomposition \eqref{eq:ewalddecompgrad}. We have
\begin{equation}
  \ewrealkergrad(\y-\x,\ewparam) = \ewparam^{2}(\y-\x)\besselK{-1}\left(\|\y-\x\|^{2}\ewparam^{2},\frac{\mhaa}{4\ewparam^{2}}\right),
\end{equation}
utilizing $\partial \besselK{0}(\rho_{1},\rho_{2})/\partial \rho_{1} = -\besselK{-1}(\rho_{1},\rho_{2})$  from \cite{HarrisFrankE2008IBgi}, and
\begin{equation}
\widehat{\ewfourierkergrad}(\fovarvec,\ewparam) = \frac{-i2\pi\fovarvec}{\mhaa+\norm{\fovarvec}^{2}}e^{-(\mhaa + \norm{\fovarvec}^{2})/4\ewparam^{2}},
\end{equation}
since $\nabla e^{-i \fovarvec\cdot\x} = -i\fovarvec e^{-i\fovarvec\cdot\x}$.

The self-interaction term $\ewselfkergrad$ is computed through the limit 
\begin{equation}
\ewselfkergrad = \lim\limits_{\norm{\mathbf{r}}\rightarrow 0} \frac{1}{2\pi}(\ewrealkergrad(\mathbf{r},\ewparam) - \nabla\kerdomain(0,\mathbf{r})) = \lim\limits_{\norm{\mathbf{r}}\rightarrow 0} \frac{-\ewparam^{2}\mathbf{r}}{2\pi}\besselK{1}\left( \frac{\mhaa}{4\ewparam^{2}},r^{2}\ewparam^{2}\right) = 0,
\end{equation}
using \eqref{eq:mhhkerdomaingrad} with
\begin{equation}
  \label{eq:relationincompmodbess2}
  \besselK{1}(\rho_{1},\rho_{2}) =   2\sqrt{\rho_{1}/\rho_{2}}\besselK{1}(2\sqrt{\rho_{1}\rho_{2}})-  \besselK{-1}(\rho_{2},\rho_{1}),
\end{equation}
and $\besselK{1}(\rho,0) = E_{2}(\rho)$ by \cite{HarrisFrankE2009MfiB}.

\subsection{Computing $\ewuf$ in the spectral Ewald summation for small $\mha$}
If the parameter $\mha$ in the modified Helmholtz equation \eqref{eq:mh}, \eqref{eq:mhbc} is small, then computing $\ewuf\fp{\x}$ in \eqref{eq:ewalddecomp} with the trapezoidal rule, accelerated with FFTs, will not be accurate. This is due to $\widehat{\ewfourierker}$ containing the factor $1/(\mhaa + \norm{\fovarvec})$ and therefore is nearly singular at $k=0$. To rectify this, we proceed as in \cite{afKlintebergLudvig2017FEsf}, based on the techniques in \cite{VicoFelipe2016Fcwf}, and cut off interactions in physical space beyond the domain of interest, i.e. $\suppbox$. Introduce
\begin{equation}
    \ewkertrun\fp{\y,\x} = \ewker\fp{\y,\x}\mathrm{rect}\fp{\frac{\norm{\x-\y}}{2\mathcal{R}}},
\end{equation}
where $\mathcal{R}$ is larger than the largest point-to-point distance in the domain $\suppbox$, and
\begin{equation}
  \label{eq:ewaldcutoff}
  \mathrm{rect}\fp{x} = \begin{cases}
    1,\quad x \leq 1/2,\\
    0,\quad x > 1/2.
    \end{cases}
\end{equation}
We now replace $\widehat{\ewfourierker}$ in \eqref{eq:ewalddecomp} with
\begin{equation}
  \widehat{\ewfourierkertrun}(\fovarvec,\ewparam) = \frac{2\pi}{\mhaa + \norm{\fovarvec}^{2}}(1 + \mha\norm{\fovarvec}\besselJ{1}\fp{\norm{\fovarvec}\mathcal{R}}\besselK{0}(\mha\mathcal{R})-\mha\mathcal{R}\besselJ{0}\fp{\norm{\fovarvec}\mathcal{R}}\besselK{1}(\mha\mathcal{R}))e^{-(\mhaa + \norm{\fovarvec}^{2})/4\ewparam^{2}},
\end{equation}
which has a well-defined limit at $k = 0$ as $\mha$ goes to zero. We can now write 
\begin{equation}
  \label{eq:ewaldsumfourierapprox}
\ewuf(\x) \approx \frac{1}{4\pi^{2}}\int\limits_{\mathbb{R}^{2}}\widehat{\ewfourierkertrun}(\fovarvec,\ewparam)\sum\limits_{n = 1}^{\Nbdry}\ewfunn e^{i\fovarvec\cdot(\y_{n}-\x)}\D \fovarvec.
\end{equation}
The error in the approximations due to cut offs can be controlled by $\mathcal{R}$. No errors will be introduced as we set $\mathcal{R} = \sqrt{2}L$  \cite{PalssonSara2019SaEs}. The real space part $\ewrealker$ and self-interaction part $\ewselfker$ in \eqref{eq:ewalddecomp} remain unchanged. Analogously, for the $\fovar$ space part for $\nabla\mhhu$ we have
\begin{equation}
  \widehat{\ewfourierkergradtrun}(\fovarvec,\ewparam) =  \frac{i 2\pi\fovarvec}{\mhaa + \norm{\fovarvec}^{2}}(1 + \mha\norm{\fovarvec}\besselJ{1}\fp{\norm{\fovarvec}\mathcal{R}}\besselK{0}(\mha\mathcal{R})-\mha\mathcal{R}\besselJ{0}\fp{\norm{\fovarvec}\mathcal{R}}\besselK{1}(\mha\mathcal{R}))e^{-(\mhaa + \norm{\fovarvec}^{2})/4\ewparam^{2}}.
\end{equation}
Note that when $\mhaa$ is large enough there is no need to modify $\ewfourierker$ and $\ewfourierkergrad$, see figure $6$ in \cite{PalssonSara2019SaEs}. This is the case when solving the advection-diffusion equation in this paper.
%%% Local Variables:
%%% mode: latex
%%% TeX-master: "manuscript.tex"
%%% End: